\documentclass[12pt,twoside]{article}
\usepackage{}
\usepackage[hmargin=0.8in,vmargin=0.9in]{geometry}
\usepackage{fancyhdr}
\usepackage{graphicx}
\usepackage{subfigure}
\usepackage{amssymb}
\usepackage{amsmath}
\usepackage{amsfonts}
\usepackage{theorem}
\usepackage{mathrsfs}
\usepackage{mathtools}
\usepackage{empheq}
\usepackage{bm}
\usepackage{bbm}
\usepackage{color}
\usepackage{setspace}
\usepackage{epstopdf}
\usepackage{float}
\usepackage{cite}
\usepackage{nicefrac}
\usepackage{booktabs} 
\usepackage{array} 
\usepackage{paralist} 
\usepackage{verbatim} 
\usepackage{subfigure} 
\usepackage[normalem]{ulem}

{\theorembodyfont{\rmfamily}\newtheorem{remark}{Remark}[section]}

\newenvironment{acknowledgments}{{\flushleft \bf Acknowledgment:}}{}

\newcommand\eref[1]{(\ref{#1})}

\newcommand*\xbar[1]{%
  \hbox{%
    \vbox{%
      \hrule height 0.5pt 
      \kern0.4ex
      \hbox{%
        \kern-0.05em
        \ensuremath{#1}%
        \kern-0.00em
      }%
    }%
  }%
}

\allowdisplaybreaks[1]

\newcommand{\kph}{{k+\frac{1}{2}}}
\newcommand{\kmh}{{k-\frac{1}{2}}}
\newcommand{\dx}{\Delta x}
\newcommand{\dy}{\Delta y}
\newcommand{\dt}{\Delta t}

\numberwithin{equation}{section}
\numberwithin{figure}{section}
\numberwithin{table}{section}

\title{A Well-Balanced Central-Upwind Scheme for the Thermal Rotating Shallow Water Equations}
\author{Alexander Kurganov\thanks{Department of Mathematics, Southern University of Science and Technology,
Shenzhen, 518055, China and Mathematics Department, Tulane University, New Orleans, LA 70118, USA; {\tt kurganov@math.tulane.edu}}, ~Yongle
Liu\thanks{Department of Mathematics, Harbin Institute of Technology, Harbin, 150001, China and Department of Mathematics, Southern
University of Science and Technology, Shenzhen, 518055, China; {\tt 11749318@mail.sustech.edu.cn}} ~and Vladimir Zeitlin \thanks{Laboratoire
de M\'et\'eorologie Dynamique, Sorbonne Universit\'e (SU), Ecole Normale Sup\'erieure (ENS), CNRS, Paris, 75231, France and Department of
Mathematics, Southern University of Science and Technology, Shenzhen, 518055, China; {\tt zeitlin@lmd.ens.fr}}}

\begin{document}

\date{}
\maketitle

\begin{abstract}
We develop a well-balanced central-upwind scheme for rotating shallow water model with horizontal temperature and/or density gradients---the
thermal rotating shallow water (TRSW). The scheme is designed using the flux globalization approach: first, the source terms are
incorporated into the fluxes, which results in a hyperbolic system with global fluxes; second, we apply the Riemann-problem-solver-free
central-upwind scheme to the rewritten system. We ensure that the resulting method is well-balanced by switching off the numerical diffusion
when the computed solution is near (at) thermo-geostrophic equilibria.

The designed scheme is successfully tested on a series of numerical examples. Motivated by future applications to large-scale motions in the
ocean and atmosphere, the model is considered on the tangent plane to a rotating planet both in mid-latitudes and at the Equator. The
numerical scheme is shown to be capable of quite accurately maintaining the equilibrium states in the presence of nontrivial topography and
rotation. Prior to numerical simulations, an analysis of the TRSW model based on the use of Lagrangian variables is presented, allowing one
to obtain criteria of existence and uniqueness of the equilibrium state, of the wave-breaking and shock formation, and of instability
development out of given initial conditions. The established criteria are confirmed in the conducted numerical experiments.
\end{abstract}



\section{Introduction}\label{sec1}
Rotating shallow water (RSW) model is a classical modeling tool in geophysical fluid dynamics, which helps to understand a variety of major 
dynamical phenomena in the atmosphere and oceans by simple analytical and computational means; see, e.g., \cite{Zeitlin2018}. It has,
however, a soft spot: the absence of horizontal gradients of temperature and/or density. The RSW model can be obtained from the full
``primitive'' hydrostatic equations of the geophysical fluid dynamics by vertical averaging under the hypothesis of columnar motion. As it
is quite common in geophysical modeling, we will be considering the fluid on the tangent to the rotating planet plane, and take into account
only the normal to the plane component of the angular velocity of planet's rotation (so-called ``traditional'' approximation). The
assumption of horizontal homogeneity of temperature/density accompanies the standard derivation. Yet, this assumption can be relaxed, which
does not substantially alter the derivation and leads to the so-called thermal shallow water (TSW) model \cite{Dellar}, which we will call
thermal rotating shallow water (TRSW) equations in the presence of rotation. This model was repeatedly rediscovered and used in the
literature both in the atmospheric and oceanic context, in particular for studies of the mixed layer in early papers
\cite{Rotunno,Lavoie,McCreary,Ripa95,Salby,Young}. It was recently applied to planetary atmospheres \cite{Cho,Warnerford}. Structurally,
while the classical shallow water equations are equivalent to those of isentropic gas dynamics with pressure depending only on density, the
TSW model corresponds to the dynamics of a gas with a specific equation of state, depending both on density and temperature; see
\cite{Zeitlin2018}.

One of the principal applications of the TRSW equation is to model atmospheric and oceanic temperature fronts; see \cite{Rotunno,YoungChen}.
This is why a numerical method capable of accurately resolving sharp temperature and pressure fronts is desirable. Development of such
method will be in the focus of this paper.

We consider a one-dimensional (1-D) TRSW equations with nonflat bottom topography, Coriolis force, buoyancy, and without dependence on the
zonal coordinate $x$. The studied system (see \cite{GLZD,GLZD2}) reads as
\begin{subequations}\label{1.1}
\begin{align}
&\frac{{\rm d}h}{{\rm d}t}+hv_y=0,\label{1.1a}\\
&\frac{{\rm d}u}{{\rm d}t}-fv=0,\label{1.1b}\\
&\frac{{\rm d}v}{{\rm d}t}+fu+bh_y+\frac{h}{2}b_y=-bZ_y,\label{1.1c}\\
&\frac{{\rm d}b}{{\rm d}t}=0,\label{1.1d}
\end{align}
\end{subequations}
where $y$ is a meridional coordinate, $t$ is time, $\frac{{\rm d}}{{\rm d}t}=\frac{\partial}{\partial_t}+v\frac{\partial}{\partial_y}$ is
the Lagrangian (material) derivative, $h(y,t)$ denotes the thickness of the fluid layer, $Z(y)$ represents the bottom topography, $u(y,t)$
and $v(y,t)$ stand for the zonal and meridional velocities, respectively, and $b(y,t)$ will be called buoyancy. In the oceanographic
applications, $b(y,t)=g\Delta\rho/\rho_0$, where $\Delta\rho:=\rho_0-\rho(y,t)$ is the density difference between the density $\rho(y,t)$
and a reference density $\rho_0$. Notice that density variations in the ocean are proportional to temperature variations, whence the name of
the model. In the atmospheric applications, $\rho$ and $\rho_0$ should be replaced with the potential temperature $\theta$ and $\theta_0$,
while $\Delta\theta:=\theta-\theta_0$. In both cases, $g$ is the acceleration due to gravity. Finally, rotation enters via $f(y)$, the
Coriolis parameter. In the simplest $f$-plane approximation, where the effects of curvature are neglected, $f$ is assumed to be constant,
that is, $f(y)\equiv f_0$. If the curvature is taken into account to the first order in the beta-plane approximation, $f$ becomes
$f(y)=f_0+\beta y,~\beta={\rm Const}$, which reduces to $f(y)=\beta y$ near the Equator, where the tangent plane is parallel to the axis of
rotation.

As one can easily show, the system \eref{1.1} can be rewritten in the form of the following system of balance laws:
\begin{subequations}\label{1.2}
\begin{align}
&h_t+\big(hv\big)_y=0,\label{1.2a}\\
&q_t+\big(huv\big)_y=fp,\label{1.2b}\\
&p_t+\Big(hv^2+\frac{b}{2}h^2\Big)_y=-fq-hbZ_y,\label{1.2c}\\
&\big(hb\big)_t+\big(hvb\big)_y=0,\label{1.2d}
\end{align}
\end{subequations}
where $q:=hu$ and $p:=hv$ are the zonal and meridional discharges or, in other words, the densities of zonal and meridional momenta,
respectively.

We notice that if the Coriolis force is not considered ($f\equiv0$) and $\rho\equiv0$ so that $b(y,t)\equiv g$, then the system \eref{1.2}
is reduced to the classical Saint-Venant system of shallow water equations, which was originally derived in \cite{SV} and is still widely
used in modeling water flows in rivers, canals, lakes, reservoirs and coastal areas. The 1-D Saint-Venant system reads as:
\begin{align}
\begin{split}
&h_t+p_y=0,\\
&p_t+\big(hv^2+\frac{1}{2}gh^2\big)_y=-ghZ_y.
\end{split}
\label{1.3}
\end{align}

It is well-known that steady-state solutions are of great particular importance because they minimize, if they are stable, the energy, and
many practically relevant waves can be viewed as small perturbations of those equilibria. The system, whenever there is an energy sink of
any nature, tends to reach a steady state. This process is called {\em adjustment}. The steady-state solutions of the Saint-Venant system
\eref{1.3} satisfy the time-independent system:
\begin{align}
\begin{split}
&\big(hv\big)_y=0,\\
&\big(hv^2+\frac{1}{2}gh^2\big)_y=-ghZ_y.
\end{split}
\label{1.4}
\end{align}
It is easy to show that the system \eref{1.4} admits a family of smooth steady-state solutions:
\begin{equation}
p\equiv\mbox{Const},\quad\frac{v^2}{2}+(h+Z)\equiv\mbox{Const}.
\label{1.5}
\end{equation}
Among which, one of the most practical relevant steady states is the so-called ``lake at rest'' equilibria, which is obtained from
\eref{1.5} when $v\equiv0$:
\begin{equation*}
v\equiv0,\quad w:=h+Z\equiv \mbox{Const}.
\end{equation*}

Capturing these steady states or their small perturbations (quasi-steady flows) accurately is a highly nontrivial task, because it may lead
to spurious oscillations when the shock capturing numerical methods are directly applied. Although such unphysical oscillations may be
attenuated once a very fine computational grid is used, this may not be affordable in practice where coarse grids have to be used.
Therefore, in the past decades, there have been many attempts to develop well-balanced numerical methods that are capable of exactly
preserving the aforementioned steady-states solutions at the discrete level; see, e.g.,
\cite{ABBKP,BCKN,BNL,CLPwb13,CCHKW_18,CK16,Kur_Acta,KLshw,KPshw,LeV98,XS14} and references therein.

Similar to the Saint-Venant system \eref{1.3}, there are some special steady states of the 1-D TRSW system \eref{1.2}. However, due to the
presence of Coriolis force and density/temperature variations, the structure of the steady state solutions becomes more complex. In fact,
the steady-state solutions of \eref{1.2} can be obtained by solving the system:
\begin{equation*}
\left\{\begin{aligned}
&p_y=0,\\
&\Big(\frac{pq}{h}\Big)_y=fp,\\
&\Big(\frac{p^2}{h}+\frac{b}{2}h^2\Big)_y=-fq-hbZ_y,\\
&\big(pb\big)_y=0.
\end{aligned}\right.
\end{equation*}
As one can easily see, in the absence of Coriolis force $(f=0)$, the system admits several particular steady-state solutions, two of which
are the following ``lake at rest'' ones:
\begin{equation}
p\equiv0,\quad b\equiv\mbox{Const},\quad w:=h+Z\equiv\mbox{Const},
\label{1.8}
\end{equation}
and
\begin{equation}
p\equiv0,\quad Z\equiv\mbox{Const},\quad \frac{b}{2}h^2\equiv\mbox{Const}.
\label{1.9}
\end{equation}
A well-balanced central-upwind scheme supplemented with interface tracking method for preserving the steady states \eref{1.8} and \eref{1.9}
was introduced in \cite{CKL}. A path-conservative approximate Riemann-problem-solver based scheme for the TSW system \eref{1.2} with $f=0$
was developed in \cite{SMC}. This scheme is well-balanced, positivity preserving and entropy dissipative in the case of flat or continuous
bottom topography $Z$. For several other recently developed well-balanced schemes for the TSW model, we refer the reader to
\cite{DZBK_Ripa,HL,Her-Due,TouKli}.

If both the Coriolis force and buoyancy are taken into account, the situation is even more complicated and it is quite challenging to design
an accurate and robust well-balanced scheme for the studied system. In this paper, we follow the idea from \cite{CCHKW_18,CCKOT_18,CHO_18},
and incorporate the source term in \eref{1.2c} into the corresponding flux term and rewrite \eref{1.2} in the following equivalent form:
\begin{equation}
\left\{\begin{aligned}
&h_t+p_y=0,\\
&q_t+\Big(\frac{pq}{h}\Big)_y=fp,\\
&p_t+L_y=0,\\
&\big(hb\big)_t+\big(pb\big)_y=0,
\end{aligned}\right.
\label{1.10}
\end{equation}
where
\begin{equation}
L:=\frac{p^2}{h}+\frac{b}{2}h^2+R,
\label{L}
\end{equation}
is a global equilibrium variable and
\begin{equation}
R(y,t):=\int\limits^y\Big[f(\xi)q(\xi,t)+h(\xi,t)b(\xi,t)Z_y(\xi)\Big]{\rm d}\xi.
\label{R}
\end{equation}

It is easy to see that the system \eref{1.10} admits the equilibria, which will be called {\em thermo-geostrophic} \cite{GLZD2}, as they are
due both to effects of rotation and temperature variation, and can be written using the global variables in an extremely simple form:
\begin{equation}
p\equiv0,\quad L\equiv{\rm Const}.
\label{1.13}
\end{equation}
On the other hand, \eref{1.10} is a hyperbolic system with a global flux and therefore it might be quite challenging to design a
Riemann-problem-solver-based upwind scheme for solving it numerically.

In this paper, our main goal is to develop a robust, stable and highly accurate finite-volume method for the TRSW system \eref{1.10}. We
note that in this system, the flux is global as it depends on a global variable $L$. This makes it hard to design upwind methods for
\eref{1.10}. We therefore derive a Riemann-problem-solver-free second-order central-upwind scheme. Central-upwind schemes have been
originally developed for hyperbolic systems of conservation laws \cite{KNP1,KLin,KTcl} and then extended and applied to hyperbolic systems
of balance laws arising in modeling shallow water flows; see, e.g., \cite{CKM,CKL,Kur_Acta,KLshw,KM,KPshw}. Here, we proceed along the lines
of previous works in \cite{CCHKW_18,CCKOT_18,CDKL,CHO_18} and develop a central-upwind scheme, which is well-balanced in the sense that it
exactly preserves the geostrophic equilibria \eref{1.13}. This is achieved by performing a piecewise linear reconstruction of the
equilibrium variables followed by well-balanced evolution, which uses modified central-upwind fluxes similar to those presented in
\cite{CCKOT_18}. Preserving the positivity of $h$ and $b$ is another crucial property a good scheme must possess. This is enforced with the
help of the so-called ``draining time step'' technique originally introduced in \cite{BNL}.

The paper is organized as follows. In \S\ref{sec2}, we present a general analysis of the TRSW system, first in the $f$-plane, then in the
equatorial beta-plane approximations. We obtain conditions of existence and uniqueness of the thermo-geostrophic equilibrium corresponding
to arbitrary initial conditions, study the process of relaxation (adjustment) to such equilibrium, and get conditions of wave-breaking and
shock formation. A well-balanced semi-discrete second-order central-upwind scheme for the TRSW system \eref{1.10} is presented in
\S\ref{sec3}. The proposed scheme is tested on a number of numerical examples reported in \S\ref{sec4}.

\section{Analysis of the One-Dimensional TRSW Equations}\label{sec2}
In this section, we analyze the studied TRSW equations with the flat bottom topography, that is, with $Z_y\equiv0$. We first obtain the
results in the case of constant Coriolis parameter $f(y)\equiv{\rm Const}$ (\S\ref{sec23}--\S\ref{sec25}). To this end, we first rewrite the
TRSW equations in the Lagrangian form (\S\ref{sec21}) and then analyze the existence of thermo-geostrophic equilibria for given initial
conditions (\S\ref{sec23}), study the adjustment process (\S\ref{sec24}), and demonstrate the breakdown of smooth solutions (\S\ref{sec25}).
We then extend these results to the equatorial case with $f(y)=\beta y$, $\beta={\rm Const}$ (\S\ref{sec26}), where the situation proves to
be different in several aspects.

\subsection{Lagrangian Formulation of the TRSW Equations}\label{sec21}
We first introduce Lagrangian positions of fluid parcels on the $y$-axis, $Y=Y(y,t)$, with, correspondingly,
$v(Y,y)=\frac{{\rm d}Y}{{\rm d}t}$ and $Y(y,0)=y$. We then use the Euler-Lagrange duality in the description of the fluid motions, that is,
the fact that at any time moment at each point in space we have a Lagrangian parcel, in order to rewrite equations \eref{1.1b} and
\eref{1.1c} with $Z_y\equiv0$ as follows:
\begin{align}
&\dot{u}-f\dot{Y}=0,\label{13}\\
&\ddot{Y}+fu+bh_Y+\frac{h}{2}b_Y=0,\label{14}
\end{align}
where we have used a dot notation $\dot{(\cdot)}:=\frac{{\rm d}}{{\rm d}t}(\cdot)$ for the Lagrangian time derivative. We notice that
equation \eref{13} expresses the Lagrangian conservation of the {\em geostrophic momentum} $u-fY$. The mass conservation equation
\eref{1.1a} in the Lagrangian formulation becomes
\begin{equation}
h(Y)\,{\rm d}Y=h_0(y)\,{\rm d}y\quad\iff\quad h(Y)Y'=h_0(y),
\label{15}
\end{equation}
where $h_0(y)$ is the initial thickness distribution and a prime notation $(\cdot)':=\frac{{\rm d}}{{\rm d}y}(\cdot)$ is used for the
partial derivative with respect to $y$.

We note that equations \eref{13} and \eref{1.1d} can be immediately integrated. This results in
\begin{equation}
u(Y,t)=f(Y(y,t)-y)+u_0(y)\quad\mbox{and}\quad b(Y,t)=b_0(y),
\label{16}
\end{equation}
where $u_0(y)$ and $b_0(y)$ are the initial distributions of zonal velocity and buoyancy, respectively. For the sake of physical
consistency, $b_0$ is assumed to be positive, that is, $b_0(y)>0$ for all $y$. Using \eref{15} and \eref{16} equation \eref{14} can be
rewritten in the following form:
\begin{equation}
\ddot{Y}+f^2(Y-y)+\frac{b_0}{Y'}\left(\frac{h_0}{Y'}\right)'+\frac{h_0b_0'}{2(Y')^2}=-fu_0.
\label{17}
\end{equation}
It is often convenient to introduce the deviation of the fluid parcels from their initial positions: $\phi(y,t):=Y(y,t)-y$. Equation
\eref{17} then takes the following form:
\begin{equation}
\ddot{\phi}+f^2\phi+\frac{b_0}{1+\phi'}\left(\frac{h_0}{1+\phi'}\right)'+\frac{h_0b_0'}{2(1+\phi')^2}=-fu_0.
\label{18}
\end{equation}

\subsection{Thermo-Geostrophic Equilibrium and Adjustment}\label{sec22}
It immediately follows from \eref{18} that its stationary solution with no displacement of fluid parcels, that is, with $\phi\equiv0$ exists
if the initial data are in thermo-geostrophic equilibrium, that is, if
\begin{equation}
b_0(y)h_0(y)'+\frac{h_0(y)}{2}b_0(y)'=-fu_0(y)
\label{22}
\end{equation}
is satisfied for all $y$. If $h_0$, $u_0$ and $b_0$ do not obey \eref{22}, that is, if there is an initial imbalance
\begin{equation}
-b_0(y)h_0(y)'-\frac{h_0(y)}{2}b_0(y)'-fu_0(y)=A_0(y)\not\equiv0,
\label{23}
\end{equation}
then the solution of the TRSW system evolves. One can show that a stationary solution minimizes the total energy of the system, which is the
sum of the kinetic and potential energies:
$$
E=\int\limits_{-\infty}^\infty\left[h\left(\frac{u^2+v^2}{2}\right)+\frac{bh^2}{2}\right]{\rm d}y.
$$
As mentioned in \S\ref{sec1}, the system out of equilibrium tends to reach an equilibrium state. As the dissipation is absent in the studied
TRSW model, this can be only achieved by evacuating energy through the emission of waves. This process will be called
{\em thermo-geostrophic adjustment}. In mathematical terms, the adjustment process describes a solution of the Cauchy problem for \eref{18}
with {\em radiation boundary conditions} at spatial infinity.

\subsection{Existence of the Adjusted State}\label{sec23}
The adjusted stationary state is described by the stationary part of equation \eref{17}, which can be rewritten in the following form:
\begin{equation}
f^2(Y-y)+fu_0+\sqrt{b_0}\left[\left(\sqrt{b_0}h_0\right)'\frac{1}{(Y')^2}+\sqrt{b_0}h_0\left(\frac{1}{2(Y')^2}\right)'\,\right]=0.
\label{33}
\end{equation}
Differentiating \eref{33} with respect to $y$ results in
\begin{equation}
f^2(Y'-1)+fu_0'+\left(\frac{1}{h_0}P^{\,\prime}\right)'=0,
\label{34}
\end{equation}
where
\begin{equation}
P:=\frac{b_0h_0^2}{2(Y')^2}
\label{35}
\end{equation}
is a generalization of the {\em Lagrangian pressure} introduced in \cite{ZMP} in the context of standard RSW equations for the TRSW case.

In order to simplify equation \eref{34}, we introduce a ``straightening'' change of the space variable $y\to\eta$ such that
$h_0(\eta)=H_0={\rm Const}$. In terms of the new variable, the mass conservation equation \eref{15} reads as
\begin{equation}
h(Y)=H_0/Y',
\label{2.12}
\end{equation}
and equation \eref{34} can be rewritten in the following form:
\begin{equation}
-\frac{{\rm d}}{{\rm d}Y}\left(\sqrt{b_0}\,\frac{{\rm d}\left(\sqrt{b_0}h_0\right)}{{\rm d}Y}\right)-Q(Y)h(Y)=f,
\label{36}
\end{equation}
where $Y=Y(\eta,t)$ and
$$
Q:=\frac{f-u'(Y)}{h(Y)}
$$
is the {\em potential vorticity}. We would like to emphasize that the potential vorticity is a Lagrangian invariant in the 1-D TRSW
equations, which is not the case in the case of two spatial dimensions; see \cite{Dellar}. 

We note that equation \eref{36} differs from the corresponding equation in the standard RSW studied in \cite{ZMP} by the presence of the
$\sqrt{b_0}$ factors.

Finally, we make two more changes of variables, $h\to\widehat h=\sqrt{b_0}h$ and $Y\to\xi=\int\sqrt{b_0}\,{\rm d}Y$, and rewrite equation
\eref{36} in the canonical form:
$$
-\frac{1}{f}\frac{{\rm d}^2\widehat h}{{\rm d}\xi^2}+Q(\xi)\widehat h=f\sqrt{b_0}.
$$
Then, in the physically relevant case of monotone $b_0$ with constant asymptotics at infinities, which corresponds to a density/temperature
front, the main configuration we are aiming at, existence and uniqueness of the decaying at infinities solution can be proved for
nonnegative potential vorticity $Q\ge0$ along the lines of \cite{ZMP}.
\begin{remark}
The case of non-monotone $b_0$ requires a special consideration.
\end{remark}

\subsection{Linear Theory of Thermo-Geostrophic Adjustment}\label{sec24}
We now analyze the adjustment process in the linear approximation, assuming that the norm of the deviations of fluid parcels is small. The
linearization of equation \eref{18} with respect to $\phi$ results in
\begin{equation}
\ddot{\phi}+f^2\phi-2\left(b_0h_0'+\frac{h_0b_0'}{2}\right)\phi'-b_0h_0\phi''=A_0.
\label{24}
\end{equation}
Let us split the dynamical variable $\phi$ into the ``slow'' and ``fast'' components: $\phi=\bar\phi+\tilde\phi$, where the slow
component is defined as the time-average: $\bar\phi=\lim_{T\to\infty}\frac{1}{T}\int_0^T\phi\,{\rm d}t$. By averaging equation \eref{24} in
time, and assuming that both $\phi$ and $\dot{\phi}$ remain bounded at all times, we obtain the following ODE for $\bar\phi$:
\begin{equation}
f^2\bar\phi-2\left(b_0h_0'+\frac{h_0b_0'}{2}\right)\bar\phi'-b_0h_0\bar\phi''=A_0.
\label{25}
\end{equation}
We then introduce a new variable $\bar\Phi=\sqrt{b_0}h_0\bar\phi$ and rewrite equation \eref{25} in the canonical form:
\begin{equation}
-\bar\Phi''+\left[\frac{\frac{f^2}{\sqrt{b_0}}+\left(\sqrt{b_0}h_0\right)''}{\sqrt{b_0}h_0}\right]\bar\Phi=\frac{A_0}{\sqrt{b_0}h_0}.
\label{26}
\end{equation}
Notice that the quantity $\frac{f^2/\sqrt{b_0}+\left(\sqrt{b_0}h_0\right)''}{\sqrt{b_0}h_0}$ reduces to the geostrophic potential vorticity
$\frac{f^2/g+h_0''}{h_0}$ in the corresponding equation obtained for the standard RSW model with a constant buoyancy $b_0=g$ studied in
\cite{ZMP}. In the physically interesting frontal case when $A_0$ is a compactly supported function and $b_0$ and $h_0$ have constant
asymptotics at infinities, it was shown in \cite{ZMP} that a solution decaying at infinities exists, and it is unique for nonnegative
geostrophic potential vorticity. Decay conditions are imposed by an obvious reason that fluid parcels should not move far away from the
front, where they are already at equilibrium. Correspondingly, decaying solutions of \eref{26} exist for
$f^2/\sqrt{b_0}+\left(\sqrt{b_0}h_0\right)''\ge0$ (also see the nonlinear existence result obtained in \S\ref{sec23}).

The PDE for the ``fast'' component of $\phi$ is obtained by subtracting \eref{24} from \eref{25}, and reads as
\begin{equation}
\ddot{\tilde\phi}+f^2\tilde\phi-2\left(b_0h_0'+\frac{h_0b_0'}{2}\right)\tilde\phi'-b_0h_0\tilde\phi''=0.
\label{27}
\end{equation}
If $b_0$ and $h_0$ are both constant, this is the 1-D Klein-Gordon equation describing gravity waves propagating at the surface of the
shallow water layer. For nonconstant $b_0$ and $h_0$, equation \eref{27} describes the wave propagation over a variable background. For
harmonic waves with $\tilde\phi=\int(\psi e^{-i\omega t}+{\rm c.c.})\,{\rm d}\omega$, where c.c. stands for complex conjugation, equation
\eref{27} becomes
\begin{equation}
(\omega^2-f^2)\psi+2\left(b_0h_0'+\frac{h_0b_0'}{2}\right)\psi'+b_0h_0\psi''=0.
\label{28}
\end{equation}
In the aforementioned frontal case, where $b_0(y)\to b_\pm={\rm Const}$ and $h_0(y)\to h_\pm={\rm Const}$ as $y\to\pm\infty$, the far
asymptotics of the solutions of \eref{28} satisfy the following constant-coefficient ODEs:
\begin{equation}
(\omega^2-f^2)\psi_\pm +b_\pm h_\pm\psi''_\pm=0,
\label{28a}
\end{equation}
and correspond to propagating inertia-gravity waves with the standard dispersion relations, obtained after applying the Fourier transform
in space $\psi=\int(\hat{\psi} e^{i k x}+{\rm c.c.})\,{\rm d}k$:
\begin{equation}
\omega_\pm^2=f^2+b_\pm h_\pm k^2,
\label{29}
\end{equation}
where $k$ is the wavenumber. Hence, the adjustment process corresponding to a solution of the Cauchy problem for \eref{24} with zero initial
$\phi$ and $\dot{\phi}$ consists of the emission of inertia-gravity waves out of the initial unbalanced front, which evacuate an excess of
energy and drive the system to an equilibrium state given by a solution of \eref{26}. The question, however, arises, whether some part of
the wave signal could be trapped at the front and, thus, prevent the front from complete equilibration. In order to answer this question, we
rewrite equation \eref{28} as follows:
\begin{equation}
\sqrt{b_0}h_0\psi''+2\left(\sqrt{b_0}h_0\right)'\psi'+\frac{\omega^2-f^2}{\sqrt{b_0}} \psi=0,
\label{30}
\end{equation}
then multiply \eref{30} by $\sqrt{b_0}h_0\psi^*$, where the star denotes the complex conjugation, and rewrite it in the following form:
\begin{equation*}
\left[\left(\sqrt{b_0}h_0\right)^2\psi^*\psi' \right]'-\left(\sqrt{b_0}h_0\right)^2{\psi^*}'\psi'+(\omega^2-f^2)h_0\psi^*\psi=0.
\end{equation*}
Integrating this equation from $-\infty$ to $\infty$ in $y$, and assuming the decay of $\psi$ far from the front, as we are looking for
trapped modes, leads to the following estimate for the eigenfrequencies:
\begin{equation*}
\omega^2=f^2+\frac{\int_{-\infty}^\infty b_0h_0^2|\psi'|^2\,{\rm d}y}{\int_{-\infty}^\infty h_0|\psi|^2\,{\rm d}y},
\end{equation*}
which shows that the frequencies are {\em suprainertial}, that is, $\omega^2\ge f^2$, while in order to have trapped modes they should be
sub-inertial, which can be shown following the lines of \cite{ZMP}. It should be kept in mind, however, that the group velocity
$c_g=\partial\omega/\partial k$ of near-inertial waves with frequencies close to $f$ is small, as follows from \eref{29}. This practically
means that the portion of the ``fast'' component of the perturbation with $\omega\approx f$ will remain at the initial location for a very
long time. We will see a manifestation of this fact in the results of numerical simulations reported below.

\subsection{Breakdown of Smooth Solutions}\label{sec25}
The standard Saint-Venant system of shallow water equations is equivalent to the isentropic Euler equations of gas dynamics. Hence, shock
formation is ubiquitous in this system. As shown in \cite{ZMP}, including rotation into the Saint-Venant system does not prevent shock
formation, but changes the breakdown conditions. In this section, we study the shock formation and properties of shocks in the TRSW
equations, which are pertinent in the context of numerical simulations using shock-capturing finite-volume methods. In order to address
these questions, we follow \cite{ZMP} and use the Lagrangian description in the case of constant initial $h_0(y)=H_0={\rm Const}$, which can
be always achieved by a ``straightening'' change of the independent spatial variable; see \S\ref{sec23}.

We rewrite the Lagrangian equations of motion in the following form:
\begin{equation}
\left\{\begin{aligned}
&\dot{v}+\frac{1}{H_0}P'=-fu,\\
&\dot{J}-v'=0,
\end{aligned}\right.
\label{38}
\end{equation}
where $P$ is given by \eref{35} with the constant $h_0=H_0$. In \eref{38}, $v$ and $J:=Y'$ are dependent variables, while $u$ is not and
needs to be determined from the conservation of the geostrophic momentum and potential vorticity $Q$:
\begin{equation*}
u-fY=u_0-fy\quad\Longrightarrow\quad u'=fJ+u_0'-f=fJ-H_0Q.
\end{equation*}
It is easy to check that the system \eref{38} can be rewritten in an explicit quasi-linear form:
\begin{equation*}
\begin{pmatrix}\dot{v}\\\dot{J}\end{pmatrix}+A\begin{pmatrix}v\\J\end{pmatrix}_y=\begin{pmatrix}-fu-\frac{H_0b_0'}{2J^2}\\0\end{pmatrix},
\quad A=\begin{pmatrix}0&-H_0b_0J^{-3}\\-1&0\end{pmatrix};
\end{equation*}
compare with \cite[equation (57)]{ZMP}.

We now take (without loss of generality) $H_0=f=1$ and proceed along the lines of \cite{ZMP}. The eigenvalues of the matrix $A$ are
$\mu_\pm=\pm\sqrt{b_0}J^{-3/2}$ and the corresponding left eigenvectors are $(1,\mp\sqrt{b_0}J^{-3/2})^\top$. Hence, the Riemann invariants
are $r_\pm=v\pm2\sqrt{b_0}J^{-1/2}$ and we have
\begin{equation}
\dot{r}_\pm+\mu_\pm(r_\pm)_y=-u+\frac{b_0'}{2J^2}.
\label{2.23}
\end{equation}
Next, we differentiate \eref{2.23} with respect to $y$ and obtain that the derivatives of the Riemann invariants, $D_\pm:=(r_\pm)_y$,
satisfy the following PDEs:
\begin{equation}
\dot{D}_\pm+\mu_\pm(D_\pm)_y+(\mu_\pm)_yD_\pm=-u_y+\left(\frac{b_0'}{2J^2}\right)_y.
\label{2.24}
\end{equation}
Using the expressions for the Riemann invariants, we obtain $r_+-r_-=4\sqrt{b_0}J^{-1/2}$, which implies that
$$
\mu_\pm=\pm\frac{1}{b_0}\left(\frac{r_+-r_-}{4}\right)^3,
$$
and thus
\begin{equation}
(\mu_\pm)_y=\frac{\partial\mu_\pm}{\partial r_+}(r_+)_y+\frac{\partial\mu_\pm}{\partial r_-}(r_-)_y+\frac{\partial\mu_\pm}{\partial b_0}b_0'
=\frac{\partial\mu_\pm}{\partial r_+}D_++\frac{\partial\mu_\pm}{\partial r_-}D_-\mp\frac{b_0'}{b_0^2}\left(\frac{r_+-r_-}{4}\right)^3.
\label{2.25}
\end{equation}
Finally, substituting \eref{2.25} into \eref{2.24} and using \eref{2.12} result in
\begin{equation}
\dot{D}_\pm+\mu_\pm(D_\pm)_y\mp\frac{b_0'}{b_0^2}\left(\frac{r_+-r_-}{4}\right)^3D_\pm
+\frac{\partial\mu_\pm}{\partial r_+}D_+D_\pm+\frac{\partial\mu_\pm}{\partial r_-}D_-D\pm=-u_y+\left(\frac{b_0'h^2}{2}\right)_y.
\label{2.26}
\end{equation}
This is a generalized Ricatti equation, which can be analyzed following \cite{Engelberg}, as it was done in \cite{ZMP} for the 1-D RSW
model. Breakdown and shock formation correspond to $D_\pm$ reaching infinite values in finite time. Compared to the corresponding equations
in \cite{ZMP}, after obvious changes $v\to-u$ and $x\to y$ in the latter, we observe the following two differences.

First, the first term on the right-hand side of \eref{2.26}, corresponding to vorticity, acquires an addition
$\left(\frac{b_0'h^2}{2}\right)_y$. Recall that if the initial vorticity is sufficiently negative, breakdown always takes place in the RSW
equations. Here, in the TRSW model, it is the vorticity plus this new term (which depends on the initial distributions of buoyancy and
thickness), should be sufficiently negative for the breakdown to take place.

Second, the breakdown is conditioned by the signs of the derivatives of Riemann invariants. As follows from their definitions, the overall
sign of the derivatives depends not only on the signs of derivatives of $v$ and $h$ (as in the standard RSW model), but also on the sign of
the derivative of $b_0$, which makes a difference. It is worth emphasizing, however, that in practice it is difficult to discriminate the
role of each contribution in the simulations, and specially designed initial conditions are required in order to do this; see \cite{BSZ}.

\subsection{Extension to the Equatorial Case}\label{sec26}
In this section, we discuss a possible extension of the results obtained in \S\ref{sec23}--\S\ref{sec25} to the equatorial case with a
variable Coriolis parameter $f(y)=\beta y$, $\beta={\rm Const}$. We first reformulate the 1-D TRSW equations on the equatorial beta-plane in
Lagrangian variables. Equations \eref{13} and \eref{14} take the following form:
\begin{align}
&\dot{u}-\beta Y\dot{Y}=0,\label{13bis}\\
&\ddot{Y}+\beta Yu+b h_Y+\frac{h}{2}b_Y=0,\label{14bis}
\end{align}
while equation \eref{15} does not change. As in the $f$-plane case, equation \eref{13bis} can be easily integrated:
\begin{equation*}
u(Y,t)=\frac{\beta}{2}(Y^2(y,t)-y^2)+u_0(y)\quad\mbox{and}\quad b(Y,t)=b_0(y),
\end{equation*}
and then by expressing $h$, $b$ and $u$ in terms of their initial values, equation \eref{14bis} reduces to 
\begin{equation}
\ddot{Y}+\beta Y\left[u_0+\frac{\beta}{2}\big(Y^2(y,t)-y^2\big)\right]+\frac{b_0}{Y'}\left(\frac{h_0}{Y'}\right)'+\frac{h_0b_0'}{2(Y')^2}=0.
\label{17bis}
\end{equation}
Already the inspection of \eref{14bis} and \eref{17bis} shows the fundamental difference from the $f$-plane case: the dependence on
meridional coordinate in the Coriolis parameter introduces higher powers of $Y$. This renders impossible the procedure used above in the
$f$-plane case in the demonstration of both existence and uniqueness of the adjusted state and breakdown. This procedure consisted of
transforming, by differentiation, the original system of Lagrangian equations into a pair of equations for the $t$- and $y$-derivatives of
$Y$ (or a single equation for the space derivative in the stationary case), while eliminating the $Y$ itself. Obviously, this is not
possible on the beta-plane, which introduces serious technical difficulties (and even principal ones, like finding Riemann invariants for a
system of three quasilinear equations). We therefore will not pursue these demonstrations below and will limit ourselves only by the linear
analysis of the thermo-geostrophic adjustment on the equatorial beta-plane.

Introducing, as before, the deviations $\phi(y,t)$ of the fluid parcels from their initial positions, we rewrite \eref{17bis} as
\begin{equation}
\ddot{\phi}+\beta(y+\phi)\left[u_0+\frac{\beta}{2}\phi(2y+\phi)\right]+\frac{b_0}{1+\phi'}\left(\frac{h_0}{1+\phi'}\right)'+
\frac{h_0b_0'}{2(1+\phi')^2}=-fu_0.
\label{18bis}
\end{equation}
The thermo-geostrophic balance of the initial conditions \eref{22} with $f=\beta y$ provides a trivial solution $\phi\equiv0$. If the
imbalance \eref{23} is small, we linearize \eref{18bis} as in \S\ref{sec24}:
\begin{equation*}
\ddot{\phi}+\beta\phi(u_0+\beta y^2)-2\left(b_0h_0'+\frac{h_0b_0'}{2}\right)\phi'-b_0h_0\phi''=A_0,
\end{equation*}
split $\phi$ into slow and fast variables, and obtain the equation for the slow motion,
\begin{equation*}
\beta\bar\phi(u_0+\beta y^2)-2\left(b_0h_0'+\frac{h_0b_0'}{2}\right)\bar\phi'-b_0h_0\bar\phi''=A_0,
\end{equation*}
which is rewritten in terms of $\bar\Phi=\sqrt{b_0}h_0\bar\phi$ in the following form:
\begin{equation*}
-\bar\Phi''+\left[\frac{\frac{\beta(u_0+\beta y^2)}{\sqrt{b_0}}+\left(\sqrt{b_0}h_0\right)''}{\sqrt{b_0}h_0}\right]\bar\Phi=
\frac{A_0}{\sqrt{b_0}h_0}.
\end{equation*}
As in \S\ref{sec24}, a solution decaying at $y\to\pm\infty$ exists if the quantity in the square brackets is nonnegative. Notice, however,
that the novelty of this expression with respect to its $f$-plane counterpart is that it has an extra possibility to become negative if the
initial flow is oriented westward (negative $u_0$) and sufficiently strong, which makes the factor $u_0+\beta y^2$ negative. We will see
that the sign of this factor, which is related to the absolute zonal momentum density on the beta-plane, also plays an important role in the
dynamics of fast motions.

Again, as in \S\ref{sec24}, we obtain for the fast component
\begin{equation*}
\ddot{\tilde\phi}+\beta\tilde\phi(u_0+\beta y^2)-2\left(b_0h_0'+\frac{h_0b_0'}{2}\right)\tilde\phi'-b_0h_0\tilde\phi''=0,
\end{equation*}
and then performing the Fourier transform in time results in
\begin{equation}
\left(\omega^2-\beta(u_0+\beta y^2)\right)\psi+2\left(b_0h_0'+\frac{h_0b_0'}{2}\right)\psi'+b_0h_0\psi''=0.
\label{28bis}
\end{equation}
In the front/jet configurations, which are of primary interest, where $b_0$ and $h_0$ have constant asymptotics and $u_0$ tends to zero at
$\pm\infty$, equation \eref{28bis} becomes at the far left and far right sides of the front, respectively:
\begin{equation}
(\omega^2-\beta^2y^2)\psi_\pm+b_\pm h_\pm\psi^{''}_\pm=0.
\label{28abis}
\end{equation}
The crucial difference between these equations  and \eref{28a} is appearance of $y^2$ as a coefficient. After rescaling $y$ with the
corresponding equatorial deformation radii $R_{e_{\pm}}:=\sqrt{\sqrt{b_\pm h_\pm}/\beta}$, and $t$ by the corresponding equatorial inertial
periods $T_{e_{\pm}}=\frac{2\pi}{\beta R_{e_\pm}}$, \eref{28abis} takes a canonical form
\begin{equation}
\psi^{''}_\pm-y^2\psi^{''}_\pm=-\bar{\omega}^2_\pm\psi^{''}_\pm,
\label{2.33}
\end{equation}
where $\bar{\omega}_\pm:=\frac{\omega}{\beta R_{e_{\pm}}}$. If the natural for the equatorial region decay boundary conditions are imposed, 
equation \eref{2.33} can be solved in terms of Gauss-Hermite (parabolic cylinder) functions obeying the equation
\begin{equation*}
\psi^{''}_n-y^2\psi^{''}_n=-(2n+1)\psi^{''}_n,\quad n = 0,1,2\ldots,
\end{equation*}
where $n$ is a number of zeroes of $\psi_n$ in $y$. Hence, even in the absence of initial inhomogeneities in the buoyancy, thickness and
zonal velocity, the fast component does not represent freely propagating inertia-gravity waves, but the waves are trapped at the Equator.
The  resulting eigenfrequencies $\bar\omega=\sqrt{2n+1}$ correspond to the infinite zonal wavelength limit of the classical spectrum of the
equatorial waves; see, e.g., \cite{Zeitlin2018}. An important conclusion, following from this analysis, is that the fast component cannot be
evacuated, like in the $f$-plane case, but remains trapped at the Equator.

Another peculiarity of the equatorial adjustment is a possible appearance  of additional trapped modes and even of instability for strong
enough westward jets. Indeed, along the same lines as in \S\ref{sec24}, equation \eref{28bis} leads to the following integral estimate for
the eigenfrequencies:
\begin{equation*}
\omega^2=\frac{\int_{-\infty}^\infty\beta(u_0+\beta y^2)h_0|\psi|^2\,{\rm d}y+\int_{-\infty}^\infty b_0h_0^2|\psi'|^2\,{\rm d}y}
{\int_{-\infty}^\infty h_0|\psi|^2\,{\rm d}y},
\end{equation*}
which shows that not only the lowest nondimensional (according to the scaling above) eigenfrequency $\bar\omega$ in \eref{28bis} can be
lower than the minimal frequency of the equatorial waves, meaning that these eigenmodes can be trapped inside the jet, but that the
eigenfrequency squared can become negative (if $u_0$ is sufficiently negative), meaning imaginary eigenfrequencies and thus {\em linear
instability}. Such instability is known in the RSW equations in the equatorial beta-plane (see \cite{RZT}, where it was analyzed along the
same lines) as a {\em symmetric inertial instability}. The present analysis shows that it exists in the TRSW model as well. Obviously, the
exponential growth of unstable modes corresponding to imaginary eigenfrequencies rapidly invalidates the linear approximation used in this
section. However, it is known from other studies \cite{BRZ} that the growth of symmetric inertial instability modes in shallow water models 
leads to their breakdown and shock formation in the negative-vorticity (anticyclonic) part of the jet. We therefore expect a similar
scenario here. The analysis of this instability in \cite{RZT} shows that it appears at Rossby numbers of the order unity and Burger numbers
below $\sqrt{1.5}$. Here, we define the Rossby number $Ro$ and Burger number $Bu$ as follows: 
\begin{equation*}
Ro=\frac{U_0}{\beta L^2},\quad Bu=\frac{\sqrt{\bar{b}H_0}}{\beta L^2},
\end{equation*}
where $U_0$ is the maximum velocity of the jet, $L$ is its typical width, $H_0$ is the mean fluid depth, and $\bar{b}$ is the mean buoyancy.

\section{Well-Balanced Semi-Discrete Central-Upwind Scheme}\label{sec3}
In this section, we describe a semi-discrete second-order central-upwind scheme for the 1-D TRSW system \eref{1.10}. To this end, we rewrite
this system in a vector form as
\begin{equation}
\bm{U}_t+\bm{G}(\bm{U},Z)_y=\bm{S}(\bm{U},f),
\label{3.1}
\end{equation}
where
\begin{equation}
\bm{U}=\begin{pmatrix}h\\
q\\
p\\
hb\end{pmatrix},\qquad\bm{G}(\bm{U},Z)=\begin{pmatrix}p\\
\frac{pq}{h}\\
L\\
pb\end{pmatrix},
\qquad\bm{S}(\bm{U},f)=\begin{pmatrix}0\\
fp\\
0\\
0\end{pmatrix}.
\label{3.2}
\end{equation}

We then derive a semi-discretization of \eref{3.1}--\eref{3.2} as follows. We divide the computational domain into a set of uniform cells
$C_k:=[y_\kmh,y_\kph]$, which are centered at $y_k=y_\kmh+\dy/2$. We denote the cell averages of the numerical solutions at time $t$ by
$\,\xbar{\bm{U}}_k(t):\approx\frac{1}{\dy}\int_{C_k}{\bm{U}}(y,t)\,{\rm d}y$ and then integrate the system \eref{3.1}, \eref{3.2} in space
to obtain the following system of ODEs:
\begin{equation}
\frac{{\rm d}}{{\rm d}t}\,\xbar{\bm{U}}_k(t)=-\frac{{\bm{{\cal G}}}_\kph(t)-{\bm{{\cal G}}}_\kmh(t)}{\dy}+\,\xbar{\bm{S}}_k(t).
\label{ODEs}
\end{equation}
Here, $\bm{{\cal G}}_\kph(t)$ are numerical fluxes, which typically depends on the reconstructed left- and right-sided point values of
$\bm{U}$ at the cell interfaces $y=y_\kph$, and
\begin{equation*}
\xbar{\bm{S}}_k(t):\approx\frac{1}{\dy}\int\limits_{C_k}\bm{S}\big(\bm{U}(y,t),f(y)\big)\,{\rm d}y,
\end{equation*}
are the approximations of the cell averages of the source term. When $f$ is taken to be a constant, we simply take
$\,\xbar{\bm{S}}_k:=f\,\xbar{p}_k$. However, when a more realistic case of $f(y)=\beta y$ is considered, we use Simpson rule to obtain
\begin{equation*}
\xbar{\bm{S}}_k=\frac{1}{6}\Big(f(y_\kmh)p_\kmh^++4f(y_k)\,\xbar{p}_k+f(y_\kph)p_\kph^-\Big),
\end{equation*}
where $p_{\kmh}^{\pm}$ and $p_{\kph}^{\pm}$ are the one-sided values of $p$ at the cell interface; these values will be defined in \S\ref{sec31} below.

Details on the computations of $\bm{{\cal G}}_\kph(t)$ are provided in \S\ref{sec32}. For the sake of brevity, we will omit the dependence
of all of the indexed finite-volume quantities on $t$ in the rest of this paper.

\subsection{Well-Balanced Reconstruction}\label{sec31}
It is quite well-known that in order to derive a well-balanced scheme one has to perform piecewise polynomial reconstruction of equilibrium
variables rather than the conservative ones; see, e.g., \cite{CCHKW_18,CK16,CCKOT_18,CDKL,CHO_18,Kur_Acta,KLshw,KPshw}. We therefore
reconstruct the equilibrium variables $\bm{V}:=(q,p,L,b)$.

To this end, we first compute the values $L$ at the cell centers $y=y_k$ as follows. If the cell averages $\{\,\xbar{\bm{U}}_k\}$ are
available at a certain time level $t$, then according to \eref{L}, one obtains
\begin{equation*}
L_k=\frac{\xbar{p}_k^2}{\xbar{h}_k}+\frac{\xbar{(hb)}_k}{2}\,\xbar{h}_k+R_k,\quad k=1,\ldots,N,
\end{equation*}
where $N$ is a total number of cells and $R_k$ can be computed using \eref{R}:
\begin{equation}
R_k=R(y_k,t)=\int\limits^{y_k}\Big[f(\xi)q(\xi,t)+h(\xi,t)b(\xi,t)Z_y(\xi)\Big]{\rm d}\xi,\quad k=1,\ldots,N.
\label{Rj1}
\end{equation}
We notice that formula \eref{Rj1} can be rewritten in the following recursive way:
\begin{equation}
R_k=R_{k-1}+\int\limits_{y_{k-1}}^{y_k}\Big[f(\xi)q(\xi,t)+h(\xi,t)b(\xi,t)Z_y(\xi)\Big]{\rm d}\xi,\quad k=2,\ldots,N,
\label{Rj2}
\end{equation}
and then we apply the following second-order quadrature to the integral in \eref{Rj2}, to obtain
\begin{equation}
R_k=R_{k-1}+\frac{1}{2}(f_{k-1}\,\xbar q_{k-1}+f_k\,\xbar q_k)\dy+\frac{1}{2}\Big(\,\xbar{(hb)}_{k-1}+\,\xbar{(hb)}_k\Big)(Z_k-Z_{k-1}),
\quad k=2,\ldots,N,
\label{Rj3}
\end{equation}
where $f_k:=f(y_k)$.

Similarly, we can use a slightly different quadrature to obtain the point values of $R$ at the cell interfaces:
\begin{equation}
R_\kph=R_\kmh+f_k\,\xbar q_k\dy+\,\xbar{(hb)}_k(Z_\kph-Z_\kmh),\quad k=1,\ldots,N.
\label{Rjph}
\end{equation}
It should be observed that the recursive formulae \eref{Rj3} and \eref{Rjph} require starting values. We first take $R_{1/2}:=0$, then
compute $R_{3/2}$ using \eref{Rjph}, and then set $R_1:=\frac{1}{2}(R_{1/2}+R_{3/2})$.

We also notice that the point values of $Z$ (which are used in \eref{Rj3} and \eref{Rjph}) are obtained as in \cite{KPshw}, namely, we take
\begin{equation*}
Z_{\kph}=\frac{Z(y_{\kph}+0)+Z({y_{\kph}-0})}{2},
\end{equation*}
which reduces to $Z_\kph=Z(y_\kph)$ when $Z$ is continuous. Then, the bottom topography is approximated using a continuous piecewise linear
interpolant
\begin{equation*}
\widetilde Z(y)=Z_\kmh+\frac{Z_\kph-Z_\kmh}{\dy}(y-y_\kmh),\quad y\in C_k,
\end{equation*}
and then take
\begin{equation*}
Z_k:=\widetilde Z(y_k)=\frac{1}{2}(Z_\kph+Z_\kmh).
\end{equation*}

Equipped with the values $\,\xbar{\bm{V}}_k:=(\,\xbar q_k,\,\xbar p_k,L_k,b_k)$, where
\begin{equation}
b_k:=\frac{\xbar{(hb)}_k}{\xbar h_k},
\label{bj}
\end{equation}
we construct a second-order piecewise linear interpolant
\begin{equation}
\widetilde{\bm{V}}(y)=\,\xbar{\bm{V}}_k+{(\bm{V}_y)}_k(y-y_k),\quad y\in C_k,
\label{poly}
\end{equation}
where ${(\bm{V}_y)}_k$ is at least first-order approximation of ${\bm{V}}_y(y_k,t)$. In order to make the reconstruction \eref{poly}
non-oscillatory, we use a nonlinear limiter to compute the slopes ${(\bm{V}_y)}_k$. In the numerical experiments reported in \S\ref{sec4},
we have used the generalized minmod limiter (see, e.g., \cite{LN,NT,Swe,vLeV}):
\begin{equation}
(\bm{V}_y)_k={\rm minmod}\left(\sigma\,\frac{\,\xbar{\bm{V}}_k-\,\xbar{\bm{V}}_{k-1}}{\dy},\,
\frac{\,\xbar{\bm{V}}_{k+1}-\,\xbar{\bm{V}}_{k-1}}{2\dy},\,\sigma\,\frac{\,\xbar{\bm{V}}_{k+1}-\,\xbar{\bm{V}}_k}{\dy}\right),
\label{Vyk}
\end{equation}
where the parameter $\sigma\in[1,2]$ helps to control the amount of numerical diffusion (larger values of $\sigma$ correspond to less
diffusive, but more oscillatory reconstruction), and the minmod function defined by
\begin{equation*}
{\rm minmod}(\alpha_1,\alpha_2,\ldots):=\left\{\begin{array}{lc}\min\limits_k\{\alpha_k\},&\mbox{if~}\alpha_k>0~\forall k,\\
\max\limits_k\{\alpha_k\},&\mbox{if~}\alpha_k<0~\forall k,\\0,&\mbox{otherwise},\end{array}\right.
\end{equation*}
is applied in \eref{Vyk} in a component-wise manner.

Next, using the piecewise polynomial reconstruction \eref{poly}, we obtain the one-sided point values of $\bm{V}$ at the cell interfaces:
\begin{equation*}
\bm{V}_\kph^+=\,\xbar{\bm{V}}_{k+1}-\frac{\dy}{2}(\bm{V}_y)_{k+1},\qquad\bm{V}_\kph^-=\,\xbar{\bm{V}}_k+\frac{\dy}{2}(\bm{V}_y)_k.
\end{equation*}
Finally, the one-sided point values of $h$ at the cell interfaces are computed by solving the following two cubic nonlinear equations, which
arise from the definition of the global variable $L$ in \eref{L}:
\begin{align}
&\Phi(h_\kph^+):=\frac{\big(p_\kph^+\big)^2}{h_\kph^+}+\frac{b_\kph}{2}\big(h_\kph^+\big)^2+R_\kph-L_\kph^+=0,\label{Phih}\\
&\Psi(h_\kph^-):=\frac{\big(p_\kph^-\big)^2}{h_\kph^-}+\frac{b_\kph}{2}\big(h_\kph^-\big)^2+R_\kph-L_\kph^-=0,\label{Psih}
\end{align}
where $R_\kph$ is defined in \eref{Rjph} and
\begin{equation*}
b_\kph:=\frac{1}{2}(b_\kph^++b_\kph^-).
\end{equation*}

Equations \eref{Phih} and \eref{Psih} are solved as in \cite{CCHKW_18}. For example, we describe how to solve \eref{Phih} (the solution of
\eref{Psih} can be obtained in a similar way, we omit here):

If $(p_\kph^+)^4>\frac{8(L_\kph^+-R_\kph)^3}{27b_{\kph}}$, it is easy to show that \eref{Phih} does not have any positive solutions. We
therefore set
\begin{equation}
h_\kph^+=w_\kph^+-Z_\kph,
\label{3.18}
\end{equation}
where $w:=h+Z$ denotes the surface and using a similar way described above to reconstruct $\bm{V}_\kph^+$, we can reconstruct $w_\kph^+$
from $\,\xbar w_k:=\,\xbar h_k+Z_k$. On the other hand, if $(p_\kph^+)^4\le\frac{8(L_\kph^+-R_\kph)^3}{27b_\kph}$, there will be two
possibilities. First, if $p_\kph^+=0$, then \eref{Phih} admits a unique positive solution, namely,
\begin{equation*}
h_\kph^+=\sqrt{\frac{2(L_\kph^+-R_\kph)}{b_\kph}}.
\end{equation*}
Otherwise, we solve \eref{Phih} exactly and obtain the following three solutions:
\begin{equation}
h_\kph^+=2\sqrt{\Upsilon}\cos\Big(\frac{1}{3}[\Theta+2\pi\ell]\Big),\quad\ell=0,1,2,
\label{root}
\end{equation}
where
\begin{equation*}
\Upsilon:=\frac{2(L_\kph^+-R_\kph)}{3b_\kph}\quad\mbox{and}\quad\Theta:=\arccos\bigg(-\frac{(p_\kph^+)^2}{b_\kph\Upsilon^{3/2}}\bigg).
\end{equation*}
One can show that one of these roots is negative, while the other two roots, which correspond to the subsonic and supersonic cases, are
positive. We single out the physically relevant solution by choosing a root in \eref{root} that is closer to the corresponding value of
$h_\kph^+$ given in \eref{3.18}.
\begin{remark}
Notice that in equations \eref{Phih} and \eref{Psih}, we have used $b_\kph$ instead of $b_\kph^\pm$ in order to guarantee that if the
solution is at the steady-state \eref{1.13}, that is, if $L_\kph^+=L_\kmh^-$ and $p_\kph^+=p_\kmh^-=0$, then the values $h_\kph^+$ and
$h_\kmh^-$ coincide. This is important for enforcing a well-balanced property of the resulting scheme; see \S\ref{sec32}.
\end{remark}

\paragraph{Desingularization.} We would like to point out that once the point values $h_\kph^\pm$ and $p_\kph^\pm$ are reconstructed, then
we can obtain the right/left-sided velocities $v_\kph^\pm$, which is needed in the computation of numerical fluxes (see \S\ref{sec32}):
\begin{equation}
v_\kph^\pm=\frac{p_\kph^\pm}{h_\kph^\pm}.
\label{3.20}
\end{equation}
One may, however, observe that the calculation in \eref{3.20} may suffer from one obvious drawback: if the point values $h_\kph^\pm$ become
too small or zero, then \eref{3.20} may not allow one to (accurately) compute $v_\kph^\pm$. In order to ameliorate this fundamental defect,
we use the same desingularization technique which was used in \cite{CCHKW_18,CKL}: we take a small positive number $\varepsilon=10^{-8}$ and set
\begin{equation}
v_\kph^\pm=\frac{2h_\kph^\pm p_\kph^\pm}{\big(h_\kph^\pm\big)^2+\max\left(\big(h_\kph^\pm\big)^2,\varepsilon^2\right)}.
\label{3.19}
\end{equation}
For consistency, we then use \eref{3.19} to recompute the point values $p_\kph^\pm=h_\kph^\pm{\cdot}v_\kph^\pm$. Furthermore, one can notice
that the same problem may occur in the computation of $b_k$ in \eref{bj} if the cell averages $\,\xbar h_k$ become very small or zero.
Therefore, we use the same desingularization technique to adjust this computation, that is, we replace \eref{bj} with
\begin{equation*}
b_k=\frac{2\,\xbar h_k\,\xbar{(hb)}_k}{\xbar h_k^{\,2}+\max\left(\,\xbar h_k^{\,2},\varepsilon^2\right)}.
\end{equation*}

\subsection{Well-Balanced Central-Upwind Numerical Fluxes}\label{sec32}
We first mention that the central-upwind numerical fluxes from \cite{KLin} are given by
\begin{equation}
\bm{{\cal G}}_\kph=\frac{a_\kph^+\bm{G}_\kph^--a_\kph^-\bm{G}_\kph^+}{a_\kph^+-a_\kph^-}+
\frac{a_\kph^+a_\kph^-}{a_\kph^+-a_\kph^-}\left[{\bm{U}_\kph^+-\bm{U}_\kph^-}-\delta\bm{U}_\kph\right],
\label{fluxes}
\end{equation}
where $\bm{U}_\kph^\pm=\big(h_\kph^\pm,q_\kph^\pm,p_\kph^\pm,h_\kph^\pm b_\kph^\pm\big)^\top$ are the left- and right-sided point values of
$\bm{U}$ at the cell interfaces, $\bm{G}_\kph^\pm=\bm{G}\big({\bm{U}}_\kph^\pm,Z_{\kph}\big)$, $a_\kph^\pm$ are the left- and right-sided
local propagation speeds, which can be estimated using the largest and smallest eigenvalues of the Jacobian
$\frac{\partial\bm{F}}{\partial\bm{U}}$ as follows:
\begin{equation*}
\begin{aligned}
a_\kph^+&=\max\left\{v_\kph^-+\sqrt{h_\kph^-b_\kph^-},\,v_\kph^++\sqrt{h_\kph^+b_\kph^+},\,0\right\},\\
a_\kph^-&=\min\left\{v_\kph^--\sqrt{h_\kph^-b_\kph^-},\,v_\kph^+-\sqrt{h_\kph^+b_\kph^+},\,0\right\},
\end{aligned}
\end{equation*}
and $\delta\bm{U}_\kph$ is a built-in ``anti-diffusion'' term:
\begin{equation}\label{d}
\delta\bm{U}_\kph={\rm minmod}\left(\bm{U}_\kph^+-\bm{U}_\kph^*,\,\bm{U}_\kph^*-\bm{U}_\kph^-\right),
\end{equation}
where, as before, the minmod function is applied in a component-wise manner and
\begin{equation}
\bm{U}_\kph^*=\frac{a_\kph^+\bm{U}_\kph^+-a_\kph^-\bm{U}_\kph^--\left\{\bm{G}_\kph^+-\bm{G}_\kph^-\right\}}
{a_\kph^+-a_\kph^-}.
\label{u*}
\end{equation}

We then rewrite the numerical fluxes \eref{fluxes} in a component-wise way:
\begin{equation}
\begin{aligned}
{\cal G}_\kph^{(1)}&=\frac{a_\kph^+p_\kph^--a_\kph^-p_\kph^+}{a_\kph^+-a_\kph^-}+
\frac{a_\kph^+a_\kph^-}{a_\kph^+-a_\kph^-}\left(h_\kph^+-h_\kph^--\delta h_\kph\right),\\
{\cal G}_\kph^{(2)}&=\frac{a_\kph^+q_\kph^-v_\kph^--a_\kph^-q_\kph^+v_\kph^+}{a_\kph^+-a_\kph^-}+
\frac{a_\kph^+a_\kph^-}{a_\kph^+-a_\kph^-}\left(q_\kph^+-q_\kph^--\delta q_\kph\right),\\
{\cal G}_\kph^{(3)}&=\frac{a_\kph^+L_\kph^--a_\kph^-L_\kph^+}{a_\kph^+-a_\kph^-}+
\frac{a_\kph^+a_\kph^-}{a_\kph^+-a_\kph^-}\left(p_\kph^+-p_\kph^--\delta p_\kph\right),\\
{\cal G}_\kph^{(4)}&=\frac{a_\kph^+p_\kph^-b_\kph^--a_\kph^-p_\kph^+b_\kph^+}{a_\kph^+-a_\kph^-}\\
&+\frac{a_\kph^+a_\kph^-}{a_\kph^+-a_\kph^-}\left(h_\kph^+b_\kph^+-h_\kph^-b_\kph^--\delta(hb)_\kph\right).
\end{aligned}
\label{G}
\end{equation}
One may now see that if the central-upwind fluxes \eref{G} are used, the steady states \eref{1.13} would not be preserved at the discrete
level. Indeed, if the discrete data satisfy $L_\kph^+=L_\kph^-$ and $p_\kph^+=p_\kph^-=0$, then it follows from \eref{Phih} and \eref{Psih}
that $h_\kph^+=h_\kph^-$ and \eref{d} and \eref{u*} imply that $\delta h_\kph\equiv0$ and $\delta p_\kph\equiv0$. Therefore, in this case,
the components ${{{\cal G}}}_\kph^{(1)}$ and ${{{\cal G}}}_\kph^{(3)}$ would vanish. However, the second  and fourth components,
${{{\cal G}}}_\kph^{(2)}$ and ${{{\cal G}}}_\kph^{(4)}$, do not necessarily vanish since $q$ and $b$ are not constant at these steady
states, so that $q_\kph^+-q_\kph^-$, $\delta q_\kph$, $h_\kph^+b_\kph^+-h_\kph^-b_\kph^-$ and $\delta(hb)_\kph$ are not in general zero.

\begin{figure}[ht!]
\centerline{\includegraphics[width=7cm]{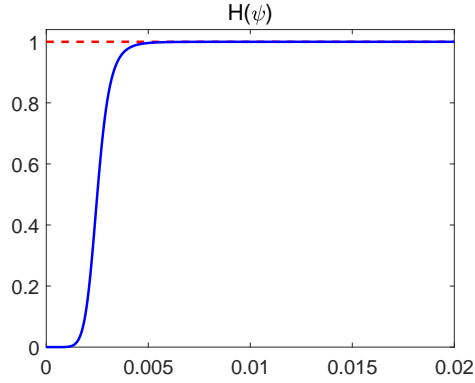}}
\caption{\sf Sketch of $H(\psi)$.\label{switch}}
\end{figure}

We therefore follow the idea from \cite{CCKOT_18} and modify ${\cal G}_\kph^{(2)}$ and ${\cal G}_\kph^{(4)}$ by adding a ``diffusion
switch'' function:
\begin{equation}
\begin{aligned}
{\cal G}_\kph^{(2)}&=\frac{a_\kph^+q_\kph^-v_\kph^--a_\kph^-q_\kph^+v_\kph^+}{a_\kph^+-a_\kph^-}+
H\big(\psi_\kph\big)\frac{a_\kph^+a_\kph^-}{a_\kph^+-a_\kph^-}\left(q_\kph^+-q_\kph^--\delta q_\kph\right),\\
{\cal G}_\kph^{(4)}&=\frac{a_\kph^+p_\kph^-b_\kph^--a_\kph^-p_\kph^+b_\kph^+}{a_\kph^+-a_\kph^-}\\
&+H\big(\psi_\kph\big)\frac{a_\kph^+a_\kph^-}{a_\kph^+-a_\kph^-}\left(h_\kph^+b_\kph^+-h_\kph^-b_\kph^--\delta(hb)_\kph\right).
\end{aligned}
\label{G2}
\end{equation}
Here, the smooth cut-off function $H(\psi)$ is defined as
\begin{equation*}
H(\psi)=\frac{(C\psi)^m}{1+(C\psi)^m},
\end{equation*}
with the constants $C=400$ and $m=8$ used in all of the numerical experiments reported in \S\ref{sec4}. We plot a sketch of this function in
Figure \ref{switch}, where one can clear see that $H(0)=0$ and if $\psi$ is small, then the value of $H$ is still very close to $0$. When
$\psi$ increases, the values of $H$ rapidly approaches $1$.

Following \cite{CCKOT_18}, we take
\begin{equation*}
\psi_\kph:=\frac{|L_{k+1}-L_k|}{\dy}\cdot\frac{y_{N+\frac{1}{2}}-y_{\frac{1}{2}}}{\max(L_k,L_{k+1})}.
\end{equation*}

Finally, the presented central-upwind scheme should preserve the positivity of $h$ and $b$. This is achieved by implementing a ``draining
time step'' technique, which was introduced in \cite{BNL}; see also \cite{BCKN,CCHKW_18,CK16}.
\begin{remark}
We note that the presented semi-discrete central-upwind scheme \eref{ODEs}, \eref{fluxes} with the modified second and fourth components
${\bm{{\cal G}}}_\kph^{(2)}$ and ${\bm{{\cal G}}}_\kph^{(4)}$ given by \eref{G2} is a system of time dependent ODEs, which should be
integrated in time by a sufficiently accurate, efficient and stable ODE solver. In our numerical experiments, we have used the three-stage
third-order strong stability preserving (SSP) Runge-Kutta method (see, e.g., \cite{GKS,GST}) with an adaptive time step computed at every
time level using the CFL number 1/2:
\begin{equation*}
\dt=\frac{\dx}{2a_{\max}},\qquad a_{\max}:=\max_k\left\{a_\kph^+,-a_\kph^-\right\}.
\end{equation*}
\end{remark}

\section{Numerical Examples}\label{sec4}
In this section, we demonstrate the performance of the proposed semi-discrete second-order well-balanced central-upwind scheme on several
numerical examples. In all of the experiments, we take the minmod parameter $\sigma=1.3$ and the boundary conditions are set to be a
zero-order extrapolation at both sides of the computational domain. For the sake of brevity, the proposed well-balanced central-upwind
scheme will be referred to as the WB-CU scheme.

In Examples 1 and 2, we test the ability of the WB-CU scheme to cope with nontrivial topography, and consider a non-rotational case, that
is, the Coriolis parameter is taken to be $f\equiv0$ there. In this case, the zonal discharge $q$ does not need to be considered and we thus
numerically solve the system
\begin{equation*}
\left\{\begin{aligned}
&h_t+p_y=0,\\
&p_t+L_y=0,\\
&(hb)_t+(pb)_y=0
\end{aligned}\right.
\end{equation*}
instead of \eref{1.10}.

In Examples 3 and 4, we test how the WB-CU scheme copes with rotation with a constant Coriolis parameter $f(y)\equiv1$. In Example 4, we
demonstrate the breakdown of smooth solutions. Finally, in Examples 5 and 6, we consider the case of a variable Coriolis parameter
$f(y)=0.1y$. In Examples 3--6, the bottom topography is taken to be flat $(Z\equiv0)$.

\subsection*{Example 1 --- Small Perturbation of a Steady-State Solution}
In the first example taken from \cite{CKL}, we study an ability of the proposed WB-CU scheme to handle a small perturbation of the following discontinuous steady states:
\begin{equation}
(h_s+Z,p_s,b_s)^\top(y,0)=\left\{\begin{aligned}&(6,0,4)^\top,\quad y<0,\\&(4,0,9)^\top,\quad y>0,\end{aligned}\right.
\label{Ex1}
\end{equation}
with the nonflat bottom topography that contains two isolated humps:
\begin{equation*}
Z(y)=\left\{\begin{array}{lc}
0.85(\cos(10\pi(y+0.9))+1),&-1\le y\le-0.8,\\
1.25(\cos(10\pi(y-0.4))+1),&0.3\le y\le0.5,\\
0,&\mbox{otherwise}.
\end{array}\right.
\end{equation*}
The initial data,
\begin{equation*}
(h+Z,p,b)^\top(y,0)=
(h_s+Z,p_s,b_s)^\top(y)+\left\{\begin{array}{lc}(0.1,0,0),&-1.5\le y\le-1.4,\\(0,0,0),&\mbox{otherwise},\end{array}\right.
\end{equation*}
are small perturbation of the steady state \eref{Ex1} and the computational domain is $[-2,2]$.
\begin{figure}[ht!]
\centerline{\includegraphics[width=6.1cm]{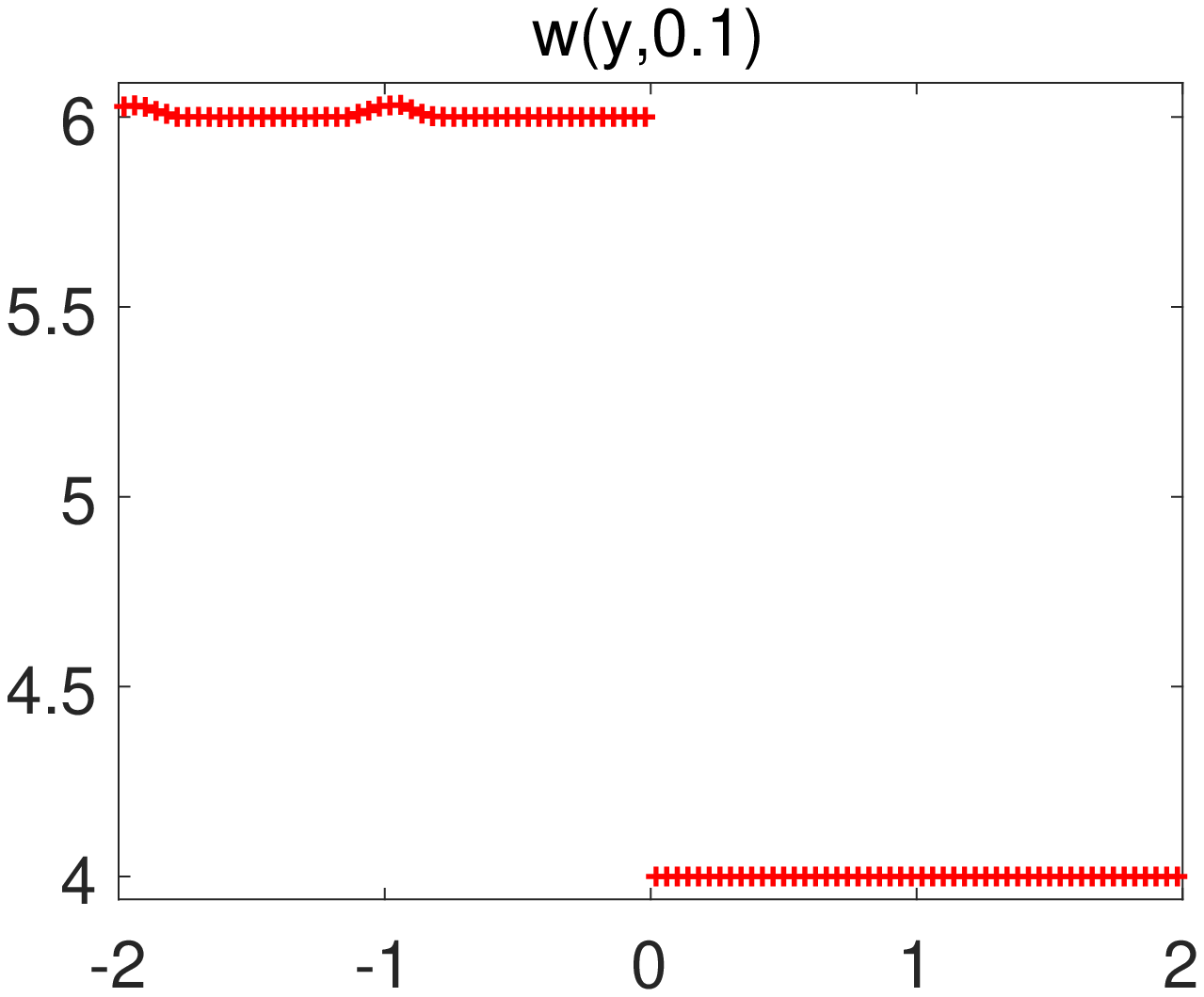}\includegraphics[width=6.1cm]{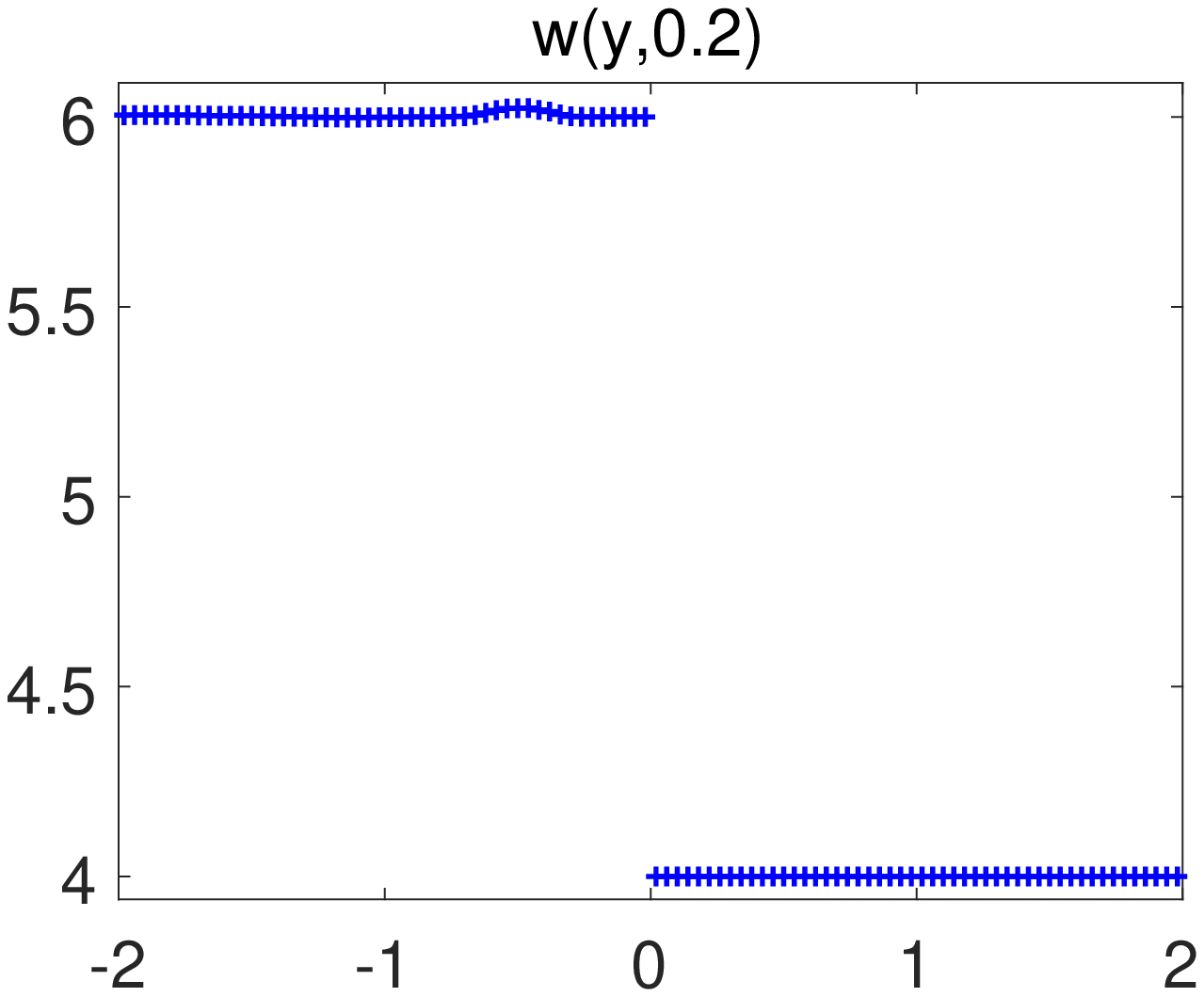}\includegraphics[width=6.1cm]{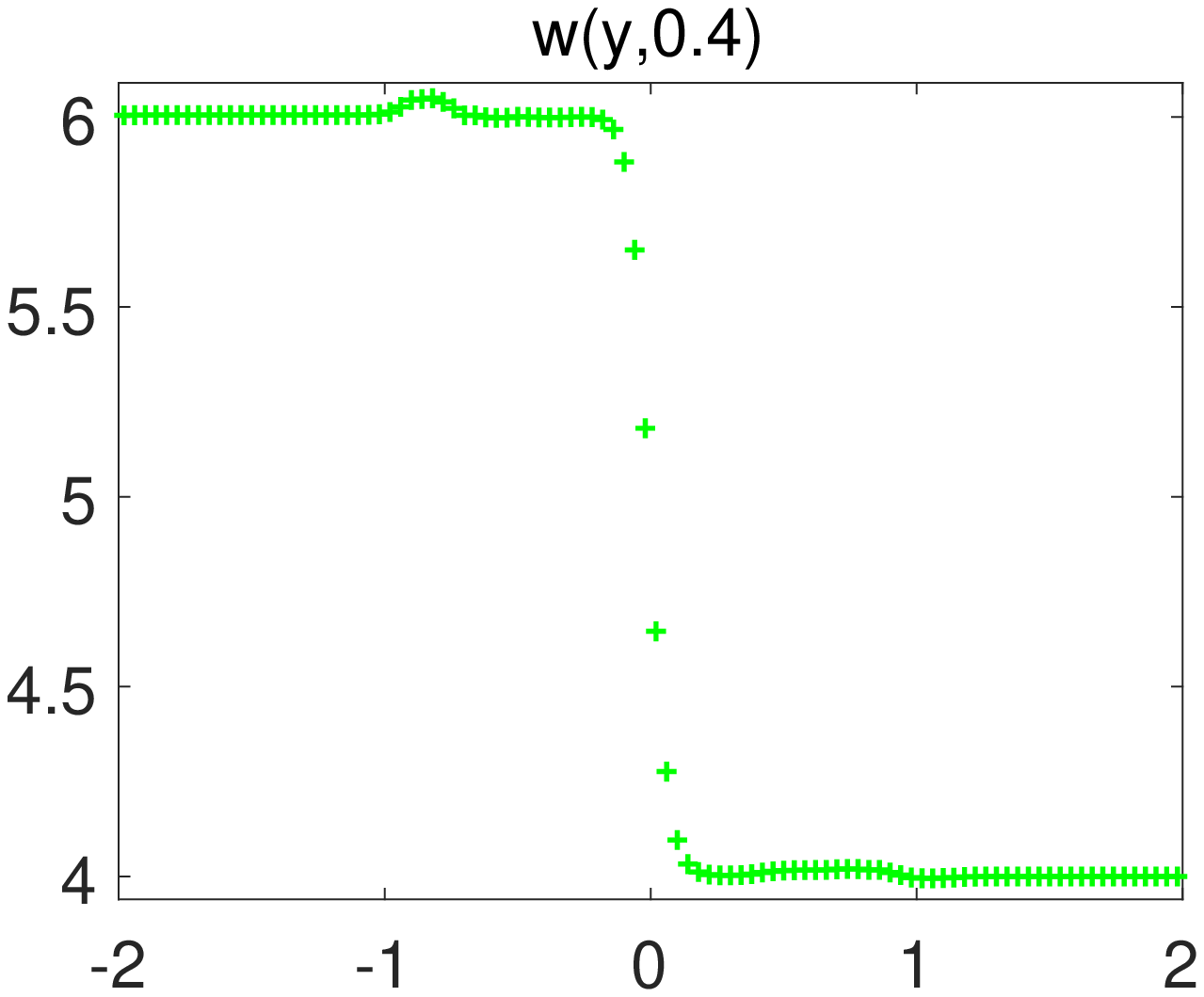}}
\caption{\sf Example 1: Water surface $w$ computed by the WB-CU scheme at different times.
\label{fig1}}
\end{figure}

We compute the solution using the uniform mesh with $\dy=0.04$ until the final time $t=0.4$. The time snapshots of the computed water
surface $w=h+Z$ are plotted in Figure \ref{fig1}. As one can see, the perturbation, initially located at $[-1.5,-1.4]$, splits into two
pulses moving into the opposite directions. The one moving to the right passes over the first and then the second hump of the bottom.
Compared to non-well-balanced results reported in \cite{CKL}, no oscillations are developed by the proposed WB-CU scheme. Moreover, at times
$t=0.1$ and $0.2$ the jump in $w$ remains almost perfectly resolved. At later time $t=0.4$, after the perturbation passed the jump at $y=0$,
the computed solution is still non-oscillatory, but the jump is smeared. In fact, it is less smeared compared to the well-balanced scheme
from \cite{CKL}, but no as sharp as the jump obtained in \cite{CKL} using a special interface tracking technique.

\subsection*{Example 2 --- Dam-Break over a Nonflat Bottom}
In the second example also taken from \cite{CKL}, we numerically solve the dam-break problem with a nonflat bottom topography given by
\begin{equation*}
Z(y)=\left\{\begin{array}{lc}
2(\cos(10\pi(y+0.3))+1),&-0.4\le y\le-0.2,\\
0.5(\cos(10\pi(y-0.3))+1),&0.2\le y\le0.4,\\
0,&\mbox{otherwise}.
\end{array}\right.
\end{equation*}
The initial data are
\begin{equation*}
(w,u,b)^\top(y,0)=\left\{\begin{aligned}&(5,0,1)^\top,&&y<0,\\&(1,0,5)^\top,&&y>0.\end{aligned}\right.
\end{equation*}
\begin{figure}[ht!]
\centerline{\includegraphics[width=7.0cm]{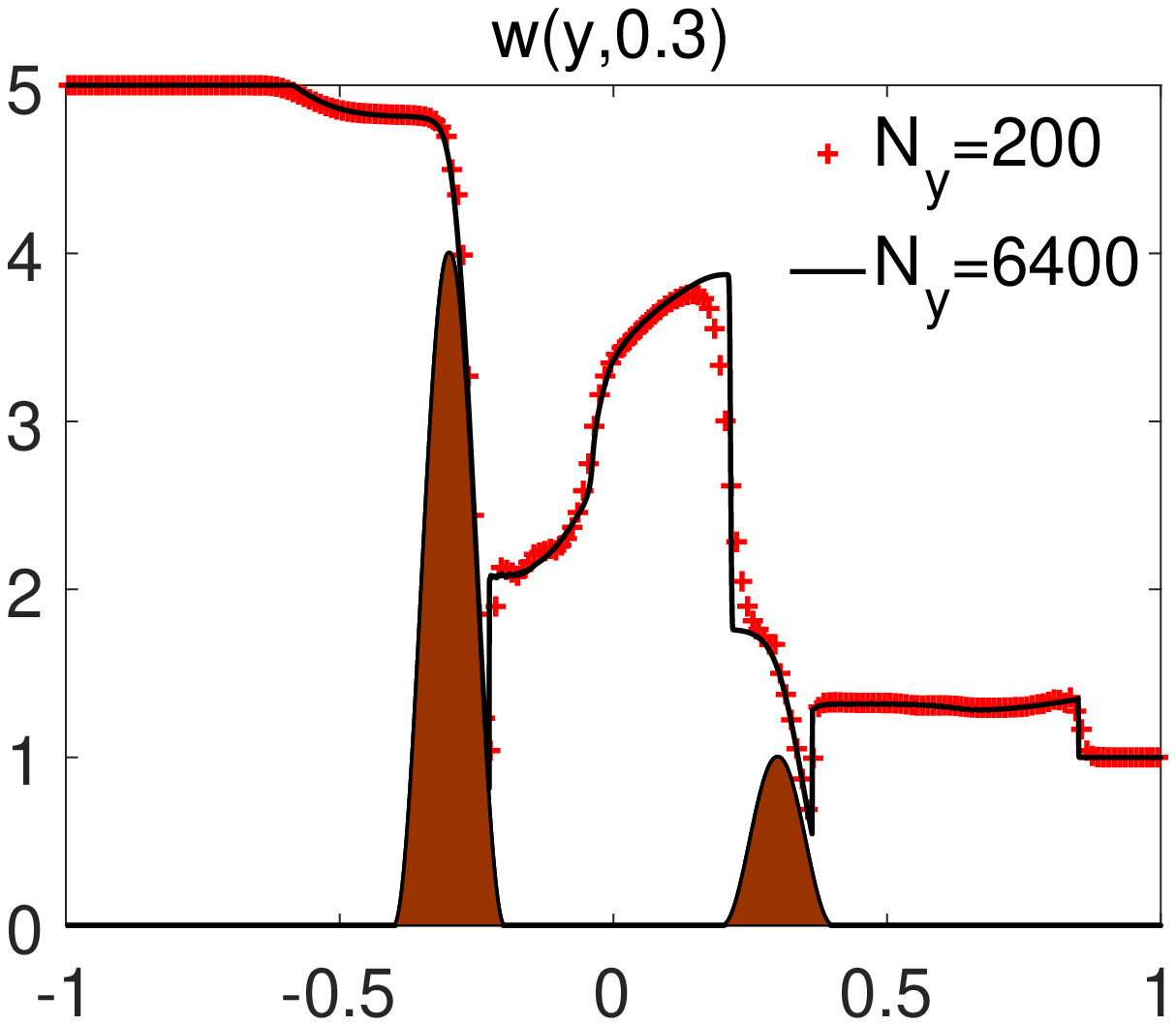}\hspace*{0.5cm}\includegraphics[width=7.0cm]{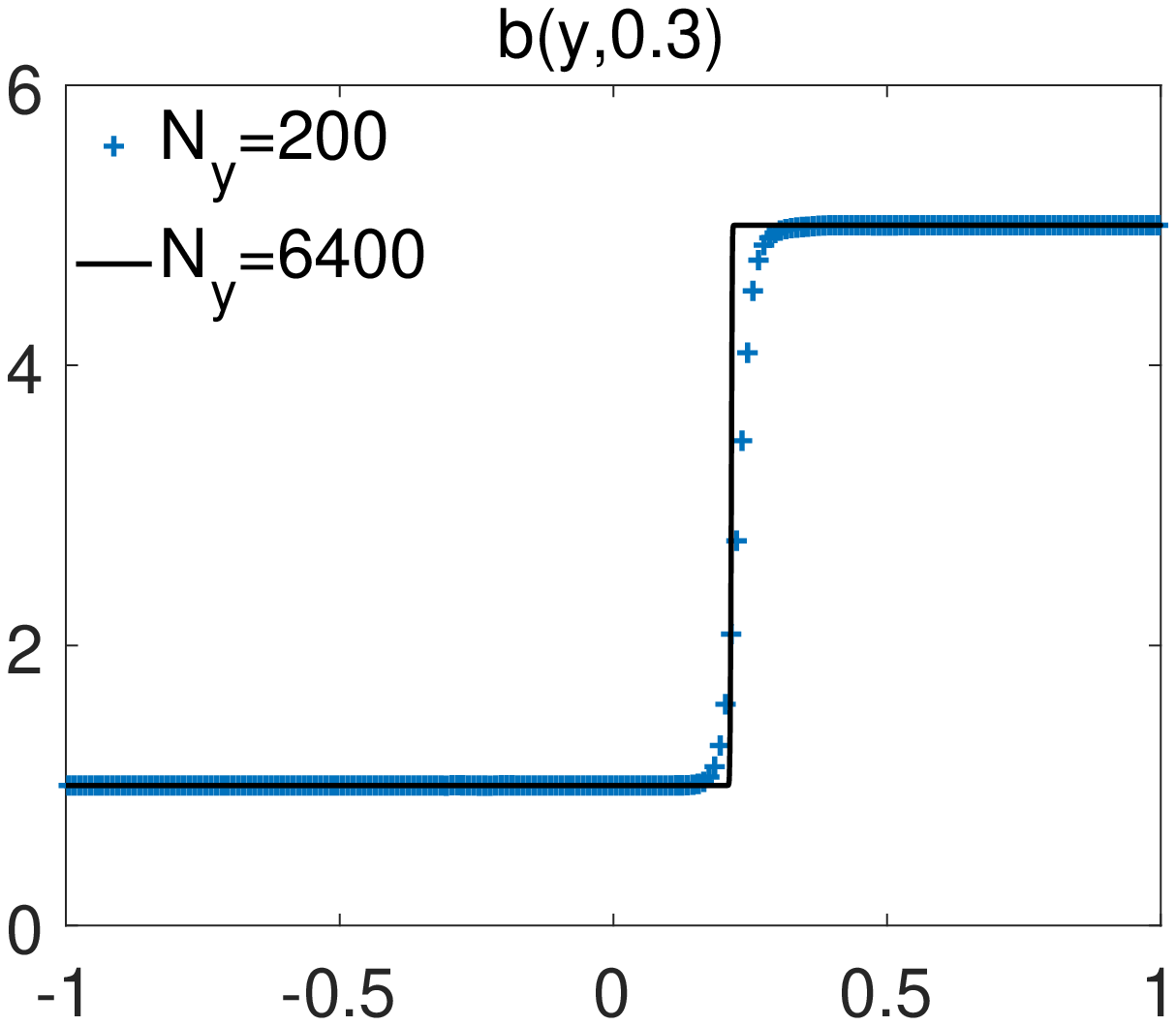}}
\caption{\sf Example 2: Water surface $w$ and buoyancy $b$ computed by the WB-CU scheme.
\label{fig2}}
\end{figure}

We solve the underlying problem using the WB-CU method with $\dy=0.01$. We also compute the reference solutions on a finer mesh with
$\dy=2/6400$. In Figure \ref{fig2}, we show the solution ($w=h+Z$ and $b$) at time $t=0.3$. As one can see, the WB-CU scheme can preserve
the positivity of $h$, which is quite small in this example. Moreover, comparing these results with those reported in \cite{CKL}, one can
conclude that the proposed WB-CU method produces more accurate results than those obtained by the well-balanced scheme used in \cite{CKL}.
However, our results are not as sharp as those obtained in \cite{CKL} using a special interface tracking technique. We note that this
technique can be incorporated into our WB-CU scheme. This may be needed when $b$ is discontinuous. It is worth emphasizing in this context
that sharp temperature fronts do occur in the atmosphere and ocean.

\subsection*{Example 3 --- Rossby Adjustment in an Open Domain in the $f$-plane}
In the third example, we numerically investigate the Rossby adjustment problem with the constant Coriolis parameter $f\equiv1$, which was
studied previously in \cite{BSZ,CDKL}. We consider the following initial conditions, which correspond to a jet over a flat $h$:
\begin{equation*}
h(y,0)\equiv1,\quad v(y,0)\equiv0,\quad u(y,0)=\frac{2(1+\tanh(2y+2))(1-\tanh(2y-2))}{(1+\tanh(2))^2},
\end{equation*}
which are prescribed in the computational domain $[-250,250]$. We take the following three sets of initial values of $b$, which correspond
to zero, positive and negative gradients of $b(y,0)$, respectively (with the first case being a ``pure'' RSW one and thus a benchmark):
$$
\mbox{(a)}~~b(y,0)\equiv1;\qquad\mbox{(b)}~~b(y,0)=1+\frac{1}{10}\tanh(0.5y);\qquad\mbox{(c)}~~b(y,0)=1-\frac{1}{10}\tanh(0.5y).
$$

\begin{figure}[ht!]
\centerline{\includegraphics[width=8cm]{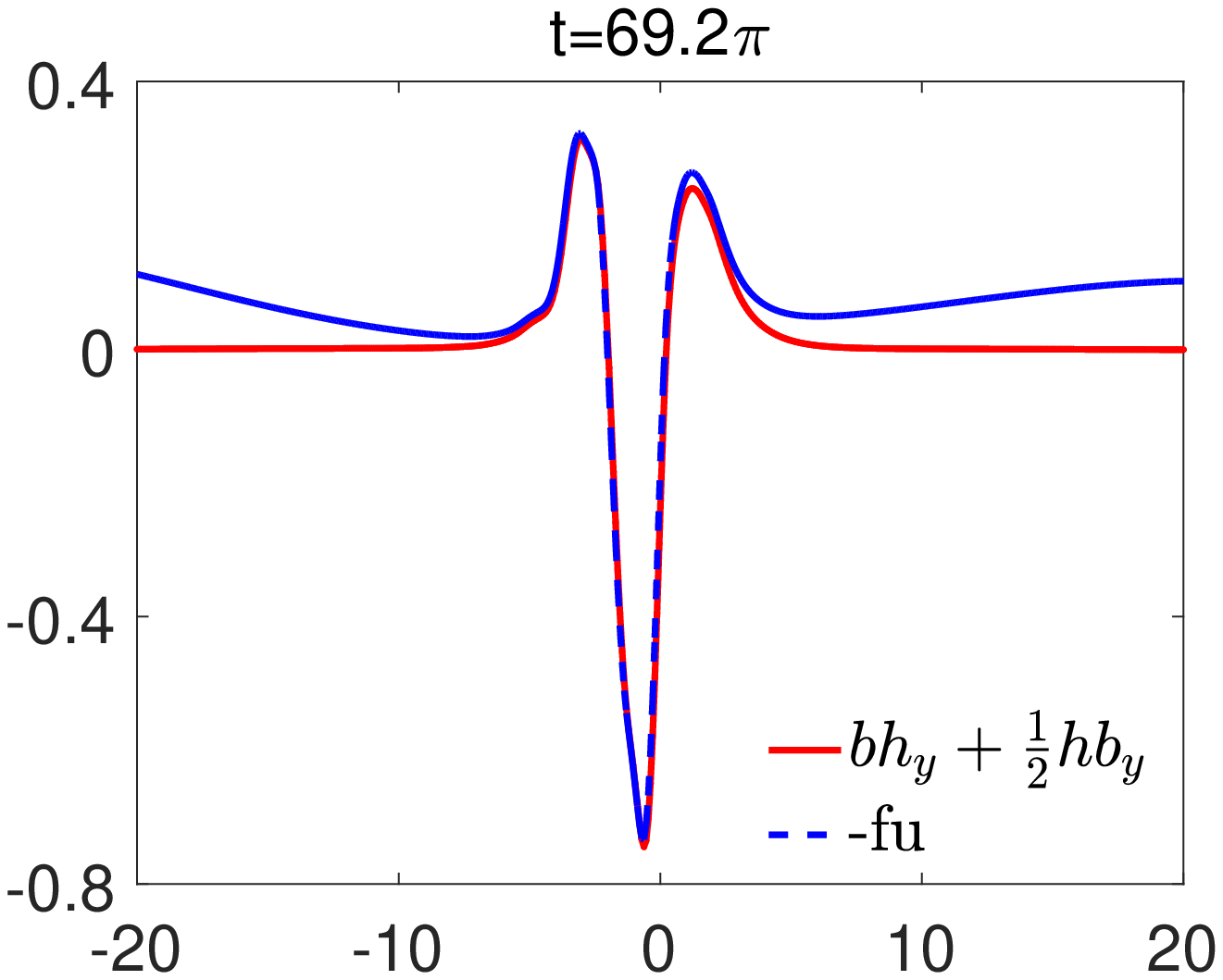}\hspace*{0.5cm}\includegraphics[width=8cm]{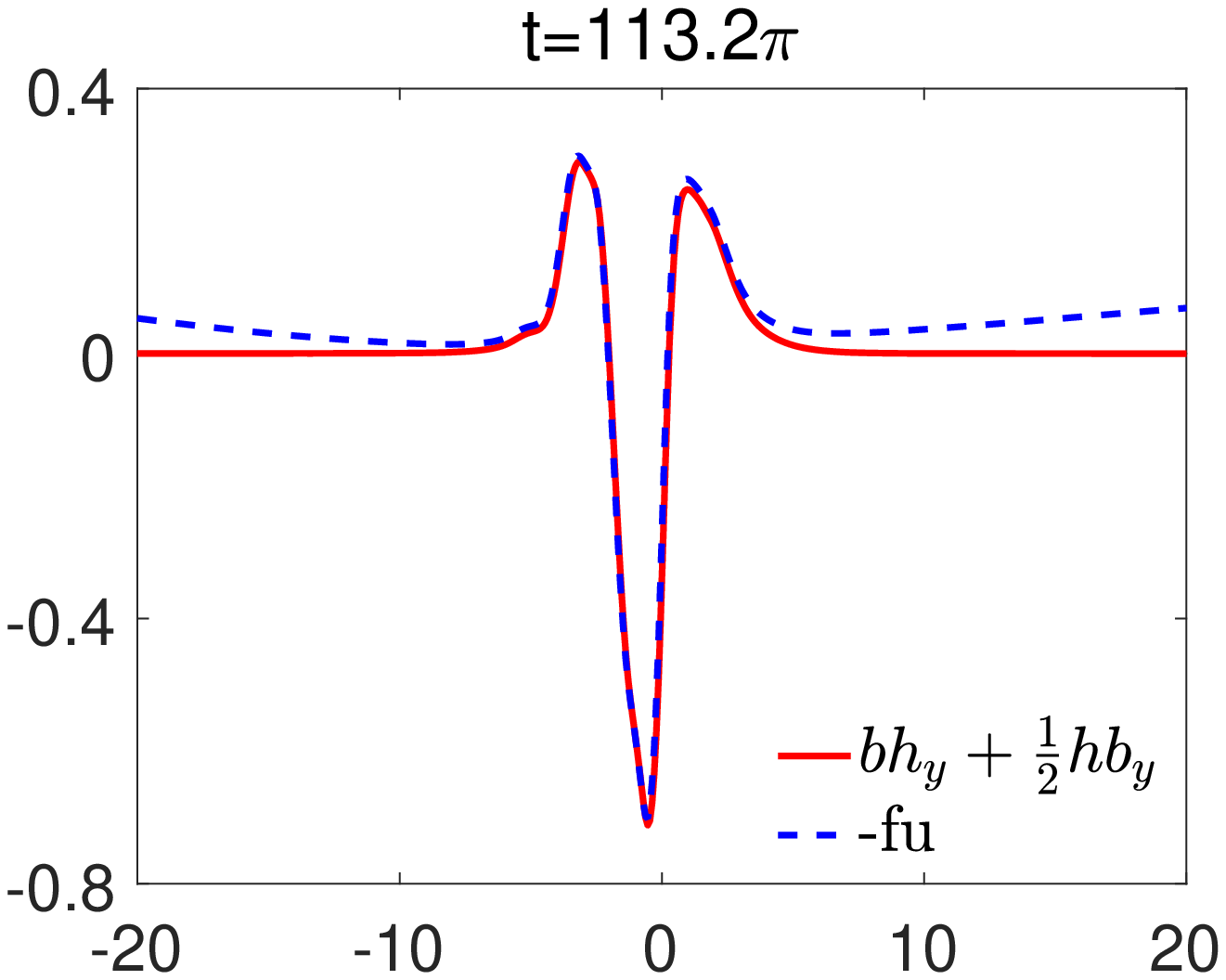}}
\caption{\sf Example 3, case (b): Snapshots of $bh_y+\frac{1}{2}hb_y$ and $-fu$ computed by the WB-CU scheme using $N=6000$ finite-volume
cells at late stages of the Rossby adjustment on the $f$-plane.\label{fig30}}
\end{figure}

\begin{figure}[ht!]
\centerline{\includegraphics[width=8cm]{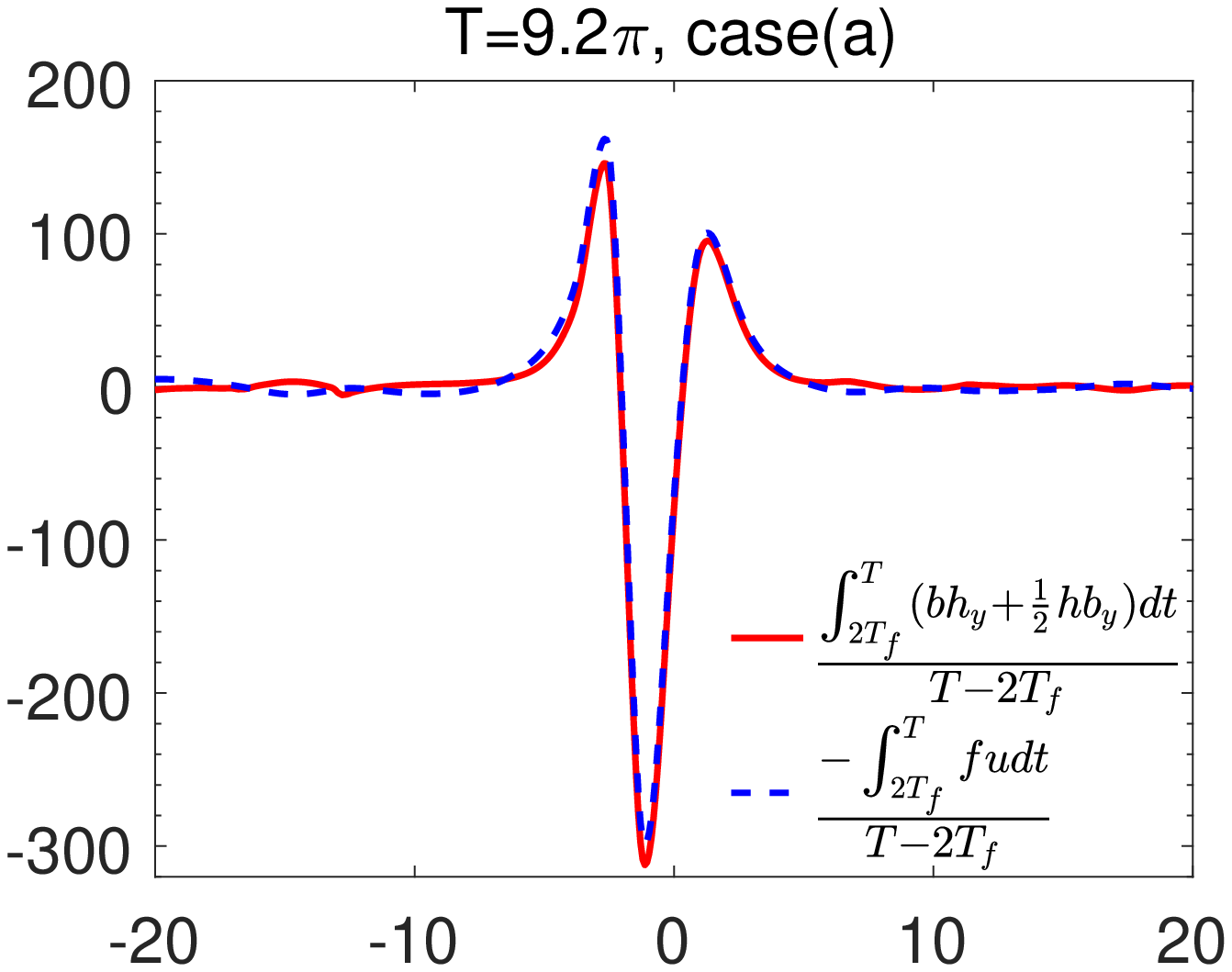}\hspace*{0.3cm}\includegraphics[width=8cm]{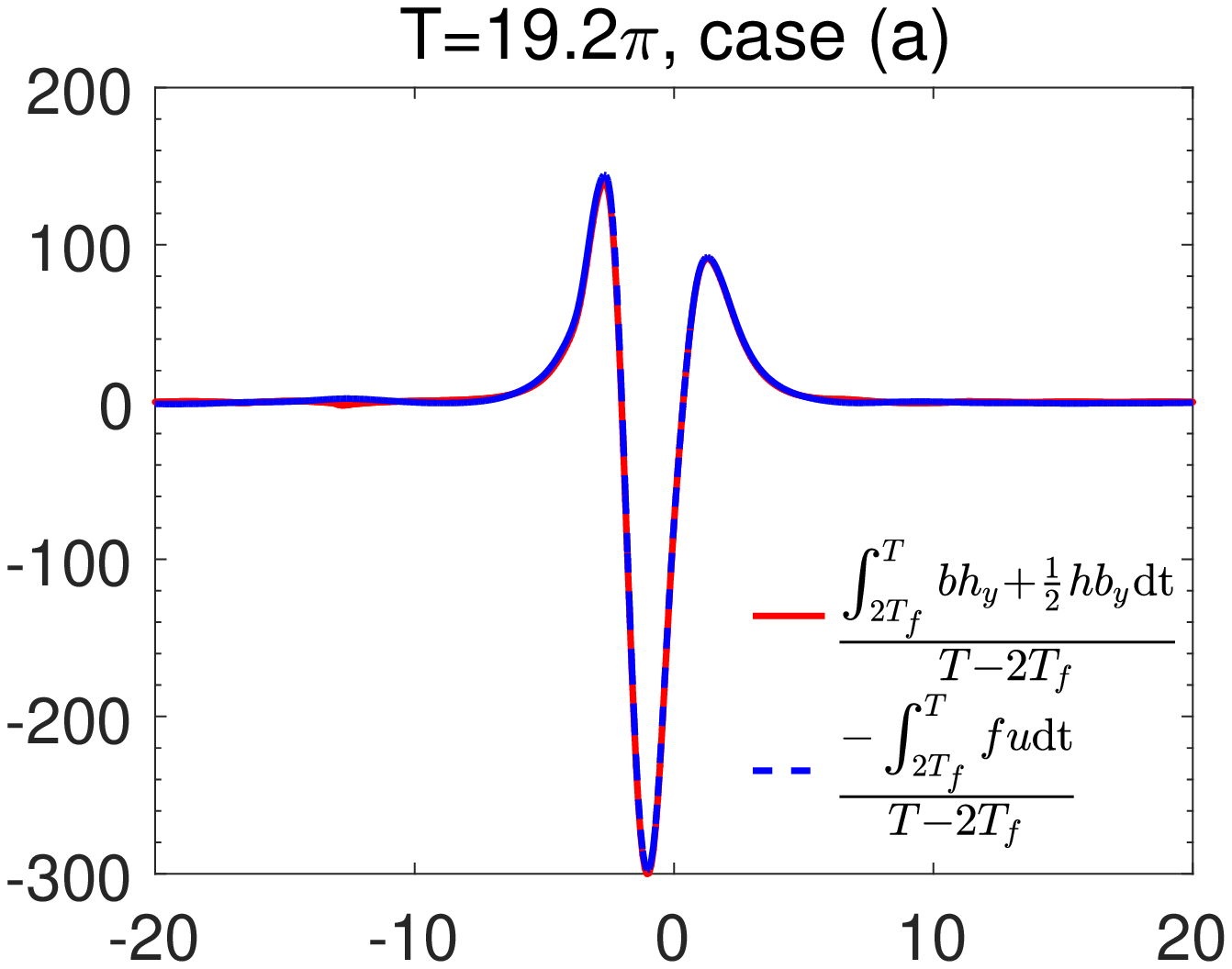}}

\vspace*{0.2cm}
\centerline{\includegraphics[width=8cm]{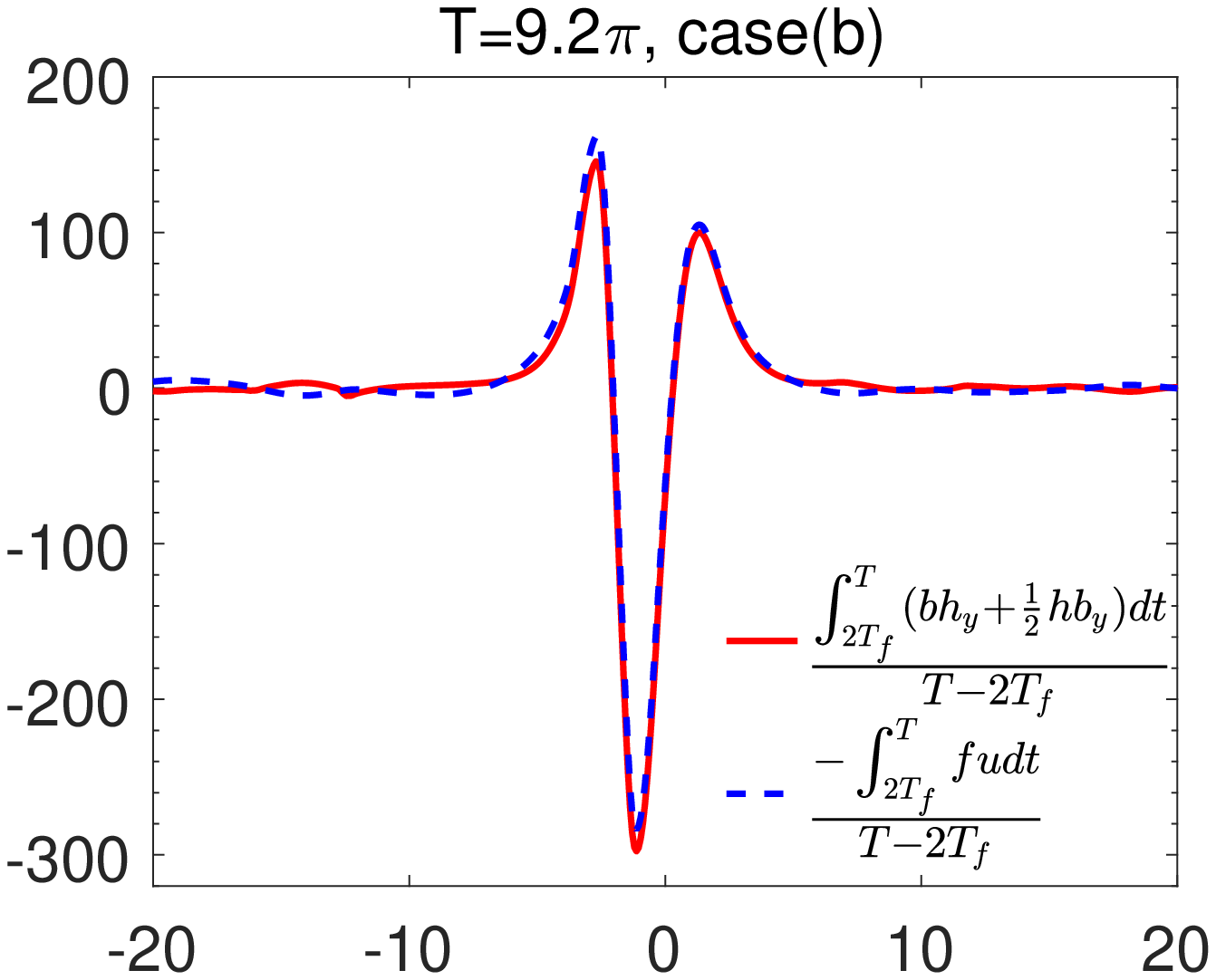}\hspace*{0.3cm}\includegraphics[width=8cm]{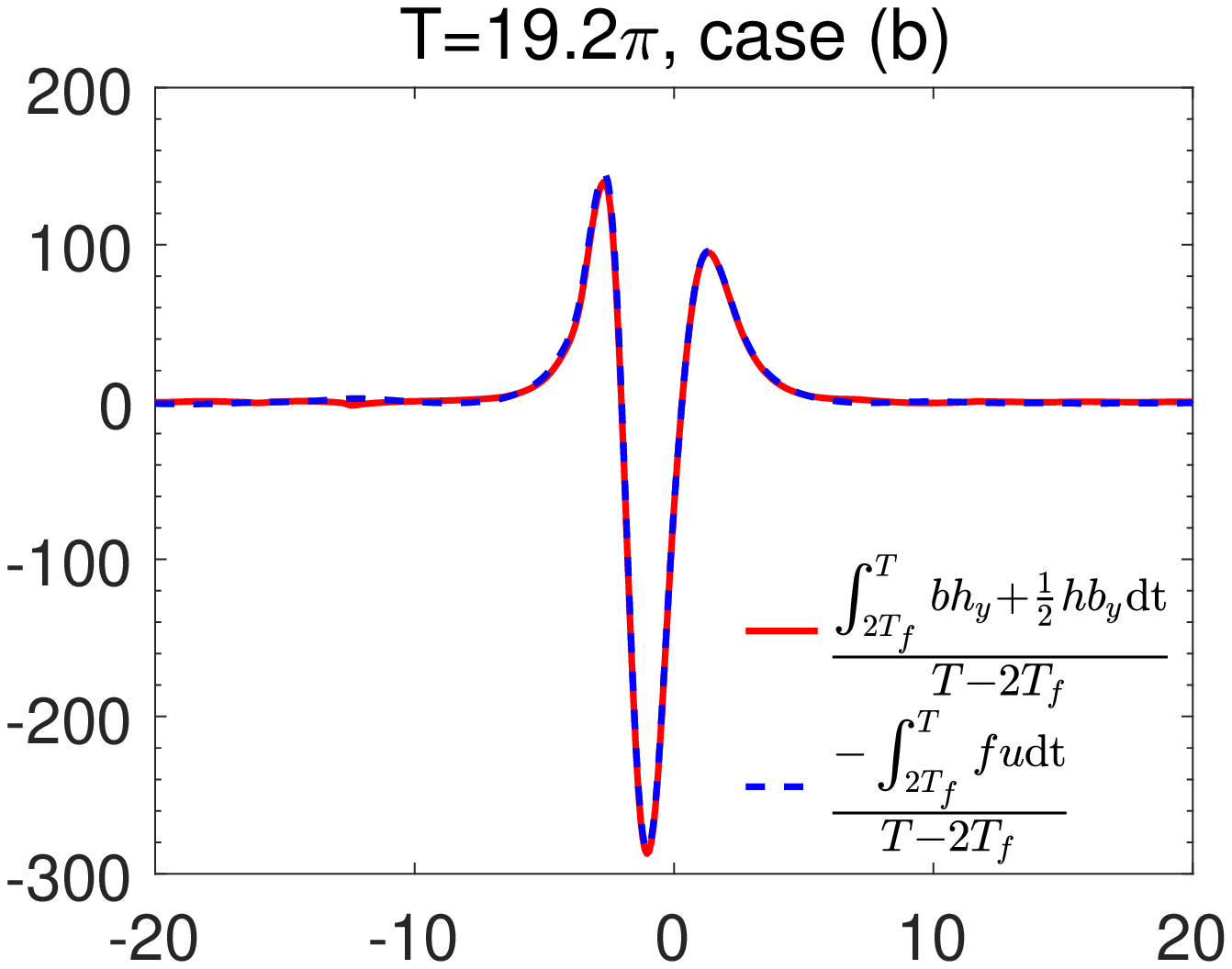}}

\vspace*{0.2cm}
\centerline{\includegraphics[width=8cm]{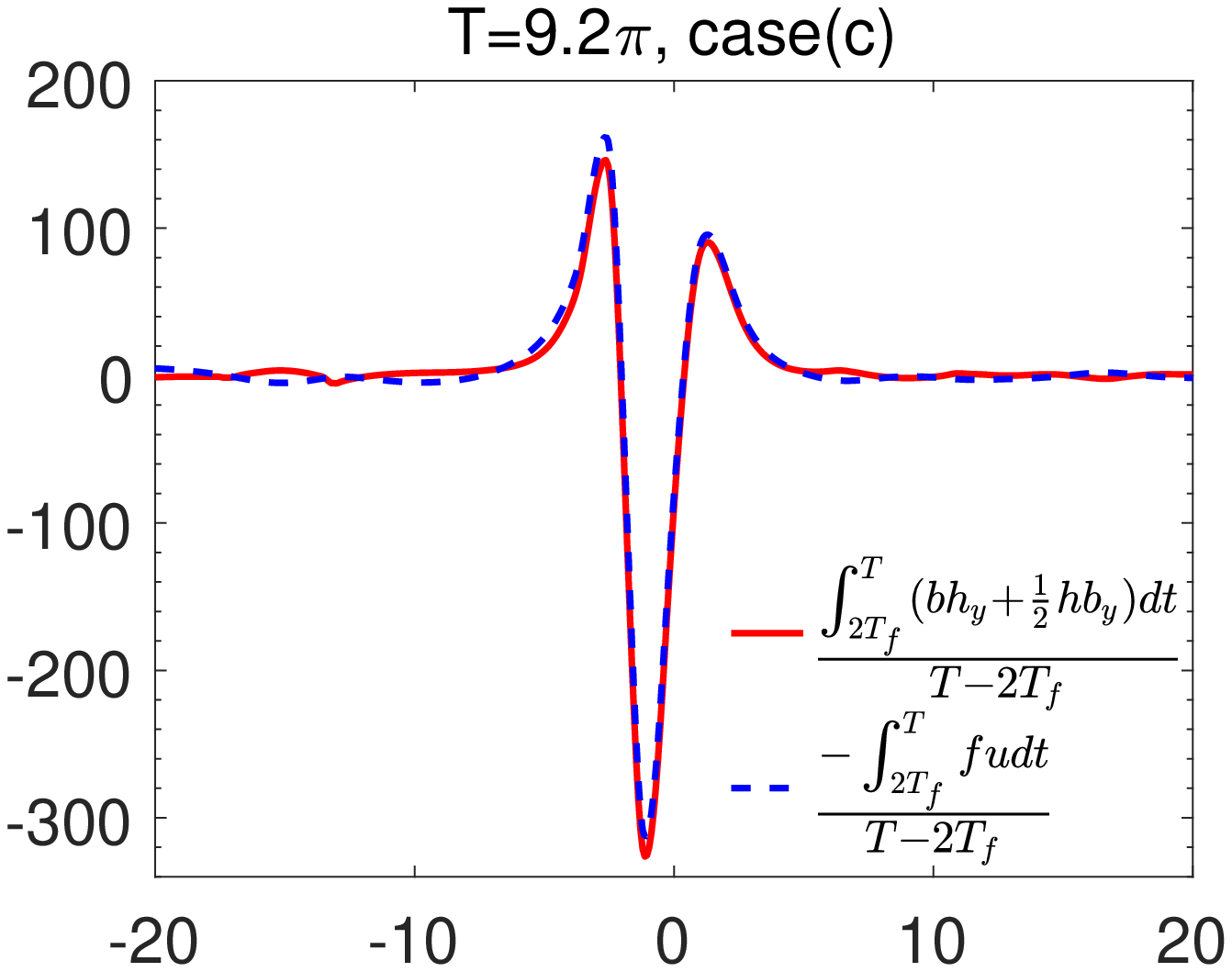}\hspace*{0.3cm}\includegraphics[width=8cm]{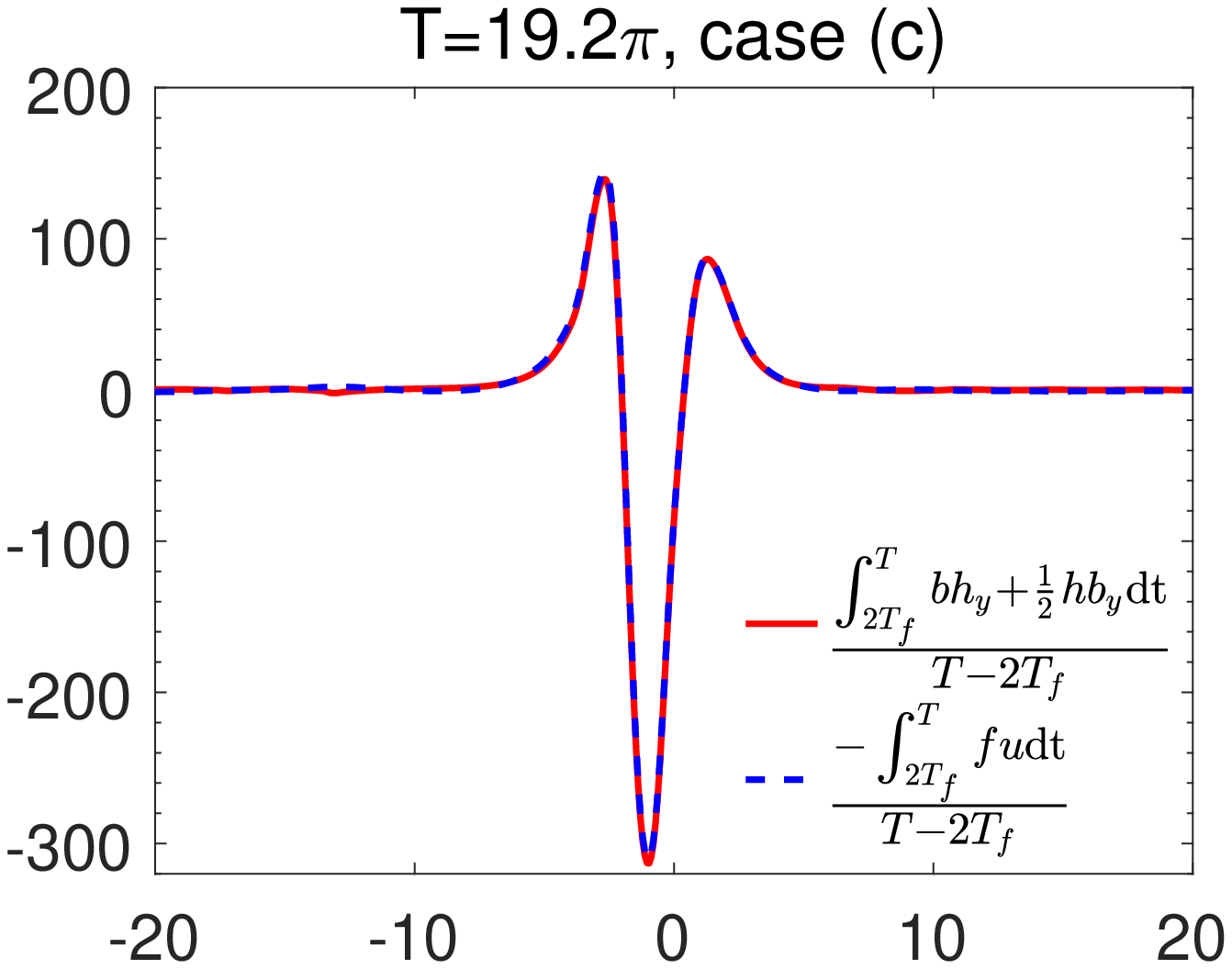}}
\caption{\sf Example 3: Evolution of $\frac{\int_{2T_f}^T\Big(bh_y+\frac{1}{2}hb_y\Big)\,{\rm d}t}{T-2T_f}$ and
$\frac{-\int_{2T_f}^Tfu\,{\rm d}t}{T-2T_f}$ computed by the WB-CU scheme using $N=6000$ finite-volume cells for different $b$ (cases (a),
(b) and (c)) during the Rossby  adjustment on the $f$-plane.\label{fig31}}
\end{figure}

In this example, the bottom topography is flat $(Z\equiv0)$ and thus the thermo-geostrophic equilibrium \eref{1.13} can be rewritten as
\begin{equation}
bh_y+\frac{1}{2}hb_y=-fu.
\label{ge}
\end{equation}
We therefore measure the quantities on the left- side right-hand sides of \eref{ge}, which are supposed to be the same at the steady state,
but remain quite different even at relatively large times $t=69.2\pi$ and $113.2\pi$; see Figure \ref{fig30}, and also the results reported
in \cite{BSZ,CDKL}. We notice that in this example, a natural time scale is the inertial period $T_f=2\pi/f=2\pi$ so that one could have
expected the solution to be very close to the thermo-geostrophic equilibrium \eref{ge} by $t=113.2\pi$. This, however, does not happen
since, as was explained in \S\ref{sec24}, some of the wave modes have almost zero group velocity and thus stay in the core of the jet for
a long time. Under such circumstances, where the fast component of motion is still present, the thermo-geostrophic balance could be
satisfied only for the slow, time-averaged component. We therefore take the time averages in \eref{ge},
\begin{equation}
\int\limits_{2T_f}^T\Big(bh_y+\frac{1}{2}hb_y\Big)\,{\rm d}t=-\int\limits_{2T_f}^Tfu\,{\rm d}t,
\label{ge2}
\end{equation}
and measure the left- and right-hand sides of \eref{ge2} at large $T$. The obtained results are shown in Figure \ref{fig31} for $T=9.2\pi$
and $19.2\pi$. As one can see, in all of the considered cases the computed solutions satisfy the relationship \eref{ge2} very well, even
when $T$ is not very large.

\begin{figure}[ht!]
\centerline{\includegraphics[width=7cm]{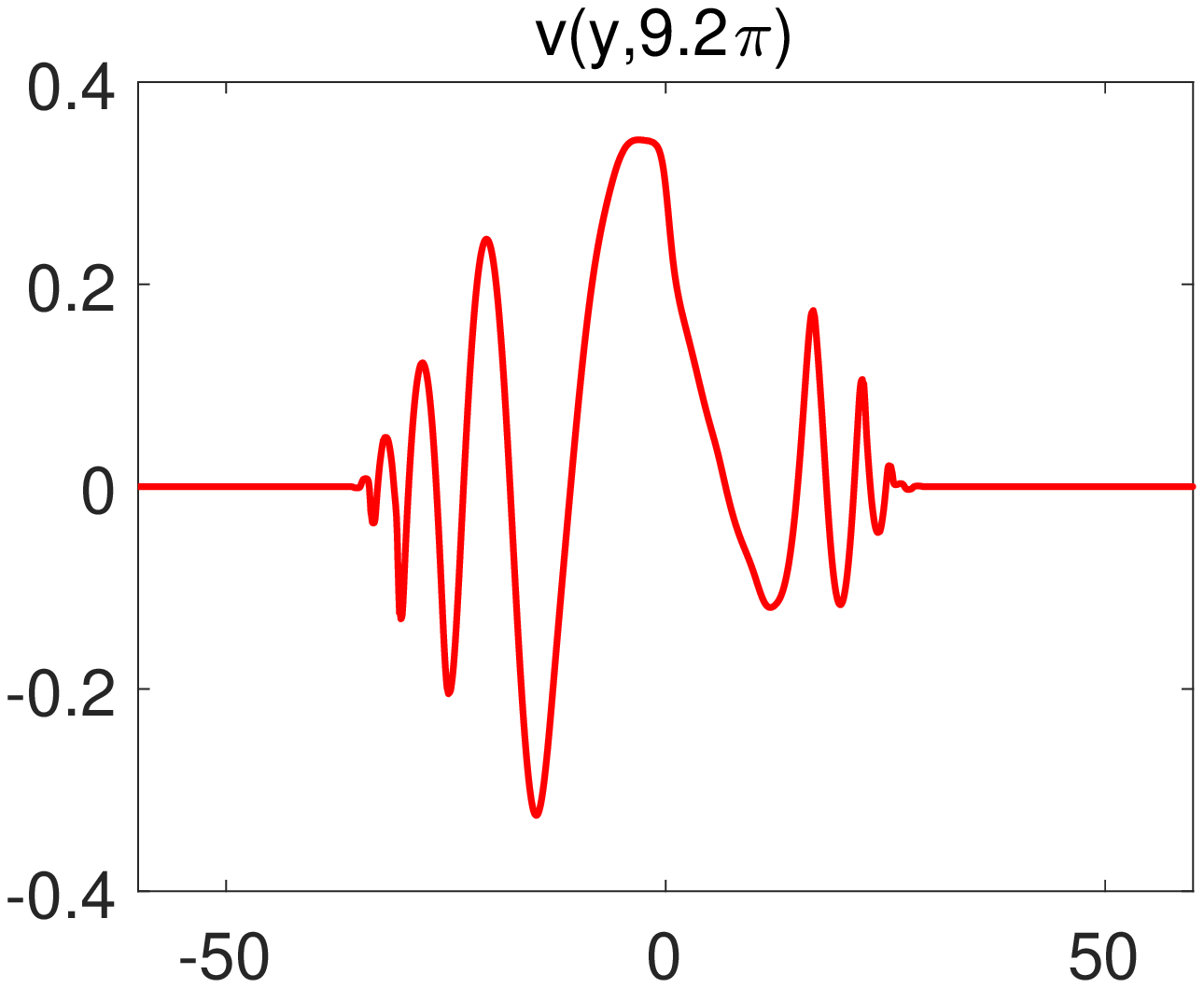}\includegraphics[width=7cm]{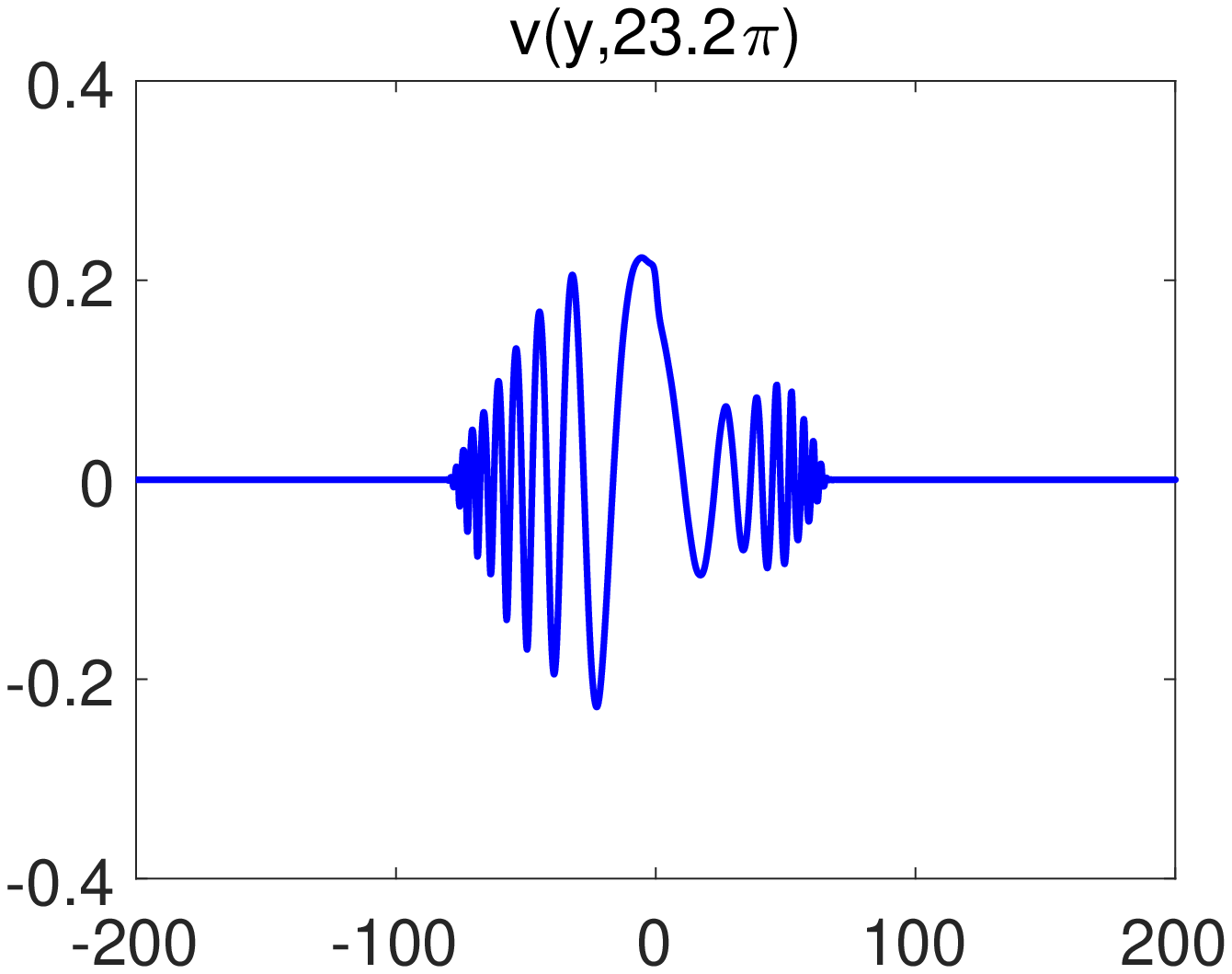}}
\centerline{\includegraphics[width=7cm]{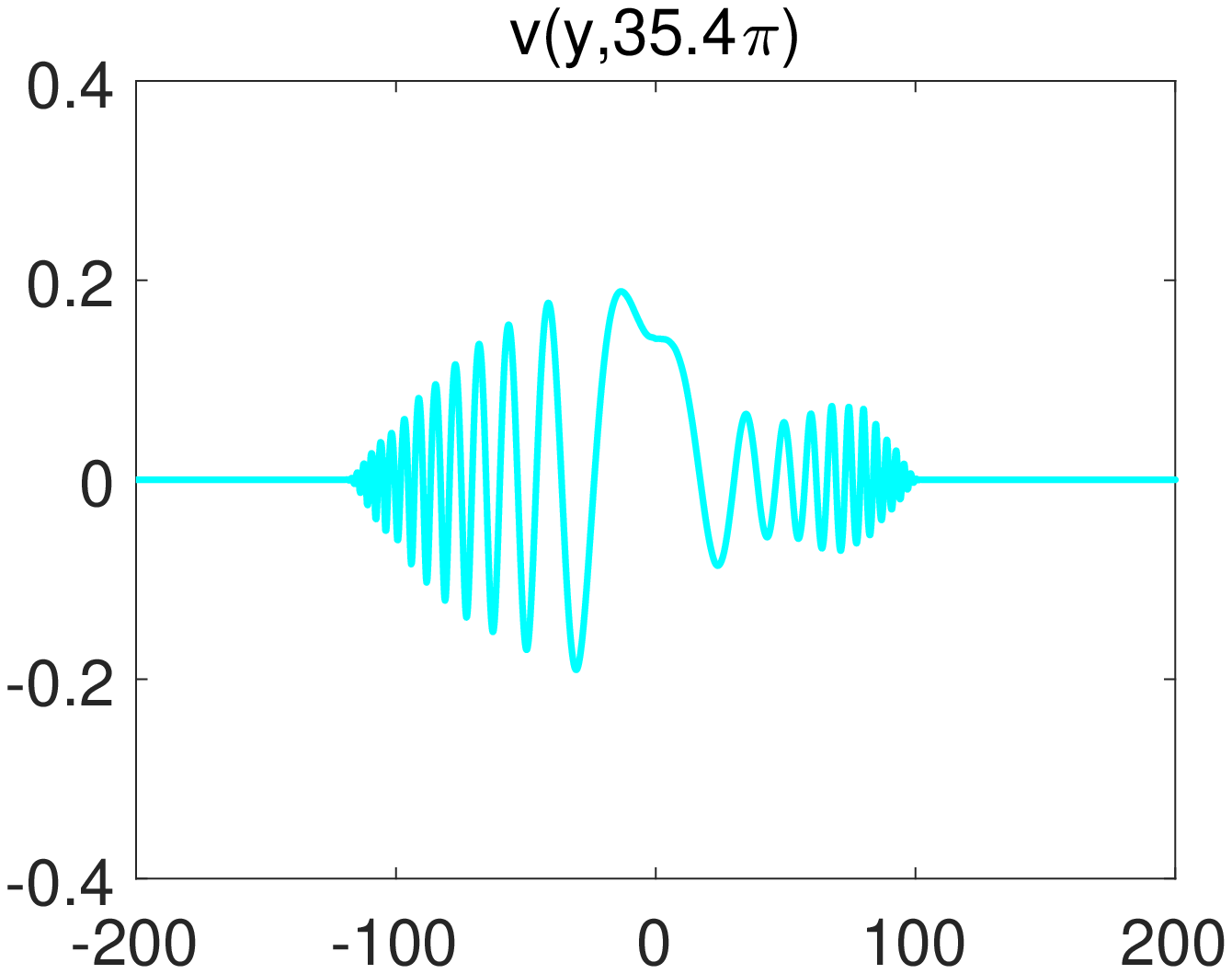}\includegraphics[width=7cm]{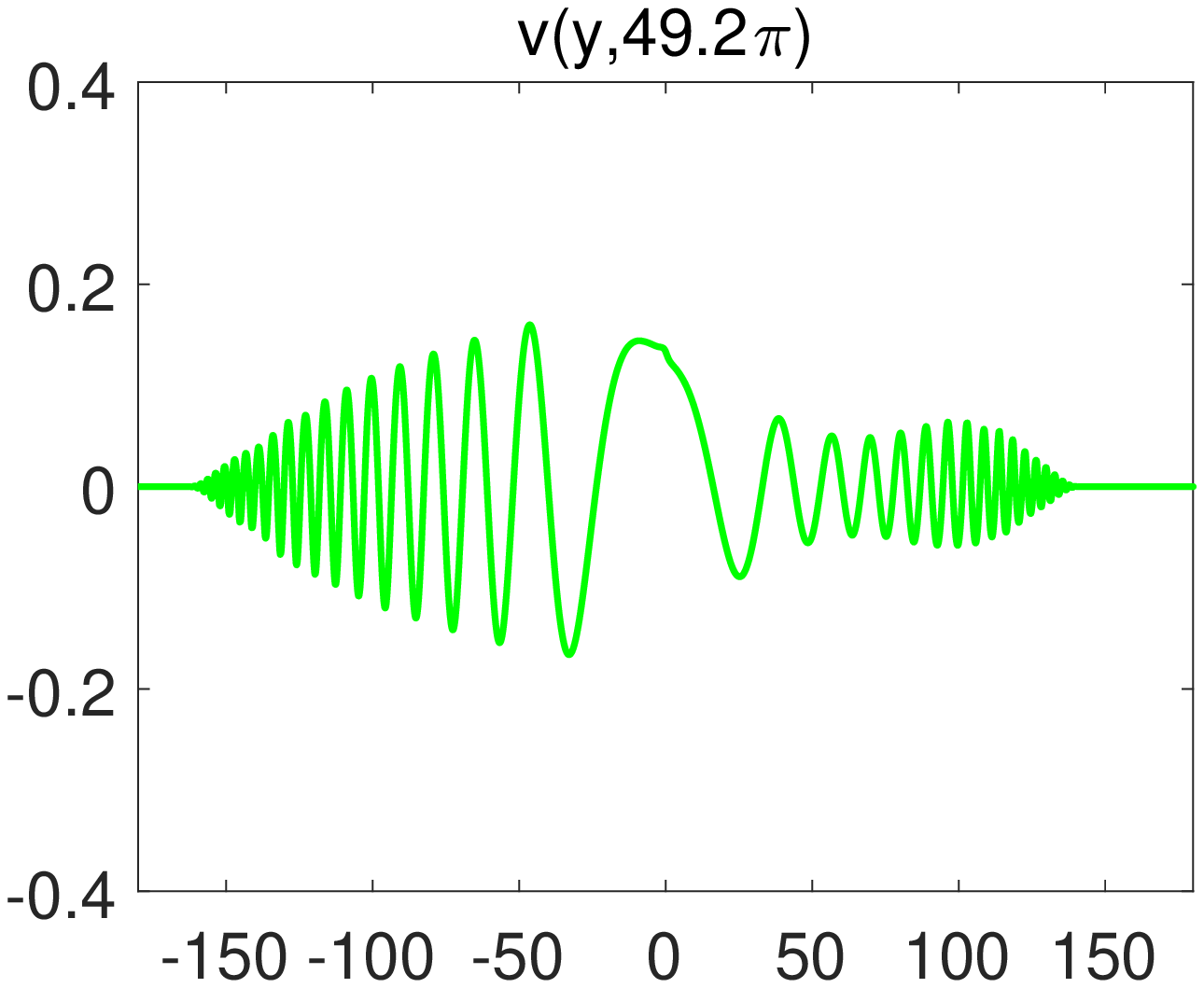}}
\caption{\sf Example 3, case (c): The snapshots of $v$ computed by the WB-CU scheme using $N=6000$ finite-volume cells during the Rossby
adjustment on the $f$-plane.\label{fig32}}
\end{figure}

\begin{figure}[ht!]
\centerline{\includegraphics[width=13.5cm]{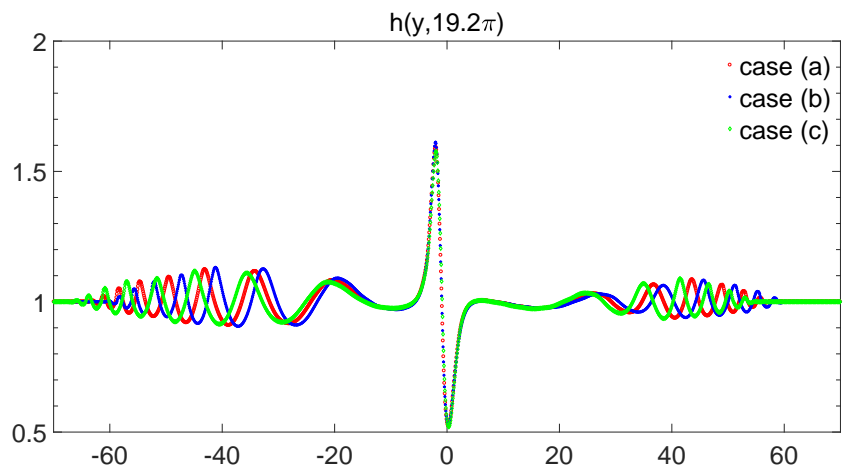}}
\caption{\sf Example 3: $h$ computed by the WB-CU scheme using $N=6000$ finite-volume cells for different $b$ at the advanced stage of the
Rossby  adjustment on the $f$-plane.\label{fig33}}
\end{figure}

In Figure \ref{fig32}, where we plot the computed meridional velocity $v$ at times $t=9.2\pi$, $23.2\pi$, $35.4\pi$ and $49.2\pi$ for case
(c), one can clearly observe the generation of inertia-gravity waves at both sides of the jet. With time increasing, the amplitude of the
signal at the center of the jet slowly decreases. The water depth $h$ computed with different initial $b$ (cases (a), (b) and (c)) is
plotted in Figure \ref{fig33}, where one can clearly see two packets of inertia-gravity waves propagating out of the jet, while $h$ rapidly
adjusts to the equilibrium profile inside the jet, as predicted by the analysis in \S\ref{sec24}.

\subsection*{Example 4 --- Breakdown of Smooth Solutions}
In this example, we numerically investigate the phenomenon of breakdown of smooth solutions described in \S\ref{sec25}. We will not be
trying to establish the relative role of different factors influencing breaking, as was discussed in \S\ref{sec25}, which would require
specially designed initial conditions (see \cite{BSZ}), and is beyond the scope of the present paper. We consider the initial conditions
consisting of a jet in a thermo-geostrophic equilibrium:
\begin{equation}
h(y,0)=1,\quad u(y,0)=3-3(\tanh(y))^2,\quad b(y,0)=10-6\tanh(y),
\label{Ex51}
\end{equation}
with a small hump in the initial meridional velocity $v$:
\begin{equation*}
v(y,0)=\left\{\begin{array}{lc}
0.1e^{-y^2}-0.1e^{-0.25},&-0.5<y<0.5,\\
0,&\mbox{otherwise}.
\end{array}
\right.
\end{equation*}
The computational domain is $[-50,50]$ and the Coriolis parameter $f\equiv1$. This setup mimics a similar simulation of a balanced jet on
the $f$-plane which was performed in \cite{BSZ}, although here the jet is ``thermal'', in a sense that $h$ is flat, and this is a gradient
of $b$ which balances the Coriolis term.

\begin{figure}[ht!]
\centerline{\includegraphics[width=7cm]{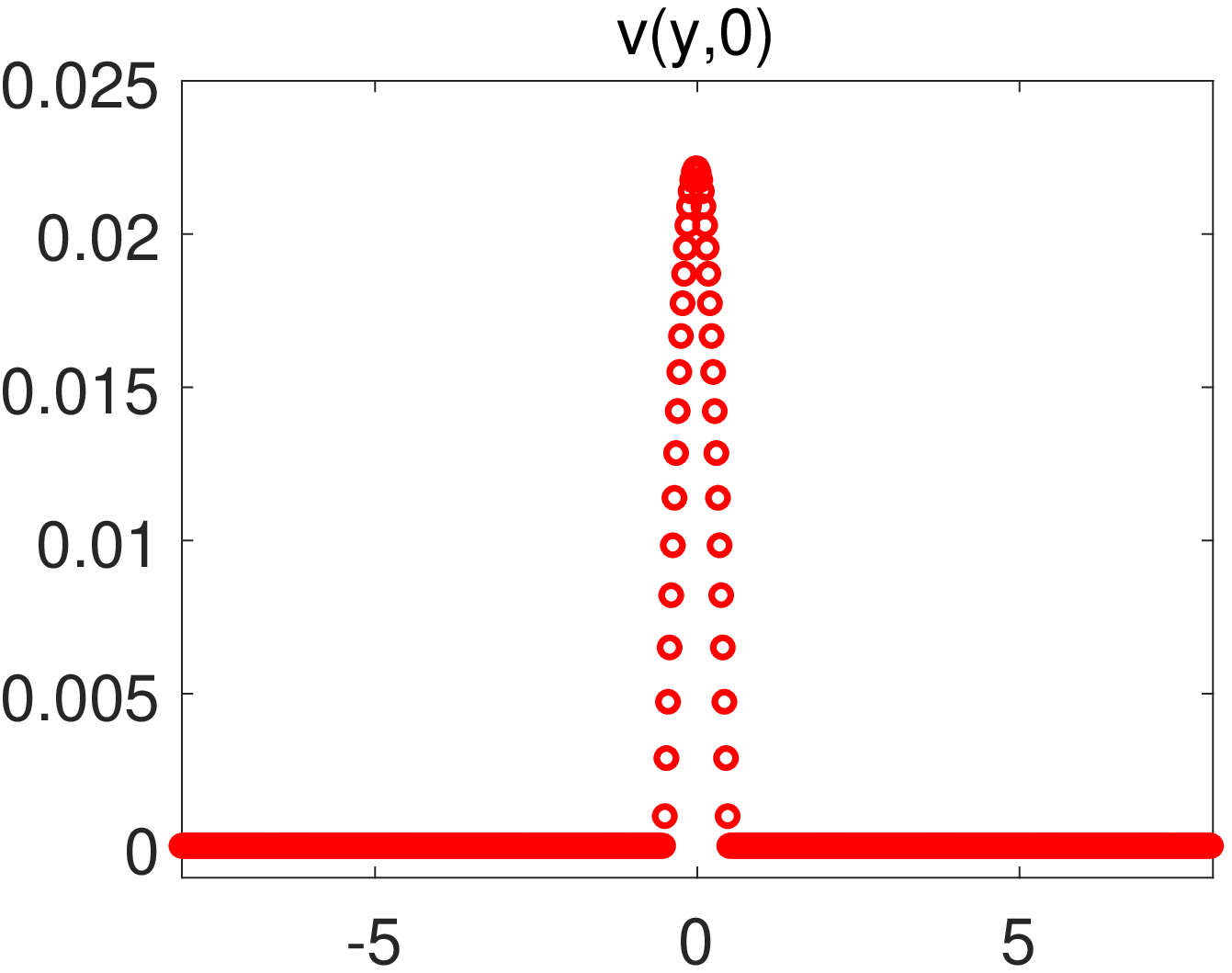}\includegraphics[width=7cm]{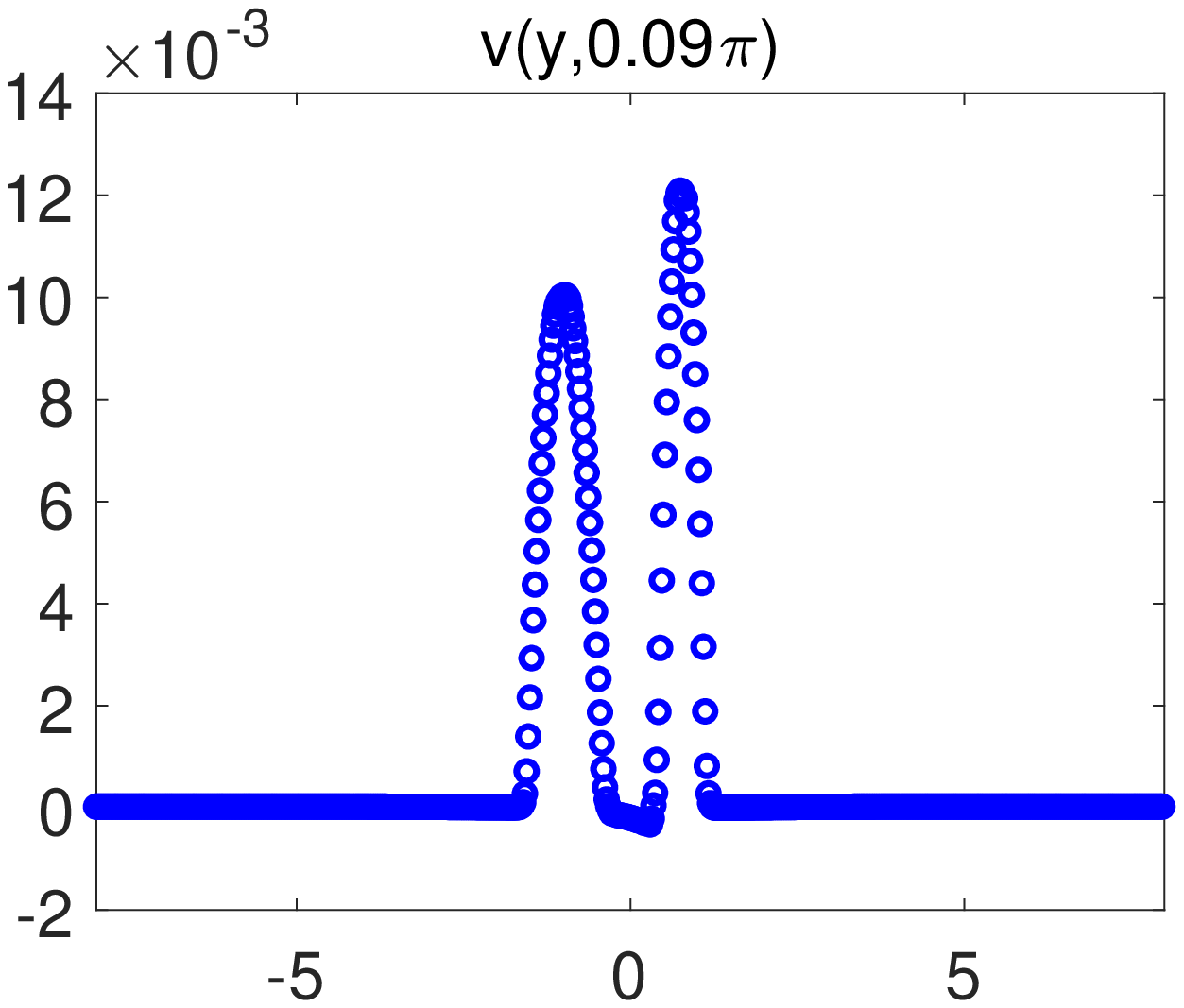}}
\centerline{\includegraphics[width=7cm]{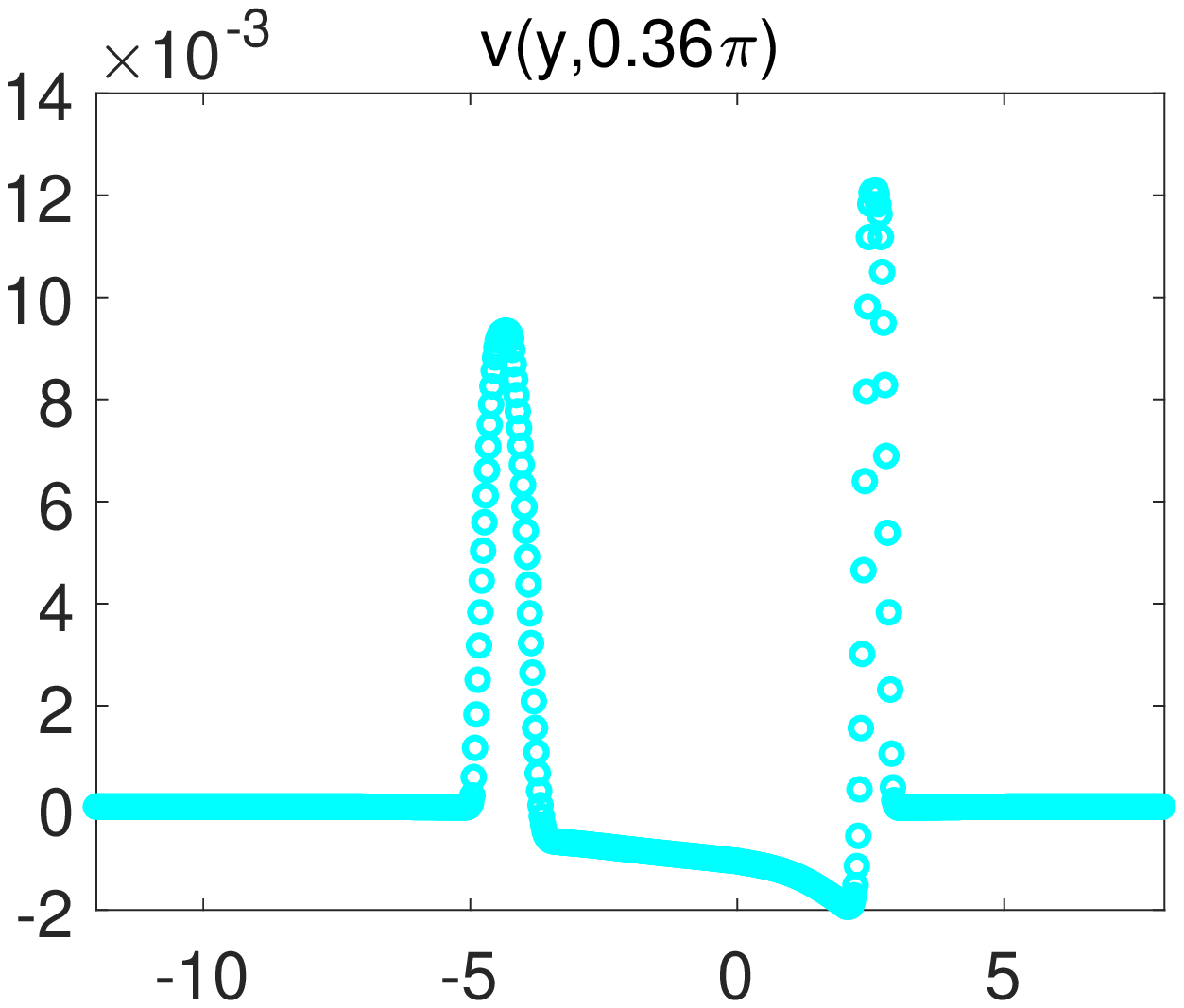}\includegraphics[width=7cm]{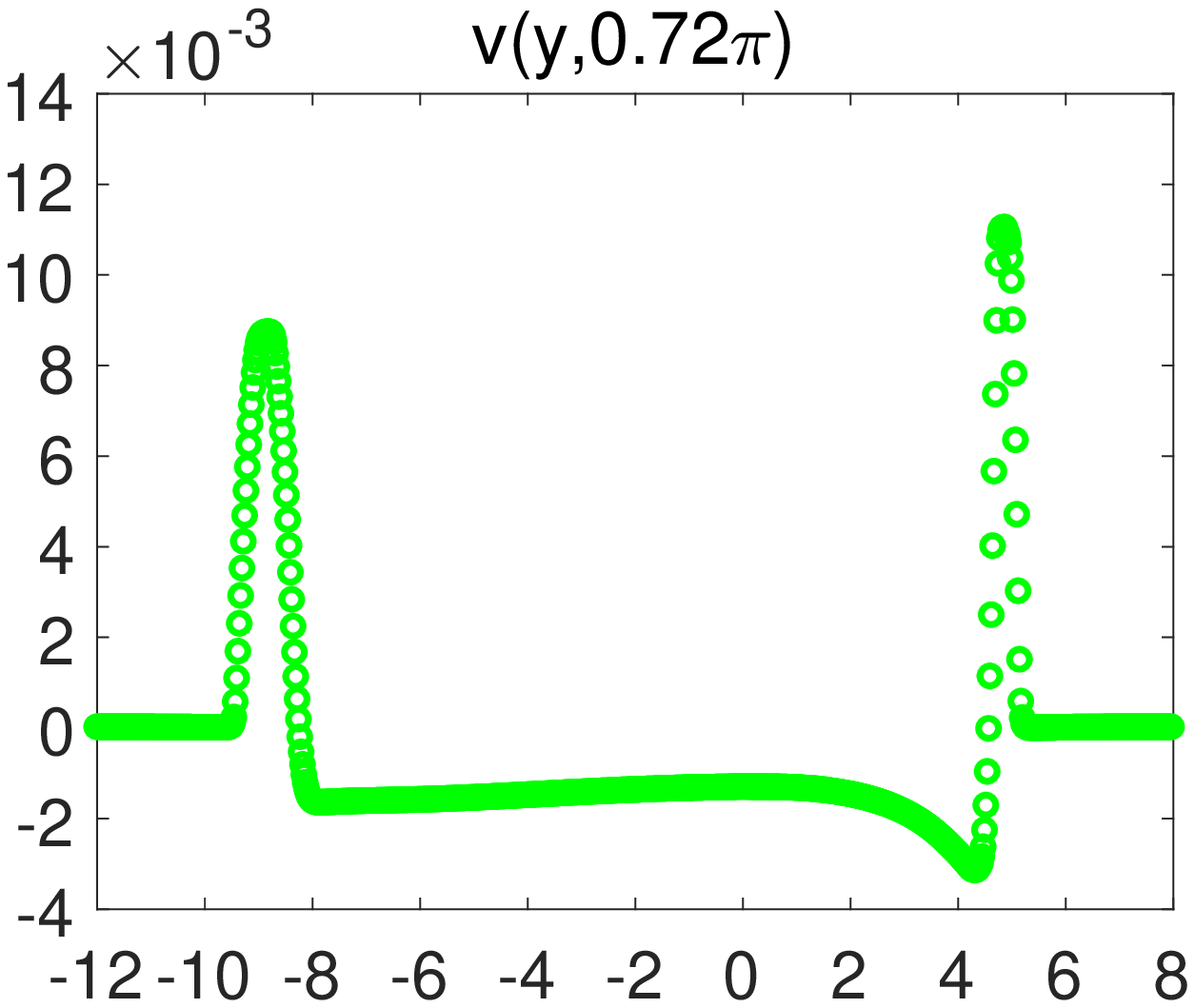}}
\centerline{\includegraphics[width=7cm]{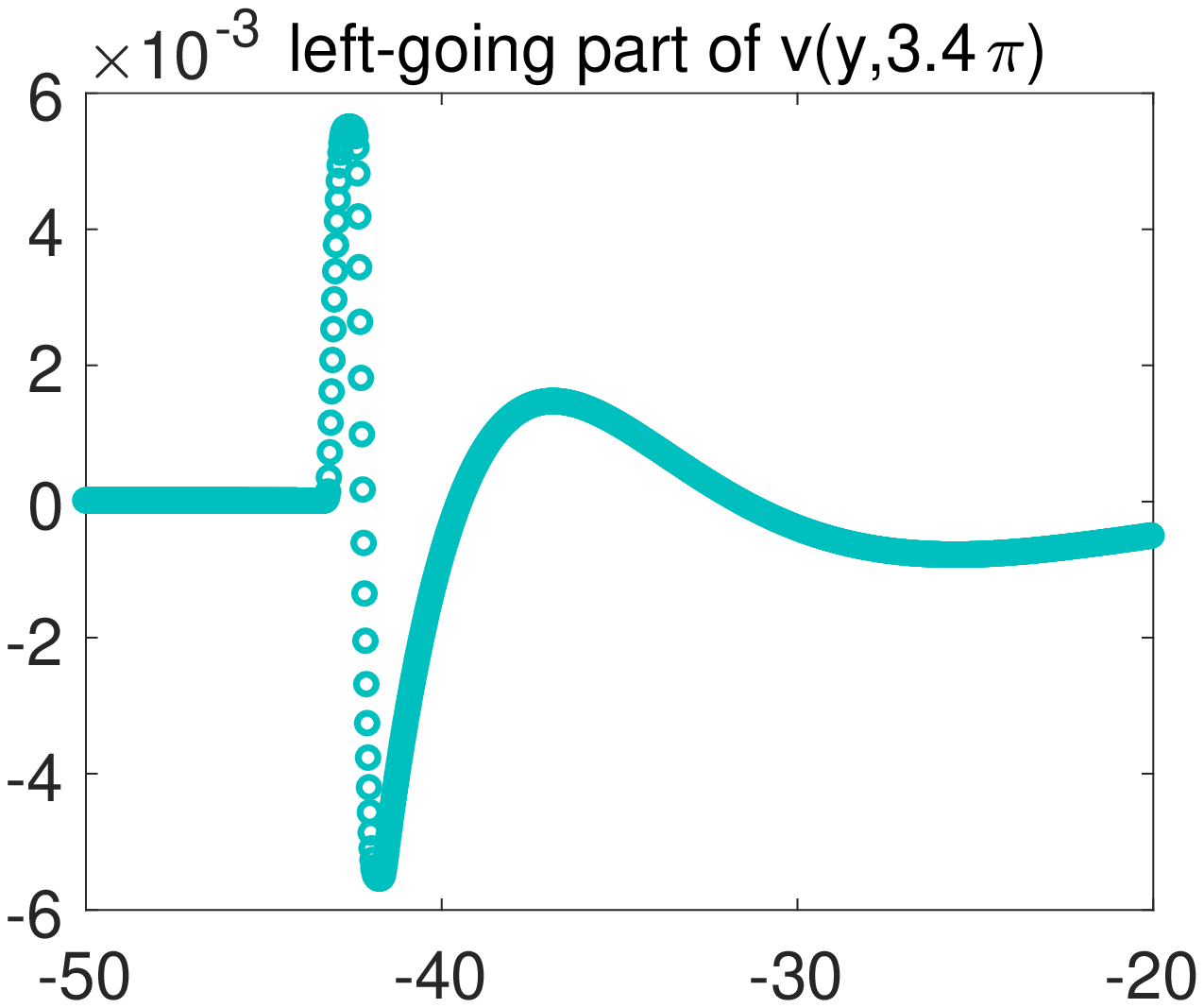}\includegraphics[width=7cm]{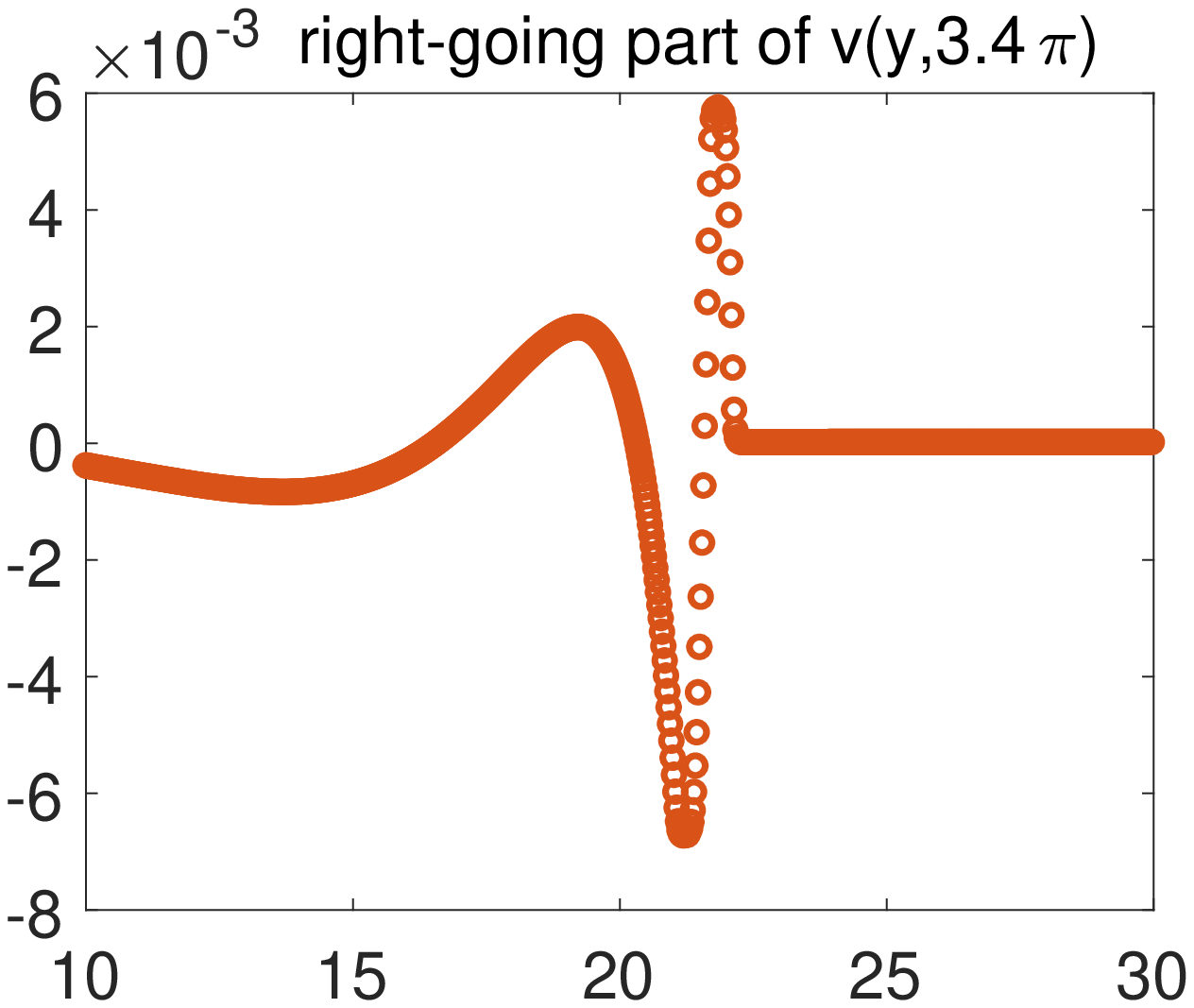}}
\caption{\sf Example 4: Evolution of the perturbation in $v$ over the balanced jet \eref{Ex51} on the $f$-plane computed by the WB-CU scheme
using $N=4000$ finite-volume cells.\label{fig51}}
\end{figure}

The evolution of the perturbation is shown in Figure \ref{fig51}. As expected, it splits in two parts propagating to the left and to the
right, both of them steepening and breaking. The right-moving part of the perturbation is breaking first, in full analogy with the
corresponding results in \cite{BSZ}, as it is moving towards the region of the decreasing background $b$ (it was a decreasing background $h$
in \cite{BSZ}). This shows that breaking over the balanced background dominated by buoyancy gradients happens in a similar way as over the
background dominated by the pressure gradients.
\subsection*{Example 5 --- Rossby Adjustment at the Equator}
\begin{figure}[ht!]
\centerline{\includegraphics[width=7cm]{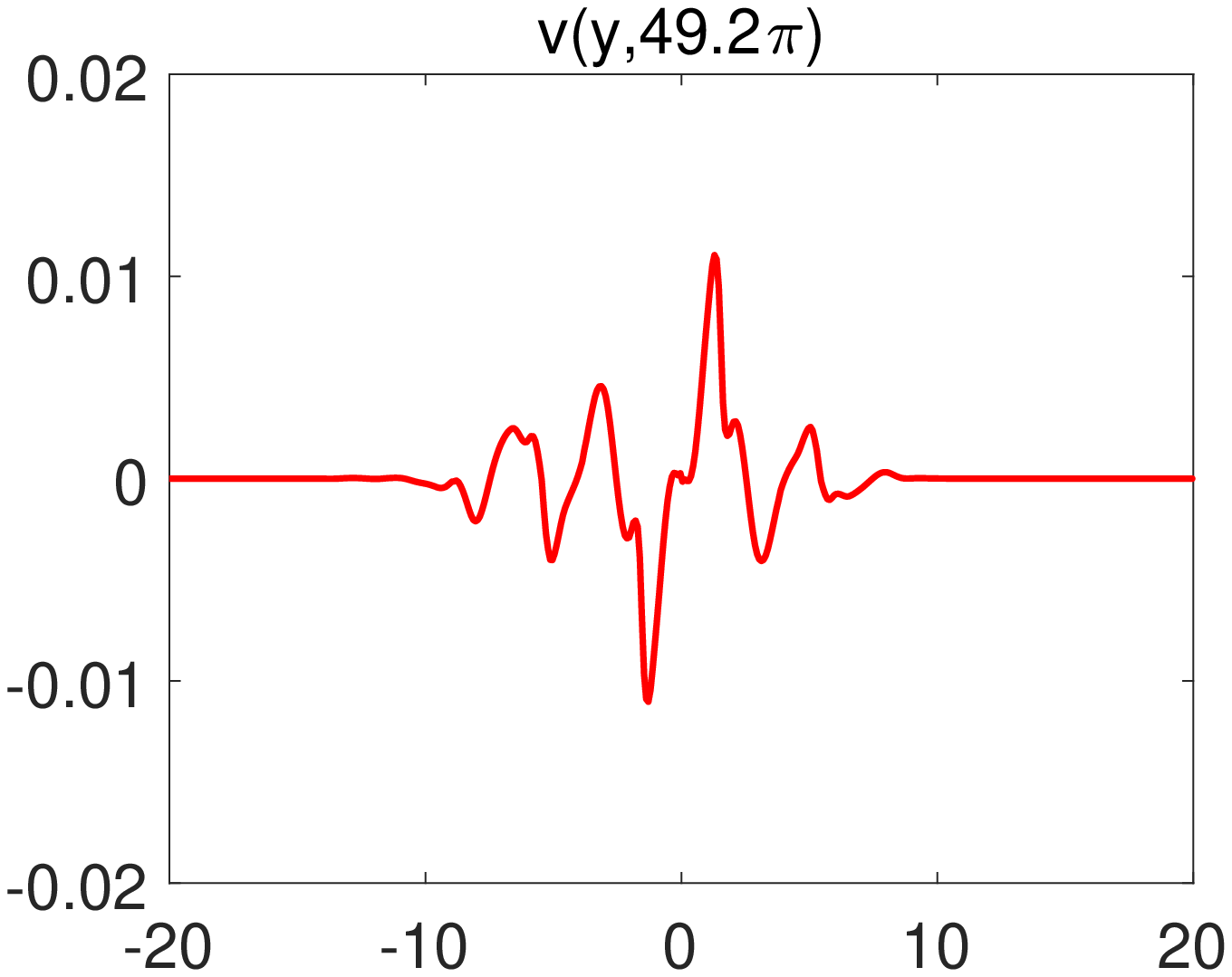}\includegraphics[width=7cm]{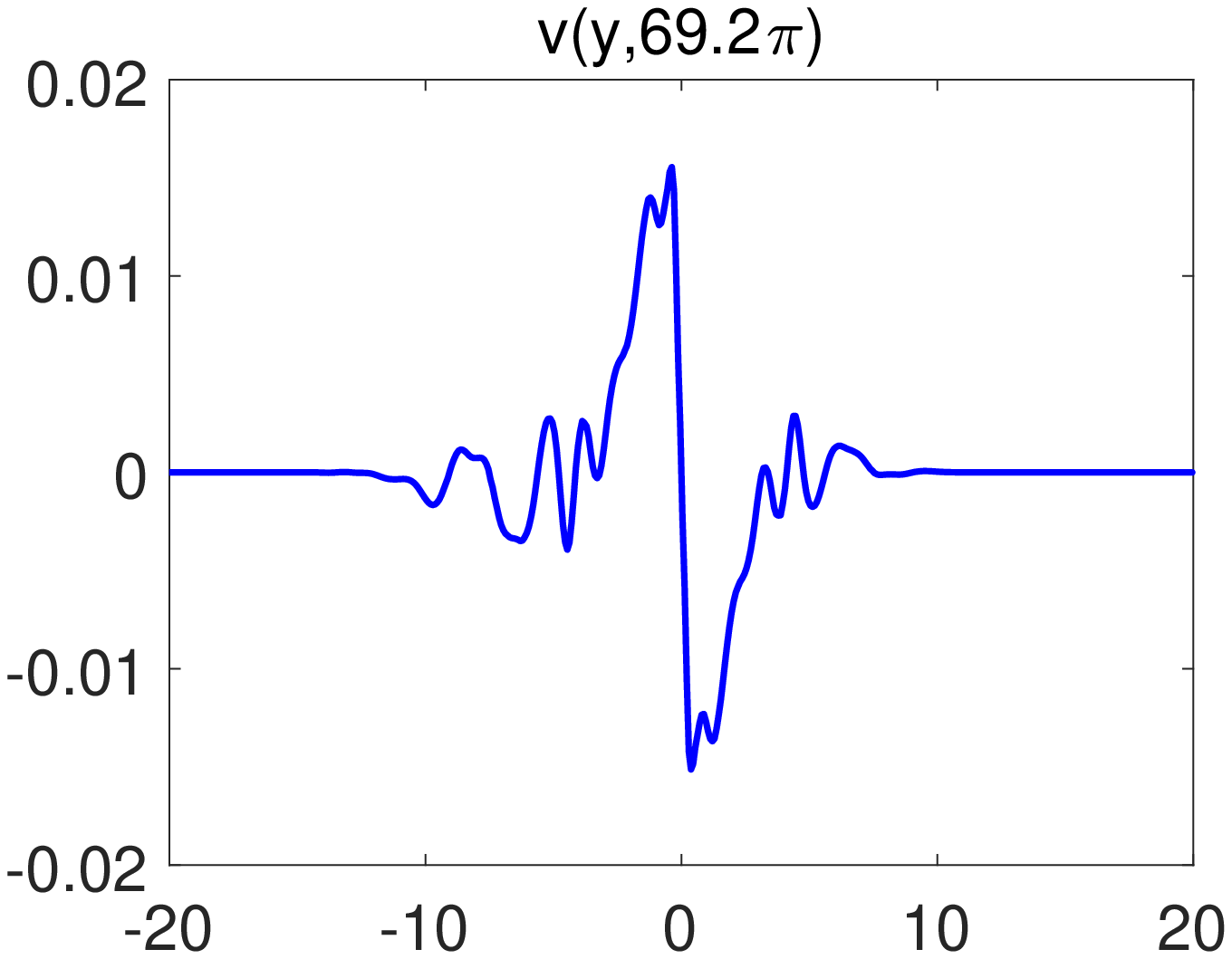}}
\centerline{\includegraphics[width=7cm]{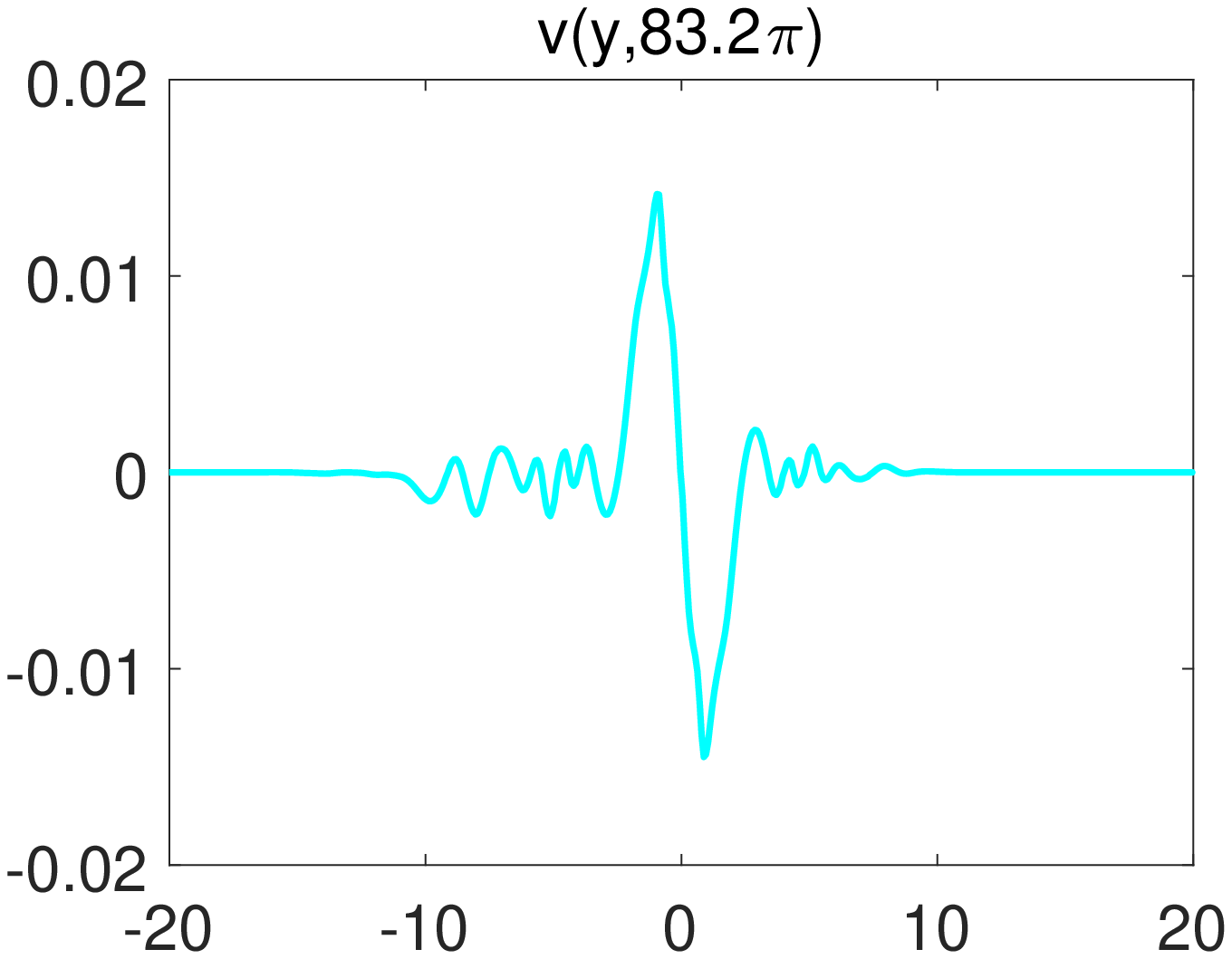}\includegraphics[width=7cm]{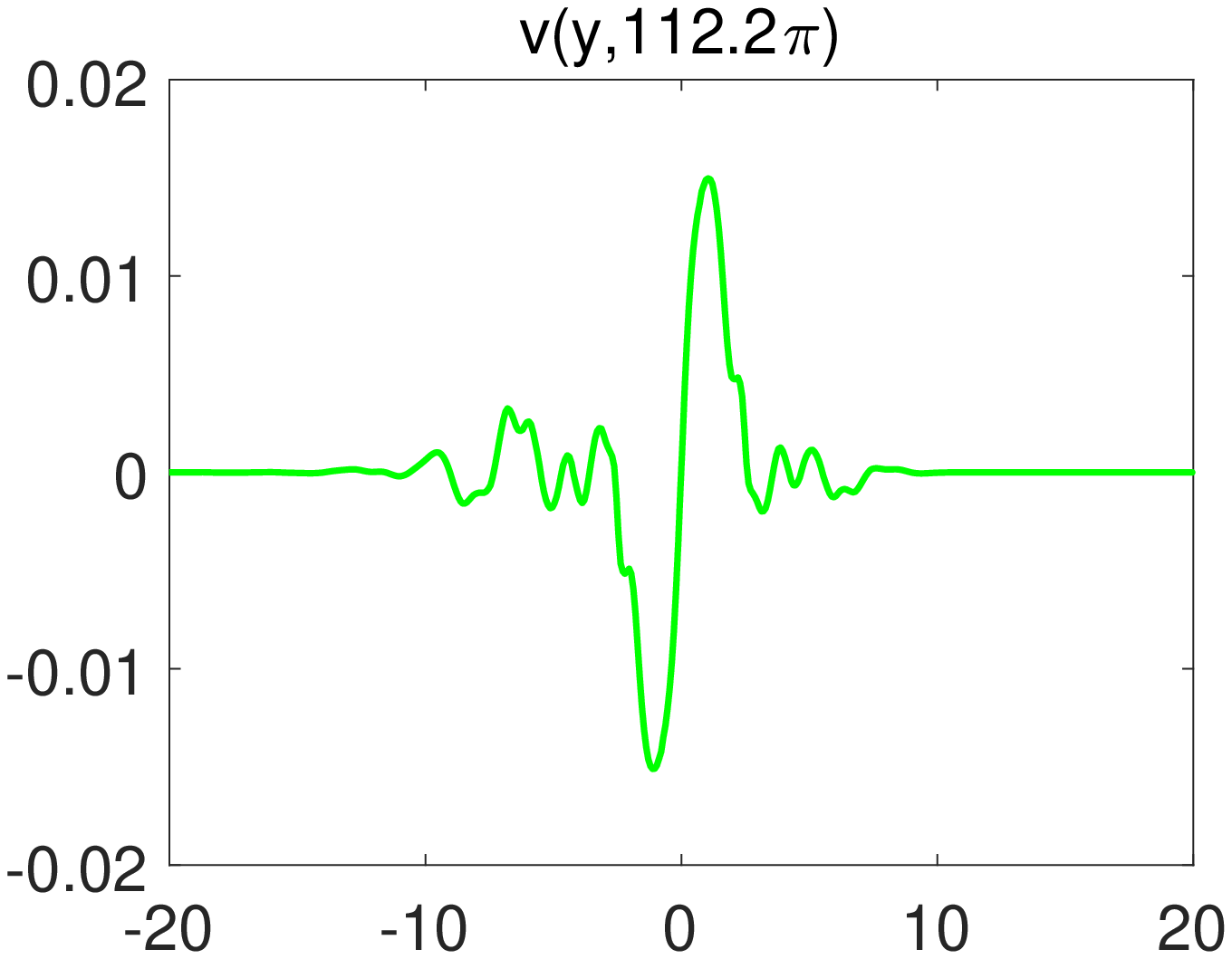}}
\caption{\sf Example 5: The snapshots of $v$ computed by the WB-CU scheme using $N=6000$ finite-volume cells with $Ro=1.0$ and $Bu=1.1$
during the equatorial adjustment.\label{fig41}}
\end{figure}

In this example, we still study the Rossby adjustment but place a jet on the equatorial beta-plane with $f(y)=0.1y$ and center it at the
Equator. As discussed in \S\ref{sec26}, the configuration with westward-oriented jet, that is, negative $u_0$, is of interest, as it may
lead to instability. We chose the following initial conditions:
\begin{equation*}
h(y,0)\equiv0.121,\quad v(y,0)\equiv0,\quad u(y,0)=-0.1e^{-y^2},\quad b(y,0)=0.1+0.01e^{-y^2},
\end{equation*}
which correspond to a small Burger number ($\mbox{Bu}=1.1$), while Rossby number is equal to 1. As in Example 3, we have also chosen a
nontrivial profile of initial $b$, as otherwise $b$ would have remained dynamically inactive. According to the results in \cite{RZT}, a jet
with these parameters is on the margin of symmetric inertial instability, and thus trapped modes of significant amplitude are expected in
this case. In Figure \ref{fig41}, we plot the meridional velocity $v$ at times $t=49.2\pi$, $69.2\pi$, $83.2\pi$ and $112.2\pi$. Trapped
modes are generated in this case in agreement with the analysis in \S\ref{sec26}. Compared to the adjustment on the $f$-plane, shown in
Figure \ref{fig32}, there are no wave-packets propagating out of the jet to the boundaries, as the waves, even in the absence of the jet,
are trapped at the Equator; see equation \eref{28abis}. We observe a large-amplitude  $v$-signal having steep gradients and localized in the
vicinity of the center of the jet. This is consistent with the marginally unstable character of the flow.

\begin{figure}[ht!]
\centerline{\includegraphics[width=8cm]{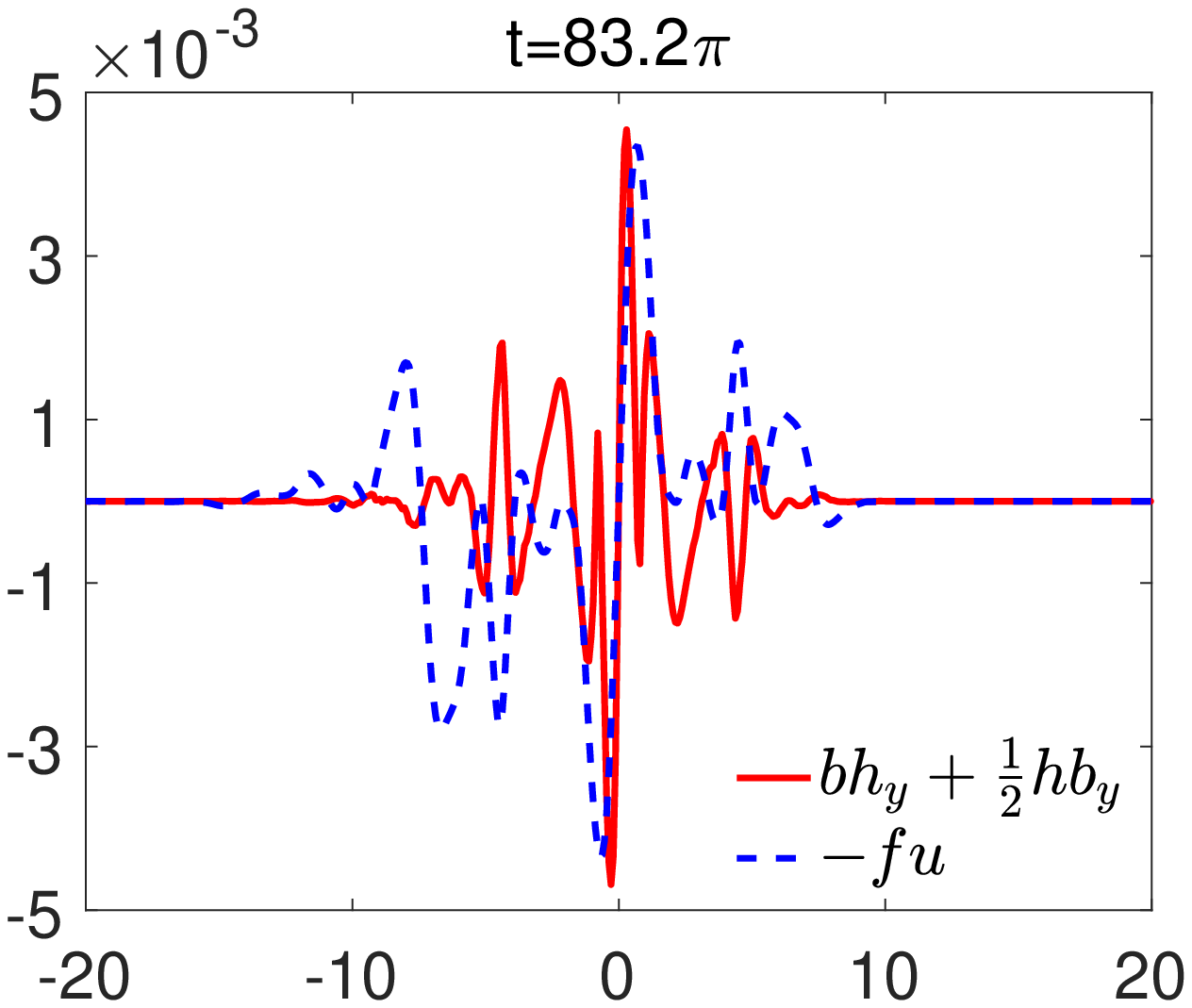}\includegraphics[width=8cm]{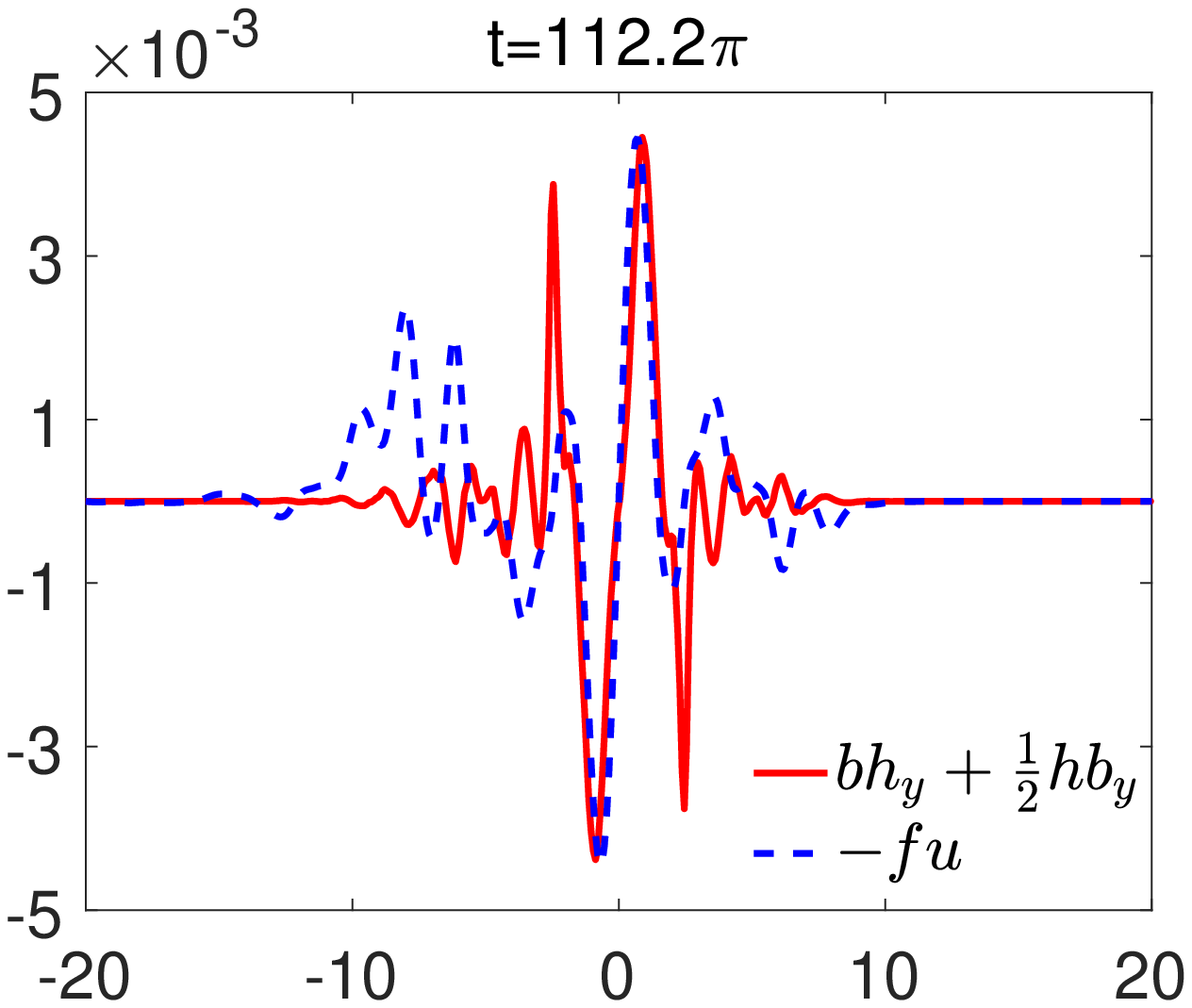}}
\caption{\sf Example 5: Evolution of $bh_y+\frac{1}{2}hb_y$ and $-fu$ computed by the WB-CU scheme using $N=6000$ finite-volume cells during
the equatorial adjustment.\label{fig402}}
\end{figure}

\begin{figure}[ht!]
\centerline{\includegraphics[width=8cm]{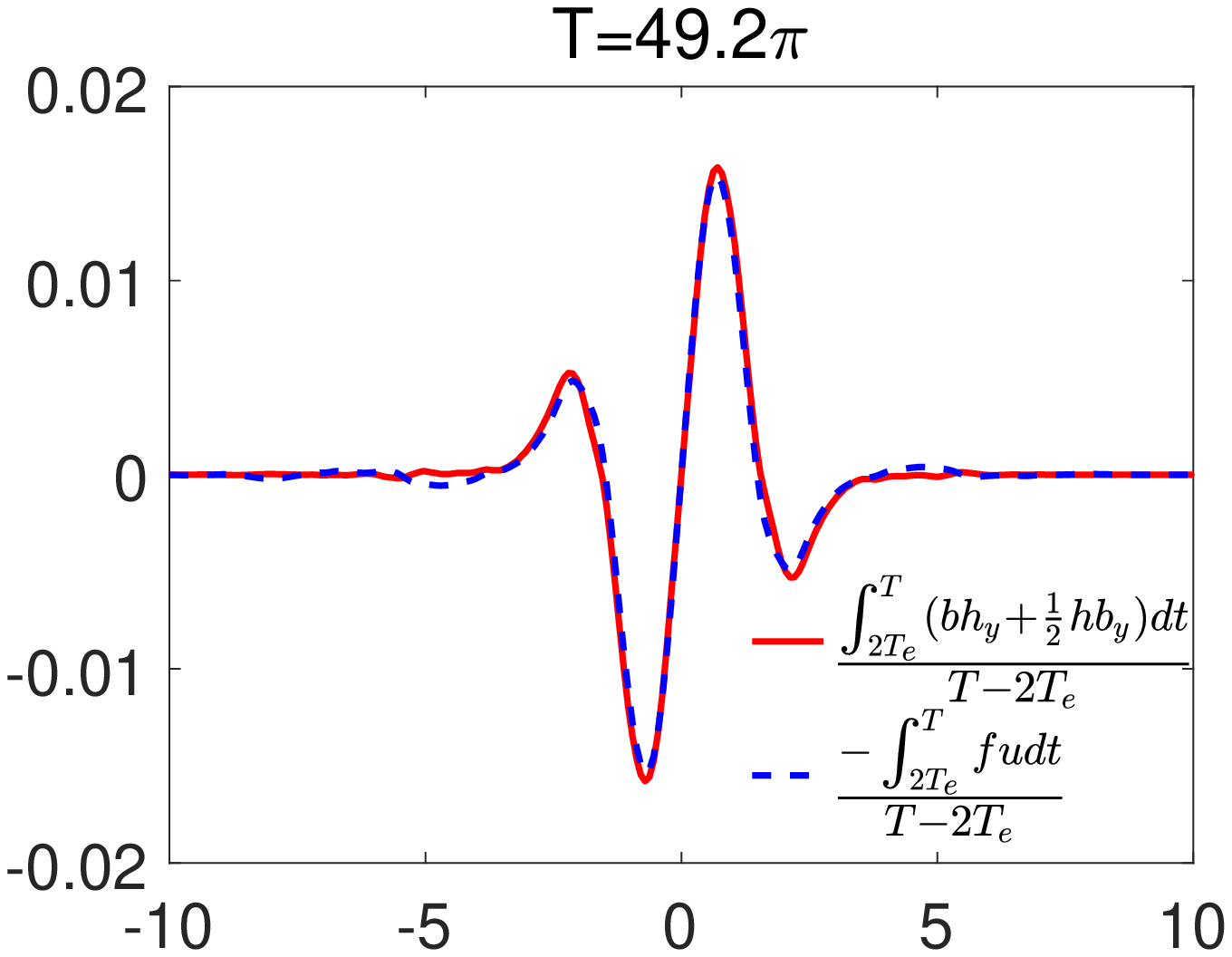}\includegraphics[width=8cm]{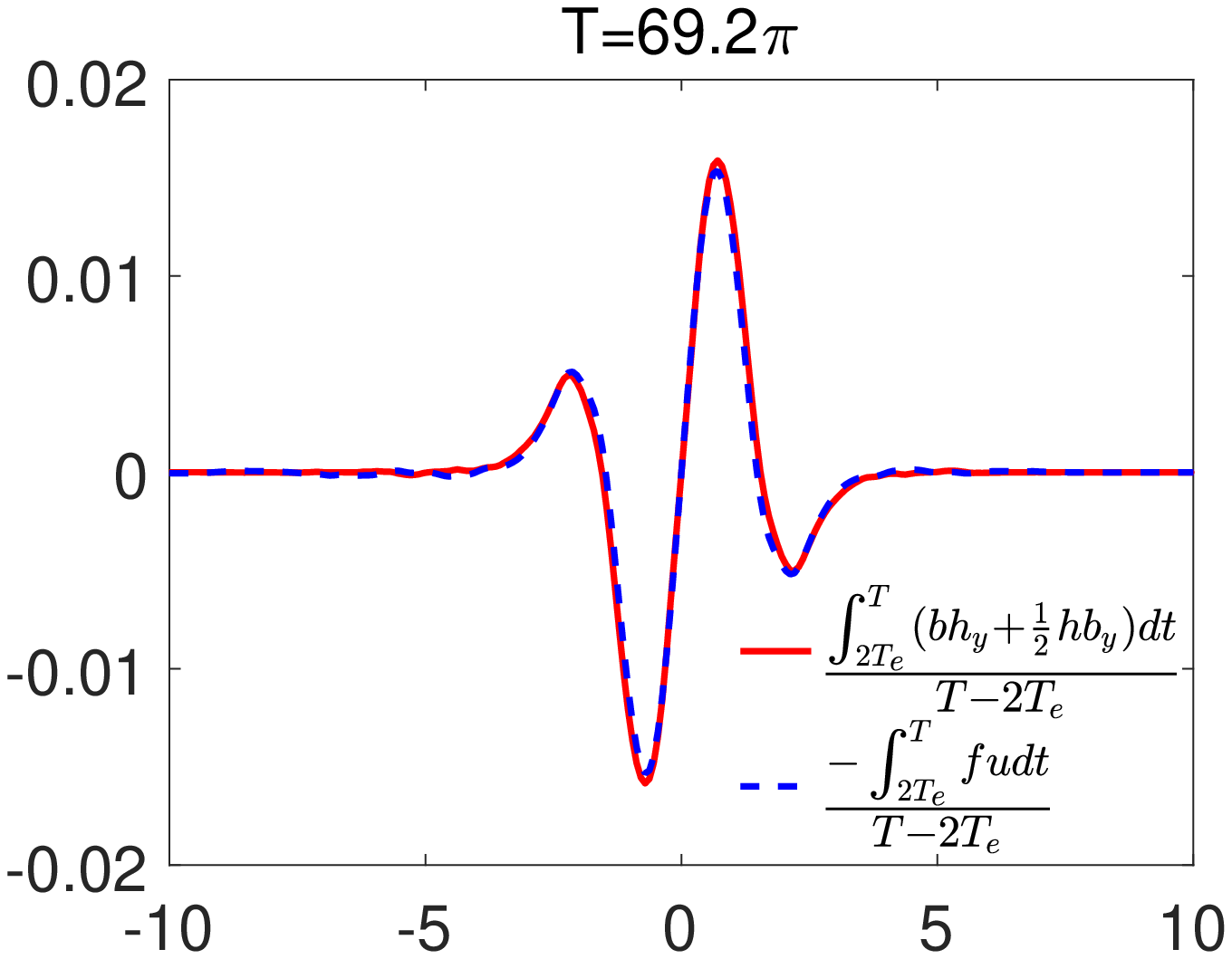}}
\caption{\sf Example 5: Evolution of $\frac{\int_{2T_e}^T\Big(bh_y+\frac{1}{2}hb_y\Big)\,{\rm d}t}{T-2T_e}$ and
$\frac{-\int_{2T_e}^Tfu\,{\rm d}t}{T-2T_e}$ computed by the WB-CU scheme using $N=6000$ finite-volume cells during the equatorial
adjustment.\label{fig402a}}
\end{figure}

The thermo-geostrophic equilibrium \eref{ge} and time averages \eref{ge2} are shown in Figures \ref{fig402} and \ref{fig402a}, respectively.
According to \S\ref{sec26}, we use the equatorial inertial period $T_e=\frac{2\pi}{\sqrt{\beta\sqrt{b_0H_0}}}$ instead of $T_f$, where
$b_0=0.1$ and $H_0=0.121$. As one can observe, even at late times, the obtained solutions do not satisfy well the relationship \eref{ge}, as
can be seen in Figure \ref{fig402}. However, from the time averages \eref{ge2} displayed in Figure \ref{fig402a}, again, like in the
$f$-plane, we conclude that even at earlier times the solutions do satisfy the relationship \eref{ge2}.

\subsection*{Example 6 --- Inertial Instability of the Balanced Equatorial Jet.}
As was shown in \S\ref{sec26}, and already discussed above, inertial instability due to the growing trapped modes arises in sufficiently
intense westward equatorial jets. In the non-thermal RSW model, this instability is expected at Rossby numbers of the order one, and Burger
number squared less than $1.5$, as was shown in \cite{BSZ}. In order to test the capability of the WB-CU scheme to capture this instability,
we performed simulations with the following initial conditions:
\begin{equation*}
h(y,0)=0.11-0.05e^{-y^2},\quad v(y,0)\equiv0,\quad u(y,0)=-0.1e^{-y^2},\quad b(y,0)\equiv0.1.
\end{equation*}
The Coriolis parameter is taken to be $f=0.1y$ and the computational domain is $[-250,250]$.

\begin{figure}[ht!]
\centerline{\includegraphics[width=7cm]{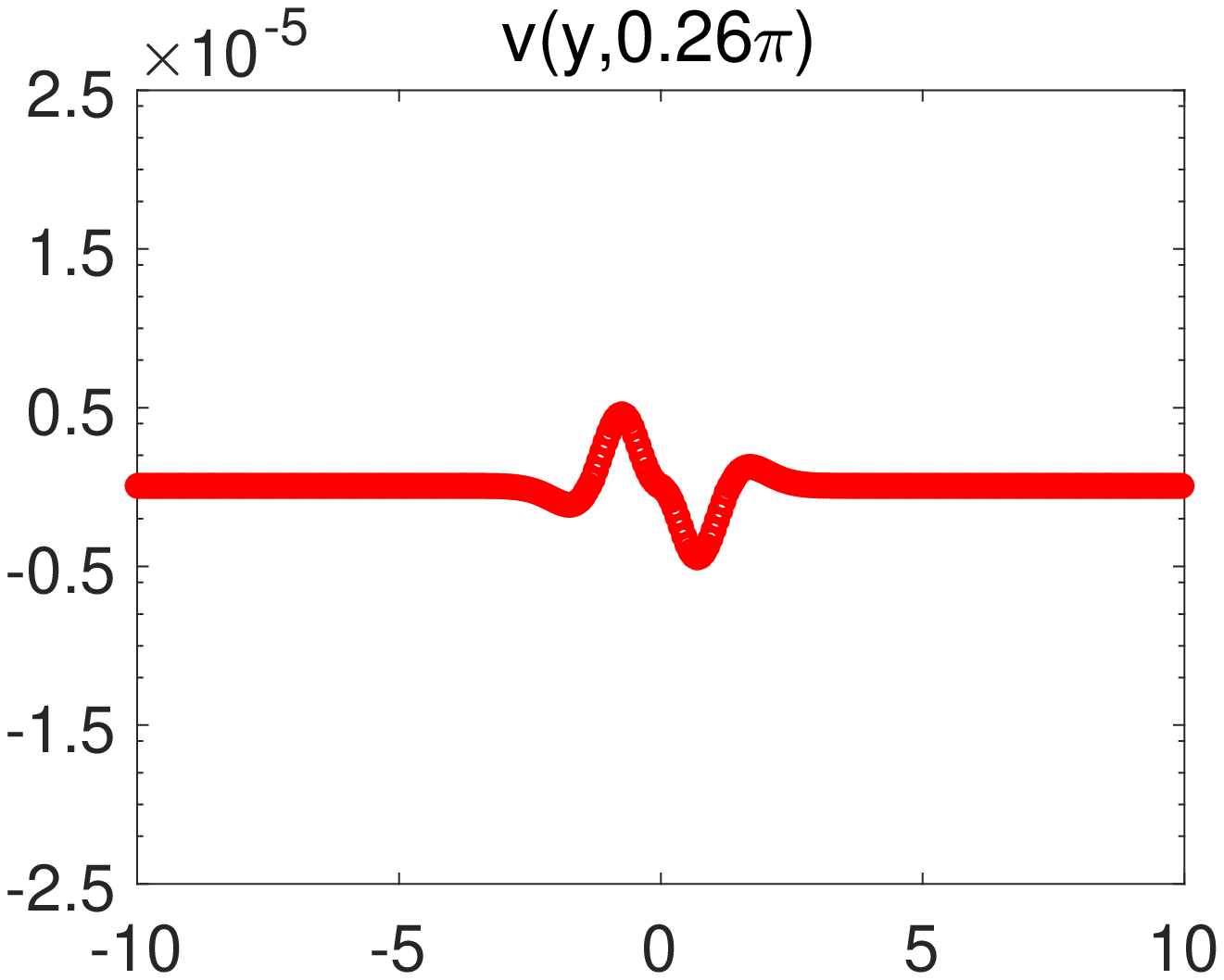}\includegraphics[width=7cm]{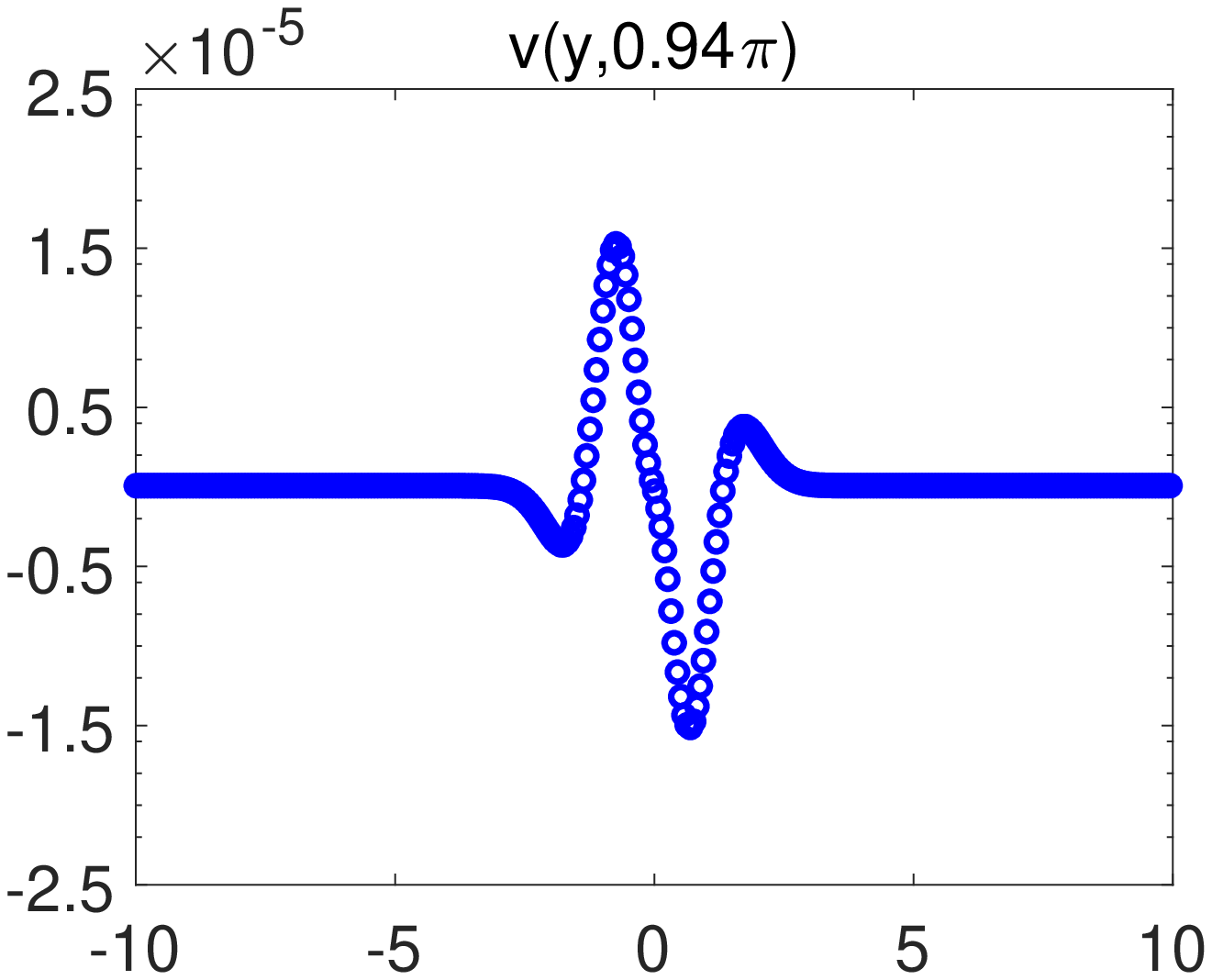}}
\centerline{\includegraphics[width=7cm]{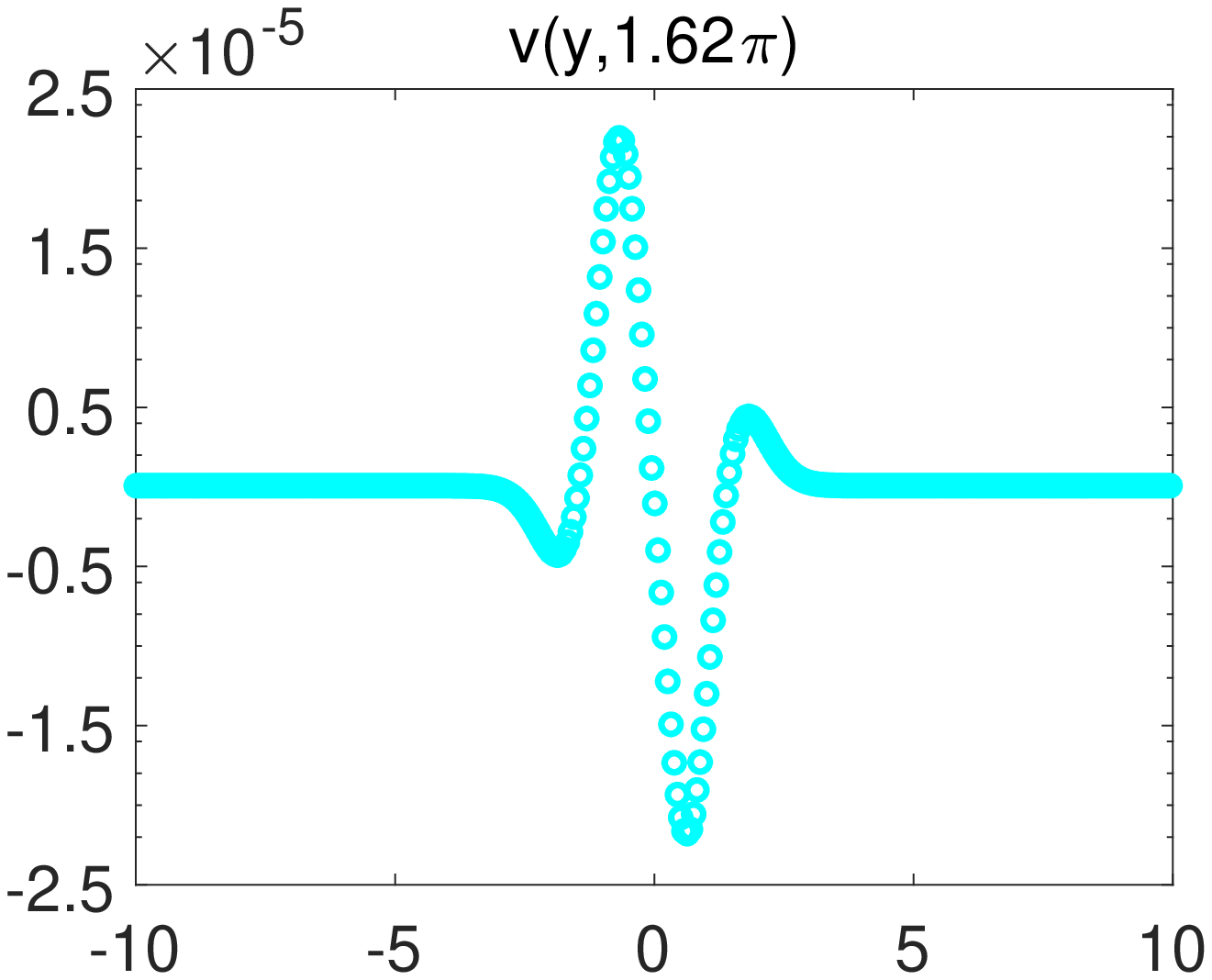}\includegraphics[width=7cm]{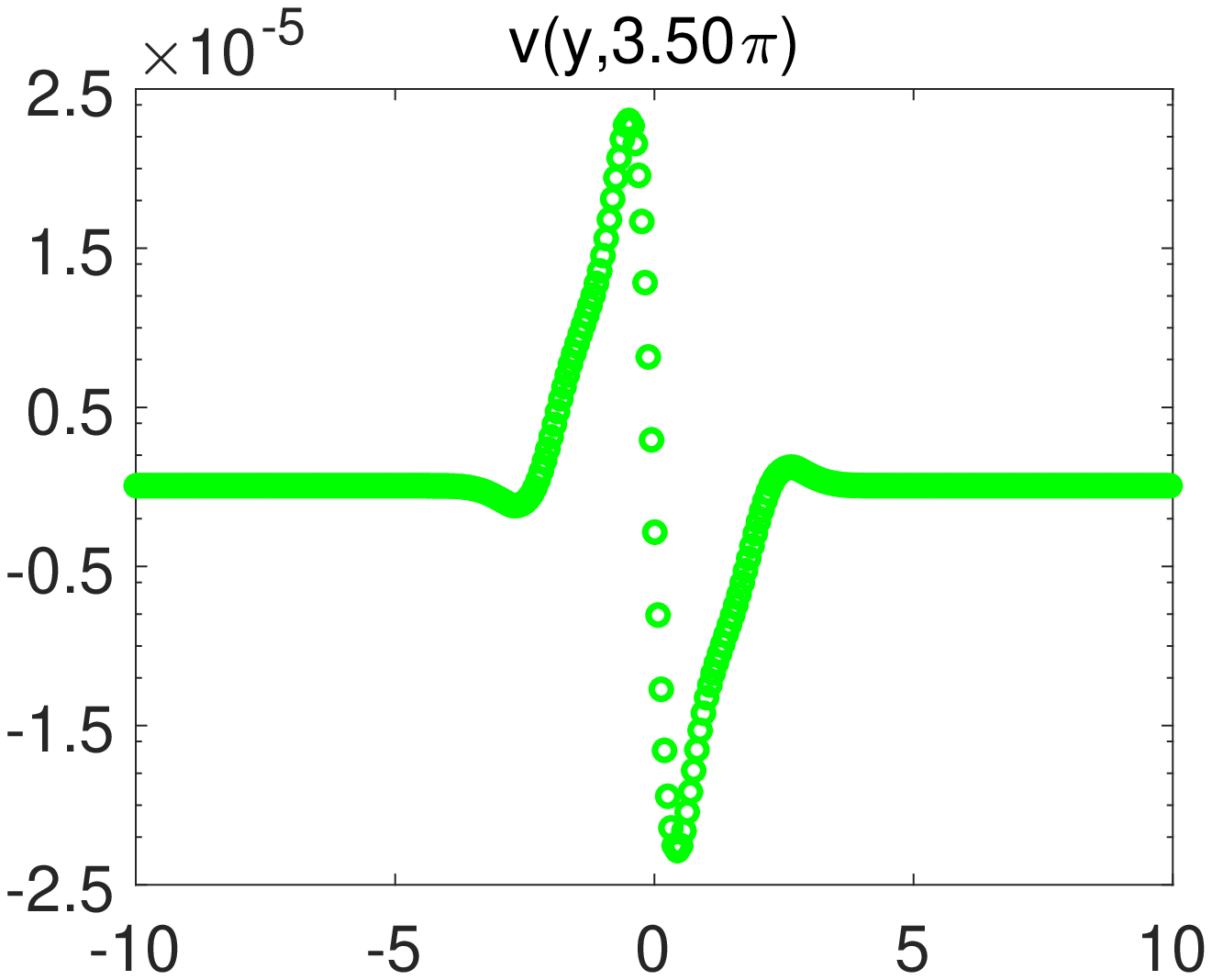}}
\caption{\sf Example 6: Evolution of the profile of $v$ computed by the WB-CU scheme using $N=8000$ finite-volume cells.\label{fig52}}
\end{figure}

This is a westward balanced equatorial jet with $\mbox{Bu}\approx1.05$ and $\mbox{Ro}=1$, which is definitely within the domain of
instability, as $b$ is flat. No perturbation is added to the jet, so the jet should remain stationary. However, the inevitable
discretization errors can be considered as a weak noise superimposed onto the ``perfect'' jet. As the spectrum of the noise is wide, a part
of it projects onto the unstable modes, which start growing. In Figure \ref{fig52}, we plot the snapshots of the meridional velocity $v$ at
times $t=0.26\pi$, $0.94\pi$, $1.62\pi$ and $3.50\pi$, which display a growing mode localized within the jet and thus confirm this scenario.

\section{Conclusion}\label{sec5}
We have constructed and tested a well-balanced central-upwind finite-volume method for the one-dimensional thermal rotating shallow water
equations. This model is a generalization of the classical rotating shallow water model, which allows one to incorporate horizontal density
and/or temperature gradients and thus make the shallow water modeling much more realistic. The results of our numerical experiments agree
with the theoretical analysis based on the Lagrangian description of the model. We have demonstrated that the proposed scheme copes well
with a nontrivial topography and effects of rotation. The thermo-geostrophic equilibria, which replace the geostrophic equilibria in the
standard rotating shallow water equations, are well maintained by the scheme. This is demonstrated on a classical example of Rossby
adjustment of a localized jet. Breakdown of smooth perturbations over the balanced configurations is demonstrated and shown to be similar to
that in the standard rotating shallow water model.

The proposed method is also tested in the equatorial region, where the Coriolis parameter does not have a constant part. This leads to
differences in the wave spectrum compared to the mid-latitude $f$-plane approximation, in which the Coriolis parameter is constant, but
also to appearance of a specific symmetric inertial instability. We have shown that the proposed scheme copes well with the adjustment of a
marginally unstable equatorial jet and captures the inertial instability of a balanced equatorial jet. Testing the scheme in the equatorial
region is of importance, as one of the main application of the scheme, once it is extended to the full plane, is to atmospheric and oceanic
dynamics in tropics.

The present work allows us to start the development of a central-upwind scheme for the two-dimensional thermal rotating shallow water
equations with confidence. This work is in progress.
\acknowledgments{The work of A. Kurganov was supported in part by NSFC grant 11771201 and NSF grants DMS-1521009 and DMS-1818666. The work
of A. Kurganov and V. Zeitlin was supported in part by the French National Program LEFE.}
\bibliographystyle{siam}
\bibliography{RTSW_ref}

\end{document}